\documentclass[12pt]{iopart}
\usepackage{amsthm,amssymb}
\usepackage{graphicx,psfrag,subfigure}
\usepackage{times}

\newif\iffiguresEPS
\newif\iffiguresPDF
\figuresPDFtrue

\newcommand{\Aset}{\mathbb{A}}

\newcommand{\Nset}{\mathbb{N}}

\newcommand{\Qset}{\mathbb{Q}}
\newcommand{\Rset}{\mathbb{R}}
\newcommand{\Sset}{\mathbb{S}}
\newcommand{\Tset}{\mathbb{T}}
\newcommand{\Zset}{\mathbb{Z}}

\theoremstyle{plain}
\newtheorem{thm}{Theorem}
\newtheorem{pro}[thm]{Proposition}
\newtheorem{lem}[thm]{Lemma}
\newtheorem{cor}[thm]{Corollary}
\newtheorem{con}{Conjecture}
\newtheorem{num}{Numerical Result}

\newtheorem{defi}{Definition}
\theoremstyle{remark}
\newtheorem{remark}{Remark}

\newcommand{\order}{\mathop{\rm o}\nolimits}
\newcommand{\acosh}{\mathop{\rm acosh}\nolimits}
\newcommand{\atan}{\mathop{\rm atan}\nolimits}

\newcommand{\Rank}{{\rm rank}}

\newcommand{\rmE}{{\rm E}}
\newcommand{\rmH}{{\rm H}}
\newcommand{\rmn}{{\rm n}}
\newcommand{\rmt}{{\rm t}}
\newcommand{\rmx}{{\rm x}}
\newcommand{\rmy}{{\rm y}}
\newcommand{\rmz}{{\rm z}}

\begin{document}

\title[The frequency map for billiards inside ellipsoids]
      {The frequency map for billiards inside ellipsoids}

\author{Pablo S. Casas and
        Rafael Ram\'{\i}rez-Ros}

\date{today}

\address{Departament de Matem\`{a}tica Aplicada I,
         Universitat Polit\`{e}cnica de Catalunya,
         Diagonal 647, 08028 Barcelona, Spain}

\begin{abstract}
The billiard motion inside an ellipsoid $Q \subset \Rset^{n+1}$
is completely integrable.
Its phase space is a symplectic manifold of dimension $2n$,
which is mostly foliated with Liouville tori of dimension $n$.
The motion on each Liouville torus becomes just a parallel
translation with some frequency $\omega$
that varies with the torus.
Besides,
any billiard trajectory inside $Q$ is tangent to $n$ caustics
$Q_{\lambda_1},\ldots,Q_{\lambda_n}$,
so the caustic parameters $\lambda=(\lambda_1,\ldots,\lambda_n)$ are
integrals of the billiard map.
The frequency map $\lambda \mapsto \omega$ is a key tool to understand
the structure of periodic billiard trajectories.
In principle,
it is well-defined only for nonsingular values of the caustic parameters.

We present two conjectures, fully supported by numerical experiments.
We obtain, from one of the conjectures, some lower bounds on the periods.
These bounds only depend on the type of the $n$ caustics.
We describe the geometric meaning, domain, and range of $\omega$.
The map $\omega$ can be continuously extended to
singular values of the caustic parameters,
although it becomes ``exponentially sharp'' at some of them.

Finally, we study triaxial ellipsoids of $\Rset^3$.
We compute numerically the bifurcation curves in the parameter space
on which the Liouville tori with a fixed frequency disappear.
We determine which ellipsoids have more periodic trajectories.
We check that the previous lower bounds on the periods are optimal,
by displaying periodic trajectories with periods four, five, and six
whose caustics have the right types.
We also give some new insights for ellipses of $\Rset^2$.
\end{abstract}




\eads{\mailto{pablo@casas.upc.es},
      \mailto{Rafael.Ramirez@upc.edu}}

\section{Introduction}

Birkhoff~\cite{Birkhoff1927} introduced the problem of
\emph{convex billiard tables} more than 80 years ago as a way to describe
the motion of a free particle inside a closed convex curve such that it
reflects at the boundary according to the law ``angle of incidence equals
angle of reflection''.
He also realized that this billiard motion can be modeled by an
area preserving map defined on an annulus.
There exists a tight relation between the invariant curves of
this billiard map and the caustics of the billiard trajectories.
Caustics are curves with the property that a trajectory,
once tangent to it, stays tangent after every reflection.
Good starting points in the literature of the billiard problem
are~\cite{KozlovTreshchev1991,Tabachnikov1995,KatokHasselblatt1995}.
We also refer to~\cite{Knill1998} for some nice figures of caustics.

When the billiard curve is an ellipse, any billiard trajectory has a caustic.
The caustics are the conics confocal to the original ellipse:
confocal ellipses, confocal hyperbolas, and the foci.
The foci are the singular elements of the family of confocal conics.
In this case, the billiard map is integrable in the sense of Liouville,
so the annulus is foliated by its invariant curves,
the billiard map becomes just a rigid rotation on its regular invariant curves,
and the rotation number varies analytically with the curve.

The billiard dynamics inside an ellipse is known.
We stress just three key results related with the search of
periodic trajectories.
First, Poncelet showed that if a billiard trajectory is periodic,
then all the trajectories sharing its caustic conic are also
periodic~\cite{Poncelet1822}. 
Second, Cayley gave some algebraic conditions to determine the
caustic conics whose trajectories are periodic~\cite{Cayley1854}.
Third, the rotation number can be expressed as a quotient of elliptic
integrals~\cite{Kolodziej1985,Tabanov1994,Waalkens_etal1997}.
We note that the search of periodic trajectories inside
an ellipse can be reduced to the search of rational rotation numbers.

A rather natural generalization of this problem is to
consider the motion of the particle inside an ellipsoid of $\Rset^{n+1}$.
Then the phase space is no longer an annulus,
but a symplectic manifold of dimension $2n$.
Many of the previous results have been extended to ellipsoids,
although those extensions are far from being trivial.
For instance,
any billiard trajectory inside an ellipsoid has $n$ caustics,
which are quadrics confocal to the original ellipsoid.
This situation is fairly exceptional, since quadrics are the only
smooth hypersurfaces of $\Rset^{n+1}$, $n\ge 2$,
that can have caustics~\cite{Berger1995}.
Then the billiard map is still completely integrable in the sense of Liouville,
being the caustics a geometric manifestation of its integrability.
In particular,
the phase space is mostly foliated with Liouville tori of dimension $n$.
The motion on each Liouville torus becomes just a parallel
translation with some frequency that varies with the torus.
Some extensions of the Poncelet theorem can be found
in~\cite{ChangFriedberg1988,Chang_etal1993a,Chang_etal1993b,Previato1999}.
Several generalized Cayley-like conditions were stated
in~\cite{DragovicRadnovic1998a,DragovicRadnovic1998b,DragovicRadnovic2006,DragovicRadnovic2008}.
Finally,
the frequency map can be expressed in terms of hyperelliptic integrals,
see~\cite{Chang_etal1993b,PopovTopalov2009}.
The setup of these last two works is $\Rset^3$,
but their formulae are effortless extended to $\Rset^{n+1}$.

From Jacobi and Darboux it is known that hyperelliptic functions play a
role in the description of 
the billiard motion inside ellipsoids and the geodesic flow on ellipsoids.
Nevertheless, we skip the algebro-geometric approach
(the interested reader is referred
 to~\cite{Moser1980,Knorrer1980,Knorrer1985,Audin1994,Audin1995})
along this paper, in order to emphasize the dynamical point of view.
Here, we just mention that the billiard dynamics inside an ellipsoid can be
expressed in terms of some Riemann theta-functions associated
to a hyperelliptic curve,
and so, one can write down explicitly the parameterizations of the invariant
tori that foliate the phase space; see~\cite{Veselov1988,Fedorov1999}.

Periodic orbits are the most distinctive special class of orbits.
Therefore, the first task to carry out in any dynamical system
should be their study,
and one of the simplest questions about them is to look for minimal periods.
In the framework of smooth convex billiards the minimal period is always two.
Nevertheless,
since all the two-periodic billiard trajectories inside
ellipsoids are \emph{singular}
---in the sense that some of their caustics are singular elements
of the family of confocal quadrics---, two questions arise.
Which is the minimal period among \emph{nonsingular}
billiard trajectories?
Which ellipsoids display such trajectories?

In order to get a flavor of the kind of results obtained in this paper,
let us consider the three-dimensional problem.
Let $Q$ be the triaxial ellipsoid given by $x^2/a + y^2/b + z^2/c = 1$,
with $0 < c < b < a$.
We assume that $a=1$ without loss of generality.
Any billiard trajectory inside $Q$ has as caustics two elements
of the family  of confocal quadrics given by
\[
Q_\lambda =
\left\{
(x,y,z) \in \Rset^3 :
\frac{x^2}{a-\lambda} + \frac{y^2}{b-\lambda} + \frac{z^2}{c-\lambda} = 1
\right\}.
\]
We restrict our attention to nonsingular trajectories.
That is, trajectories whose caustics are ellipsoids: $0 < \lambda < c$,
1-sheet hyperboloids: $c < \lambda < b$,
or 2-sheet hyperboloids: $b < \lambda < a$.
The singular values $\lambda \in \{a,b,c\}$ are discarded.
It is known that there are only four types of couples of
nonsingular caustics: EH1, H1H1, EH2, and H1H2.
The notation is self-explanatory.
It is also known that any nonsingular periodic billiard trajectory inside $Q$
has three so-called \emph{winding numbers} $m_0,m_1,m_2 \in \Nset$
which describe how the trajectory folds in $\Rset^3$.
For instance, $m_0$ is the period.
The geometric meanings of $m_1$ and $m_2$ depend on the type
of the couple of caustics, see section~\ref{sec:3D}.
We stand out two key observations about winding numbers.
First, some of them must be even.
Namely, the ones that can be interpreted
as the number of crossings with some coordinate plane. 
Second, we conjecture that they are ordered as follows: $m_2 < m_1 < m_0$.
This unexpected behaviour is supported by extensive numerical experiments.
In fact, we believe that it holds in any dimension.
The combination of both observations crystallizes in the following lower bounds.

\begin{thm}\label{thm:LowerBounds3D}
If the previous conjecture on the winding numbers holds,
any periodic billiard trajectory inside a triaxial ellipsoid
of $\Rset^3$ whose caustics are of type EH1, H1H1, EH2, and H1H2
has period at least five, four, five, and six, respectively.
\end{thm}

All the billiard trajectories of periods two and three are singular.
The two-periodic ones are contained in some coordinate axis,
so they have two singular caustics.
The three-periodic ones are contained in some coordinate plane,
so they have one singular caustic.

\begin{figure}
\iffiguresPDF
\begin{center}
\includegraphics[width=6in]{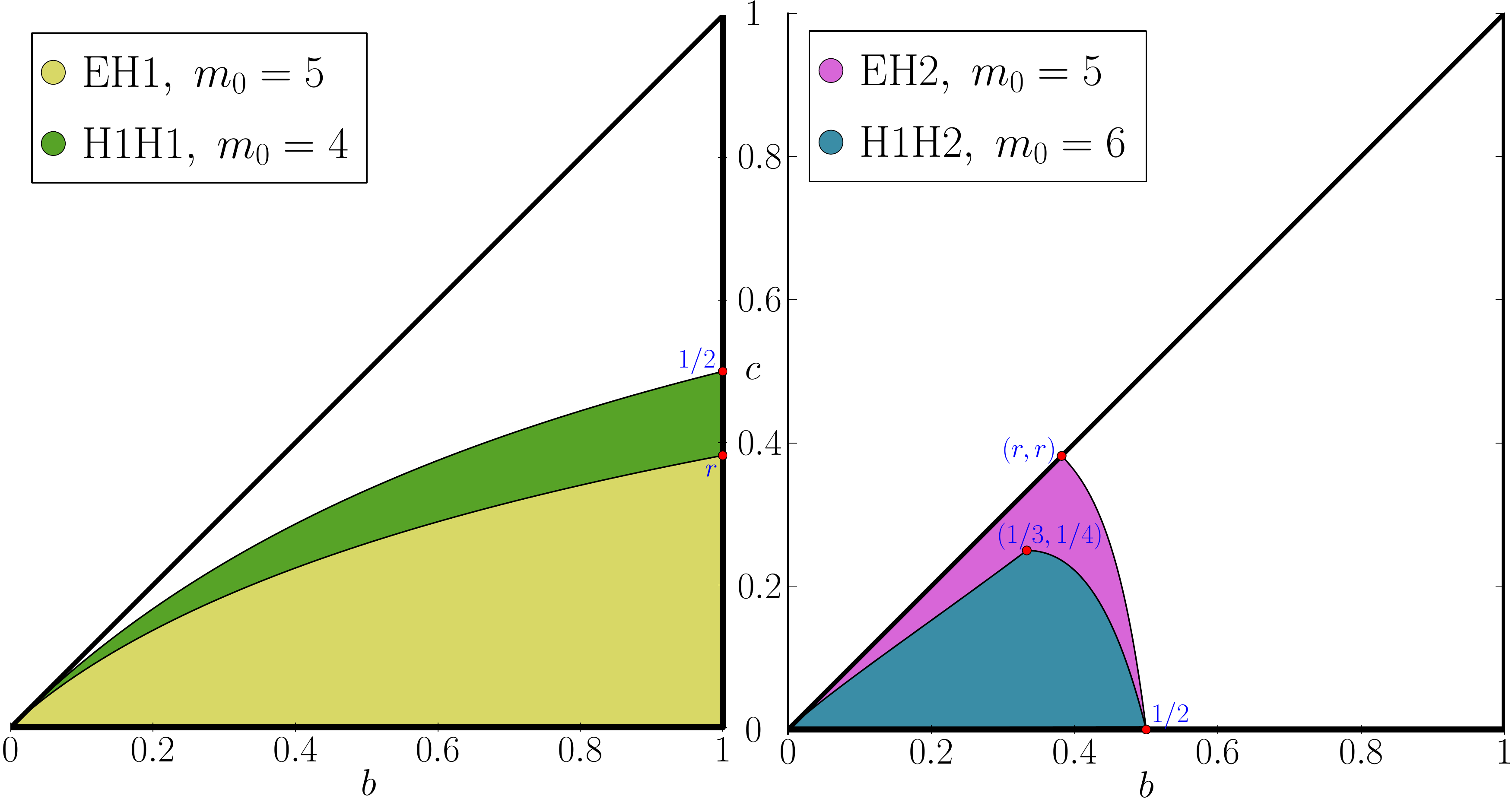}
\end{center}
\else
\vspace{3.275in}
\fi
\caption{The four regions of ellipsoids with minimal periodic trajectories.
Left: The yellow region (type EH1, period 5) is contained
      in the green one (type H1H1, period 4).
Right: The blue region (type H1H2, period 6) is contained
       in the magenta one (type EH2, period 5).}
\label{fig:BifurcationsMinimal}
\end{figure}

We shall prove in section~\ref{sec:FrequencyMap}
the generalization of these lower bounds to any dimension,
see theorem~\ref{thm:LowerBounds}.
Samples of periodic trajectories with minimal periods are
shown in figure~\ref{fig:PeriodicTrajectories3D}.
Hence, these lower bounds are optimal.
Next, we look for ellipsoids with minimal periodic trajectories.
We recall that $a = 1$, so each ellipsoid $Q$ is represented by a point
in the triangle $P = \{ (b,c) \in \Rset^2 : 0 < c < b < 1 \}$.
Let $P^\ast_1$, $P^\ast_2$, $P^\ast_3$, and $P^\ast_4$ be the four regions
of $P$ that correspond to ellipsoids with minimal periodic trajectories
of type EH1, H1H1, EH2, and H1H2, respectively.
They are shown in figure~\ref{fig:BifurcationsMinimal}.
Their shapes are described below.

\begin{num}\label{num:BifurcationsP}
Let $r = (3-\sqrt{5})/2 \approx 0.382$, $b^\ast_1 = b^\ast_2 = 1$,
and $b^\ast_3 = b^\ast_4 = 1/2$.
Then
\[
P^\ast_j =
\left \{ (b,c) \in P :  b <  b^\ast_j, \; c <  g^\ast_j(b) \right \},
\quad  1 \le j \le 4,
\]
for some continuous functions $g^\ast_j: [0,b^\ast_j] \to \Rset$ such that
\begin{enumerate}
\item
$g^\ast_1$ is concave increasing in $[0,1]$, $0 < g^\ast_1(b) < b$ for all $b\in(0,1)$,
and $g^\ast_1(1) = r$;
\item
$g^\ast_2$ is concave increasing in $[0,1]$,
$g^\ast_1(b) < g^\ast_2(b) < b$ for all $b\in(0,1)$,
and $g^\ast_2(1) = 1/2$;
\item
$g^\ast_3$ is the identity in $[0,r]$, concave decreasing in $[r,1/2]$,
and $g^\ast_3(1/2) = 0$; and
\item
$g^\ast_4$ is increasing in $[0,1/3]$,
concave decreasing in $[1/3,1/2]$,
$3b /4 < g^\ast_4(b) < b$ for all $b\in(0,1/3)$,
$0<g^\ast_4(b) < g^\ast_3(b)$ for all $b\in(1/3,1/2)$,
$g^\ast_4(1/3)=1/4$, and $g^\ast_4(1/2) = 0$.
\end{enumerate}
\end{num}

The functions $g^\ast_j$ can be explicitly expressed by means of
algebraic formulae.
We shall prove that $g^\ast_2(b) = b/(1+b)$ in
proposition~\ref{pro:BifurcationsH1H1}.
We shall study the other three functions in another
paper~\cite{RamirezRos1}, because the techniques change drastically.
A generalized Cayley-like condition is the main tool.
For instance,
\[
g^\ast_4(b) =
\cases{\big(1-b/2-\sqrt{b(1-3b/4)}\big)b/(1-b)^2 & \mbox{for $0 \le b \le 1/3$}\\
       (1-2b)b/(1-b)^2 & \mbox{for $1/3 \le b \le 1/2$}}.
\]
This function $g^\ast_4$ is not concave in the interval $[0,1/3]$.

We shall describe in section~\ref{sec:3D} the regions corresponding
to ellipsoids that have periodic trajectories with given winding numbers
(or quasiperiodic trajectories with given frequencies)
for the four caustic types.
Those general regions have the same shape as these ``minimal'' regions.
That is, they are below the graphs of some functions with properties
similar to the ones stated previously.
Therefore, we discover a general principle.
\emph{The more spheric is an ellipsoid,
      the poorer are its four types of billiard dynamics.}
Here, spheric means $(b,c) \approx (1,1)$.
We quantify this principle in
propositions~\ref{pro:Criterions} and~\ref{pro:Criterions2}.

The key step for the numerical computation of these regions
is to explicitly extend the frequency map for singular values
of the caustic parameters.
The extension is ``exponentially sharp'' at some points,
which implies another general principle.
\emph{The billiard trajectories with some almost singular caustic are ubiquitous.}
We shall enlighten it in subsection~\ref{ssec:Singular3D}
by giving a quantitative sample.
The minimal periodic trajectories shown in
figure~\ref{fig:PeriodicTrajectories3D} reinforce it.

Finally, we want to mention that there exist many remarkable results
about periodic trajectories in other billiard and geodesic problems.
For instance,
several nice algebraic closed geodesics on a triaxial ellipsoid can be seen
in~\cite{Fedorov2005},
and a Cayley-like condition for billiards \emph{on} quadrics was
established in~\cite{AbendaFedorov2006}.
Some results stray from any integrable setup.
For example,
some general lower bounds on the number of periodic billiard
trajectories inside strictly convex smooth hypersurfaces can be found
in~\cite{Babenko1992,Farber2002,FarberTabachnikov2002}.
The planar case was already solved by Birkhoff~\cite{Birkhoff1927}.
Of course, these lower bounds are useless for integrable systems where
the periodic trajectories are organized in continuous families.

We complete this introduction with a note on the organization of the paper.
In section~\ref{sec:Preliminaries} we review briefly some well-known results
about billiards inside ellipsoids in order to fix notations that will be
used along the rest of the paper.
Next, the frequency map is introduced and interpreted in
section~\ref{sec:FrequencyMap}.
This section, concerning ellipsoids of $\Rset^{n+1}$,
also contains two conjectures and the lower bounds on the periods.
Billiards inside ellipses of $\Rset^2$ and
inside triaxial ellipsoids of $\Rset^3$
are thoroughly studied in sections~\ref{sec:2D} and \ref{sec:3D},
respectively.
Billiards inside nondegenerate ellipsoids of $\Rset^{n+1}$ are
revisited in section~\ref{sec:nD}.
Some technical lemmas have been relegated to the appendices.

\section{Preliminaries}
\label{sec:Preliminaries}

In this section details are scarce and technicalities are avoided.
Experts can simply browse this section.
We will list several basic references for the more novice readers.

\subsection{Confocal quadrics and elliptic billiards}

The following results go back to Jacobi, Chasles, Poncelet,
and Darboux.

The starting point of our discussion is the $n$-dimensional
nondegenerate ellipsoid
\begin{equation}
\label{eq:Q}
Q =
\left\{
x=(x_1,\ldots,x_{n+1}) \in \Rset^{n+1} :
\sum_{i=1}^{n+1} \frac{x_i^2}{a_i} = 1
\right\},
\end{equation}
where $a_1,\ldots,a_{n+1}$ are some fixed real constants such that
$0 < a_1 < \cdots < a_{n+1}$.
The degenerate cases in which the ellipsoid has some symmetry
of revolution are not considered here.
This ellipsoid is an element of the family of \emph{confocal quadrics}
given by
\[
Q_\mu =
\left\{
x=(x_1,\ldots,x_{n+1}) \in \Rset^{n+1} :
\sum_{i=1}^{n+1} \frac{x_i^2}{a_i - \mu } = 1
\right\},\qquad \mu \in \Rset.
\]
The meaning of $Q_\mu$ is unclear in the singular cases
$\mu \in \{a_1,\ldots,a_{n+1} \}$.
In fact, there are two natural choices for the
singular confocal quadric $Q_\mu$ when $\mu=a_j$.
The first choice is to define it as the $n$-dimensional coordinate hyperplane
\[
H_j =
\left\{ x=(x_1,\ldots,x_{n+1}) \in \Rset^{n+1} : x_j = 0 \right\},
\]
but it also makes sense to define it as the $(n-1)$-dimensional
\emph{focal quadric}
\[
F_j =
\left\{
x=(x_1,\ldots,x_{n+1}) \in \Rset^{n+1} :
x_j = 0 \mbox{ and } \sum_{i \neq j} \frac{x_i^2}{a_i - a_j} = 1
\right\},
\]
which is contained in the hyperplane $H_j$.
Both choices fit in the framework of elliptic billiards,
but we shall use the notation $Q_{a_j} = H_j$ along this paper.

\begin{thm}[\cite{Moser1980,KozlovTreshchev1991,Audin1994,Tabachnikov1995}]
\label{thm:JacobiChasles}
Once fixed a nondegenerate ellipsoid $Q$,
a generic line $\ell \subset \Rset^{n+1}$ is tangent to exactly $n$
distinct nonsingular confocal quadrics $Q_{\lambda_1},\ldots,Q_{\lambda_n}$
such that $\lambda_1 < \lambda_2 < \cdots < \lambda_n$,
$\lambda_1 \in (-\infty,a_1)  \cup (a_1,a_2)$,
and $\lambda_i \in (a_{i-1},a_i)  \cup (a_i,a_{i+1})$, for $i=2,\ldots,n$.
\end{thm}

Set $a_0=0$.
If a generic line $\ell$ has a transverse intersection with the ellipsoid $Q$,
then $\lambda_1 > 0$, so $\lambda_1 \in (a_0,a_1) \cup (a_1,a_2)$.
The value $\lambda_1 = 0$ is attained just when $\ell$ is tangent to $Q$.
A line is generic in the sense of the theorem if and only if it is
neither tangent to a singular confocal quadric\footnote{By abuse
of notation, it is said that a line
is tangent to the singular confocal quadric $Q_{a_j}$
when it is contained in the coordinate hyperplane $H_j$ or
when it passes through the focal quadric $F_j$.}
nor contained in a nonsingular confocal quadric.

If two lines obey the reflection law at a point $x \in Q$,
then both lines are tangent to the same confocal quadrics~\cite{Tabachnikov1995}.
This shows a tight relation between elliptic billiards and confocal quadrics:
all lines of a billiard trajectory inside the ellipsoid $Q$
are tangent to exactly $n$ confocal quadrics $Q_{\lambda_1},\ldots,Q_{\lambda_n}$,
which are called \emph{caustics} of the trajectory.
We will say that $\lambda = (\lambda_1,\ldots,\lambda_n) \in \Rset^n$
are the \emph{caustic parameters} of the trajectory.

\begin{defi}
A billiard trajectory inside a nondegenerate ellipsoid of the Euclidean space
$\Rset^{n+1}$
is \emph{nonsingular} when it has $n$ distinct nonsingular caustics;
that is,
when its caustic parameter belongs to the \emph{nonsingular caustic space}
\begin{equation}
\label{eq:NonsingularSet}
\Lambda =
\left\{
(\lambda_1,\ldots,\lambda_n) \in \Rset^n :
   \begin{array}{l}
    0 < \lambda_1 < \lambda_2 < \cdots < \lambda_n \\
     \lambda_i \in (a_{i-1},a_i)  \cup (a_i,a_{i+1}) \mbox{ for } 1 \le i \le n
   \end{array}
\right\}.
\end{equation}
\end{defi}

We will only deal with nonsingular billiard trajectories along this paper.
We denote the $2^n$ open connected components of the nonsingular caustic space
as follows:
\[
\Lambda_\sigma =
\left\{
(\lambda_1,\ldots,\lambda_n) \in \Rset^n :
   \begin{array}{l}
    0 < \lambda_1 < \lambda_2 < \cdots < \lambda_n \\
    \lambda_i \in (a_{i+\sigma_i-1}, a_{i+\sigma_i}) \mbox{ for } 1 \le i \le n
   \end{array}
\right\},
\]
for $\sigma = (\sigma_1,\ldots,\sigma_n) \in \{0,1\}^n$.
For instance, the first caustic $Q_{\lambda_1}$ is an ellipsoid if and only if
$\lambda_1 \in (a_0,a_1)$; that is, if and only if $\lambda \in \Lambda_\sigma$
with $\sigma_1=0$.
We will draw the space $\Lambda$ for ellipses and triaxial ellipsoids
of $\Rset^3$ in sections~\ref{sec:2D} and~\ref{sec:3D}, respectively.

\begin{thm} \label{thm:Poncelet}
If a nonsingular billiard trajectory closes after $m_0$ bounces,
all trajectories sharing the same caustics close after $m_0$ bounces.
\end{thm}

Poncelet proved this theorem for conics~\cite{Poncelet1822}.
Darboux generalized it to triaxial ellipsoids of $\Rset^3$.
Later on, this result was generalized to any dimension
in~\cite{ChangFriedberg1988,Chang_etal1993a,Chang_etal1993b,Previato1999}.

\begin{thm}\label{thm:Cayley}
The nonsingular billiard trajectories sharing the caustics
$Q_{\lambda_1},\ldots,Q_{\lambda_n}$ close after $m$ bounces
---up to the action of the group of symmetries
   $G = (\Zset/2\Zset)^{n+1}$ of $Q$---,
if and only if $m\ge n+1$ and
\begin{equation}
\label{eq:Cayley}
\Rank
\left(
\begin{array}{ccc}
h_{m + 1}  & \cdots & h_{n + 2} \\
\vdots   &        & \vdots \\
h_{2m - 1} & \cdots & h_{n + m}
\end{array}
\right) < m - n,
\end{equation}
where
$\sqrt{(a_1-s)\cdots(a_{n+1}-s)(\lambda_1-s)\cdots(\lambda_n-s)} =
 h_0 + h_1 s + h_2 s^2 + \cdots$.
\end{thm}

The group $G$ is formed by the $2^{n+1}$ reflections
---involutive linear transformations---
with regard to coordinate subspaces.
The phrase
``a billiard trajectory closes after $m$ bounces up to the action of $G$''
means that if $(q_k)_{k \in \Zset}$ is the sequence of impact points of the
trajectory, then there exists a reflection $g \in G$
such that $q_{k+m} = g(q_k)$ for all $k\in\Zset$.
Hence, billiard trajectories that close after $m$ bounces up to
the action of the group $G$, close after $m_0 = m$ or $m_0 = 2m$ bounces,
because $q_{k+2m} = g(q_{k+m}) = g^2(q_k) = q_k$.

Cayley proved this theorem for conics~\cite{Cayley1854}.
Later on, this result was generalized to any dimension by Dragovi\'{c}
and Radnovi\'{c} in~\cite{DragovicRadnovic1998a,DragovicRadnovic1998b}.

\subsection{Complete integrability of elliptic billiards}

We recall some results obtained by Liouville, Arnold,
Moser, and Kn\"orrer.

A symplectic map $f:M \to M$ defined on a
$2n$-dimensional symplectic manifold
is \emph{completely integrable} if there exist
some smooth $f$-invariant functions $I_1,\ldots,I_n:M \to \Rset$
(the \emph{integrals}) that are in involution
---that is, whose pair-wise Poisson brackets vanish---
and that are functionally independent almost everywhere on
the phase space $M$.
In this context,
the map $I=(I_1,\ldots,I_n) : M \to \Rset^n$ is called
the \emph{momentum map}.
A point $m\in M$ is a \emph{regular point} of the momentum map
when the $n$-form $\rmd I_1  \wedge \ldots \wedge \rmd I_n$ does not vanish at $m$.
A vector $\lambda \in \Rset^n$ is a \emph{regular value} of the momentum map
when every point in the level set $I^{-1}(\lambda)$ is regular,
in which case the level set is a Lagrangian submanifold of $M$
and we say that $I^{-1}(\lambda)$ is a \emph{regular level set}.

The following result is a discrete version of the Liouville-Arnold
theorem.

\begin{thm}[\cite{Veselov1991}]
Any compact connected component of a regular level set $I^{-1}(\lambda)$
is diffeomorphic to $\Tset^n$, where $\Tset = \Rset/\Zset$.
In appropriate coordinates the restrictions of the map
to this torus becomes a parallel translation $\theta \mapsto \theta + \omega$.
The map $\lambda \mapsto \omega$ is smooth at the regular values of
the momentum map.
\end{thm}

Thus the phase space $M$ is almost foliated by Lagrangian invariant tori
and the map on each torus is simply a parallel translation.
These tori are called \emph{Liouville tori},
the shift $\omega$ is the \emph{frequency} of the torus,
and the map $\lambda \mapsto \omega$ is the \emph{frequency map}.
The dynamics on a Liouville torus with frequency $\omega$
is $m_0$-periodic if and only if $m_0 \omega \in \Zset^n$.
Liouville tori become just invariant curves when $n=1$,
in which case the shift is usually called the
\emph{rotation number} of the invariant curve,
and denoted by $\rho$, instead of $\omega$.

Now, let $Q$ be a (strictly) convex smooth hypersurface of $\Rset^{n+1}$
diffeomorphic to the sphere $\Sset^n$, not necessarily an ellipsoid.
The billiard motion inside $Q$ can be modelled by means of a
symplectic diffeomorphism defined on the $2n$-dimensional phase space
\begin{equation}\label{eq:PhaseSpace}
M =
\left\{
(q,p) \in Q \times \Sset^n : \mbox{$p$ is directed outward $Q$ at $q$}
\right\}.
\end{equation}
We define the billiard map $f: M \to M$, $f(q,p)=(q',p')$, as follows.
The new velocity $p'$ is the reflection of the old velocity $p$
with respect to the tangent plane $T_q Q$.
That is,
if we decompose the old velocity as the sum of its
tangent and normal components at the surface:
$p = p_\rmt + p_\rmn$ with $p_\rmt \in T_q Q$ and $p_\rmn \in N_q Q$,
then $p' = p_\rmt - p_\rmn = p - 2p_\rmn$.
The new impact point $q'$ is the intersection of the ray
$\{ q + \mu p' :  \mu > 0\}$ with the surface $Q$.
This intersection is unique and transverse by convexity.

Elliptic billiards fit in the frame of the Liouville-Arnold
theorem.

\begin{thm}[\cite{Moser1980,Knorrer1985,Veselov1988,Audin1994}]
\label{thm:MK}
The billiard map associated to the nondegenerate ellipsoid~(\ref{eq:Q})
is completely integrable and the caustic parameters
$\lambda_1,\ldots,\lambda_n$ are the integrals.
The set of regular values of the corresponding momentum map
is given by~(\ref{eq:NonsingularSet}).
\end{thm}

\section{The frequency map}
\label{sec:FrequencyMap}

\subsection{Definition and interpretation}

The rotation number for the billiard inside an ellipse is a quotient of
elliptic integrals; see~\cite{Kolodziej1985,ChangFriedberg1988}.
Explicit formulae for the frequency map of the billiard inside a triaxial
ellipsoid of $\Rset^3$ can be found in~\cite[\S III.C]{Chang_etal1993b}.
An equivalent formula is given in~\cite[\S 5]{PopovTopalov2009}.
Both formulae contain hyperelliptic integrals and they can be effortless
generalized to any dimension.
Since we want to avoid as many technicalities as possible,
we will not talk about Riemann surfaces, basis of holomorphic differential
forms, basis of homology groups, period matrices, or other objects that
arise in the theory of algebraic curves.

The following notations are required to define the frequency map.
Once fixed the parameters $a_1,\ldots,a_{n+1}$ of the ellipsoid,
and the caustic parameters $\lambda_1,\ldots,\lambda_n$, we set
\[
T(s) = \prod_{i=1}^{2n+1} (c_i -s),\quad
\{c_1,\ldots,c_{2n+1}\} =
\{a_1,\ldots,a_{n+1}\} \cup \{\lambda_1,\ldots,\lambda_n\}.
\]
If $\lambda \in \Lambda$,
then $c_1,\ldots,c_{2n+1}$ are pair-wise distinct and positive,
so we can assume that
\begin{equation}\label{eq:Inequalities}
c_0 := 0 < c_1 < \cdots < c_{2n+1}.
\end{equation}
Hence, $T(s)$ is positive in the $n+1$ open intervals $(c_{2j},c_{2j+1})$,
and the improper integrals
\begin{equation}\label{eq:KIntegrals}
K_{ij} = \int_{c_{2j}}^{c_{2j+1}} \frac{s^i \rmd s}{\sqrt{T(s)}},
\qquad  i = 0,\ldots,n-1,
\quad j = 0,\ldots,n
\end{equation}
are absolutely convergent, real, and positive.
We also consider the $n+1$ column vectors
\[
K_j = (K_{0j},\ldots,K_{n-1,j})^t \in \Rset^n.
\]
It is known that vectors $K_1,\ldots,K_n$ are linearly
independent; see~\cite[\S III.3]{Griffiths1989}.

\begin{defi}\label{defi:FrequencyMap}
The \emph{frequency map} $\omega:\Lambda \to \Rset^n$ of the billiard
inside the nondegenerate ellipsoid $Q$ is the unique solution of the
system of $n$ linear equations
\begin{equation}\label{eq:Frequency}
K_0 + 2 \sum_{j=1}^n (-1)^j \omega_j K_j = 0.
\end{equation}
\end{defi}

\begin{remark}\label{remark:Homogeneous}
Sometimes it is useful to think that the frequency $\omega$
depends on the parameter $c=(c_1,\ldots,c_{2n+1}) \in \Rset_+^{2n+1}$,
and not only on the caustic parameter
$\lambda=(\lambda_1,\ldots,\lambda_n) \in \Lambda$.
In such situations, we will write $\omega=\varpi(c)$.
This map $c \mapsto \varpi(c)$ is homogeneous of degree zero and analytic
in the domain defined by inequalities~(\ref{eq:Inequalities}).
Homogeneity is deduced by performing a change of scale in
the integrals~(\ref{eq:KIntegrals}).
Hence, we can assume without loss of generality that $c_{2n+1} = a_{n+1} = 1$.
Analyticity follows from the fact that the integrands
in~(\ref{eq:KIntegrals}) are analytic with respect to
the variable of integration in all the intervals of integration and
with respect to $c$, as long as condition~(\ref{eq:Inequalities}) takes place.
\end{remark}

This definition coincides with the formulae contained
in~\cite{Chang_etal1993b,PopovTopalov2009} for $n=2$.
It is motivated by the characterization of periodic billiard
trajectories contained in the next theorem.
The factor $2$ has been added to simplify the interpretation
of the components of the frequency map,
which are all positive, due to the factors $(-1)^j$.

\begin{thm}[\cite{DragovicRadnovic2006,DragovicRadnovic2008}]
\label{thm:DR}
The nonsingular billiard trajectories inside the nondegenerate ellipsoid $Q$
are periodic with exactly $m_j$ points at $Q_{c_{2j}}$ and $m_j$ points
at $Q_{c_{2j+1}}$ if and only if
$m_0 K_0 + \sum_{j=1}^n (-1)^j m_j K_j = 0$.
\end{thm}

The numbers $m_0,m_1,\ldots,m_n$ that appear in theorem~\ref{thm:DR}
are called \emph{winding numbers}.
The nonsingular billiard trajectories with caustic parameter $\lambda$
are periodic with winding numbers $m_0,m_1,\ldots,m_n$ if and only if
\begin{equation}
\label{eq:FrequencyWinding}
\omega_j(\lambda) = \frac{m_j}{2 m_0} \in \Qset_+,\qquad j=1,\ldots,n.
\end{equation}
We note that $m_0$ is the number of bounces in the ellipsoid $Q = Q_{c_0}$,
so it is the period.

\begin{remark}
\label{rem:WindingNumbers}
The sequence of winding numbers of a nonsingular  periodic billiard
trajectory contains information about how the trajectory folds
in the space $\Rset^{n+1}$.
The following properties can be deduced from results contained
in~\cite[\S 4.1]{DragovicRadnovic2006}.
Here, ``number of \underline{\phantom{m}}'' means
``number of times that any periodic billiard trajectory with those
  caustic parameters do \underline{\phantom{m}} along one period''.
The intervals $(c_{2j},c_{2j+1})$ with $j\neq 0$ can adopt exactly
four different forms, each one giving rise to its own geometric picture.

\begin{enumerate}
\item If $(c_{2j},c_{2j+1}) = (a_j,\lambda_{j+1})$,
      then $m_j$ is the number of crossings with $H_j$,
      so it is even and $m_j/2$ is the number of oscillations around
      the hyperplane $H_j$;
\item If $(c_{2j},c_{2j+1}) = (\lambda_j,a_{j+1})$,
      then $m_j$ is the number of crossings with $H_{j+1}$,
      so it is even and $m_j/2$ is the number of oscillations around
      the hyperplane $H_{j+1}$; 
\item \label{item:PlaneOscillations}
      If $(c_{2j},c_{2j+1}) = (a_j,a_{j+1})$,
      then $m_j$ is the number of (alternate) crossings with $H_j$ and $H_{j+1}$,
      so it is even and $m_j/2$ is the number of rotations that the
      trajectory performs when projected onto the
      $(x_{j},x_{j+1})$-coordinate plane $\pi_j$; and
\item \label{item:CausticOscillations}
      If $(c_{2j},c_{2j+1}) = (\lambda_j,\lambda_{j+1})$,
      then $m_j$ is the number of (alternate) tangential touches with
      $Q_{\lambda_j}$ and $Q_{\lambda_{j+1}}$,
      so it can be even or odd, and it is the number of oscillations
      between both caustics.
\end{enumerate}
\end{remark}

These four properties suggest us the following definitions,
which establish the geometric meaning of the components
of the frequency map.
They change with the open connected components of the
nonsingular caustic space.

\begin{defi}
Let $\omega = (\omega_1,\ldots,\omega_n) : \Lambda  \to \Rset^n$
be the frequency map.
\begin{enumerate}
\item
If $(c_{2j},c_{2j+1}) = (a_j,\lambda_{j+1})$,
then $\omega_j = m_j/2 m_0$ is the \emph{$H_j$-oscillation number};
\item
If $(c_{2j},c_{2j+1}) = (\lambda_j,a_{j+1})$,
then $\omega_j = m_j/2 m_0$ is the \emph{$H_{j+1}$-oscillation number};
\item
If $(c_{2j},c_{2j+1}) = (a_j,a_{j+1})$,
then  $\omega_j = m_j/2 m_0$ is the
\emph{$\pi_j$-rotation number}; and
\item
If $(c_{2j},c_{2j+1}) = (\lambda_j,\lambda_{j+1})$,
$2\omega_j = m_j/m_0$ is the
\emph{$(Q_{\lambda_j},Q_{\lambda_{j+1}})$-oscillation number}.
\end{enumerate}
\end{defi}

\begin{remark}\label{rem:Factor_2}
It is important to notice that (only) when a $m_0$-periodic billiard
trajectory has two caustics of the same type ---that is,
when some interval $(c_{2j},c_{2j+1})$ falls into the fourth case---,
it is possible that $m_0 \omega \notin \Zset^n$,
although then $2 m_0 \omega \in \Zset^n$.
This is due to the factor $2$ that we have added in
the definition of the frequency map.
\end{remark}

Finally, if $Q_{\lambda_1}$ is not an ellipsoid
---that is, if $\lambda_1 > a_1$---, then $c_1=a_1$,
and $m_0$ is the number of crossings with $H_1$, so it is even.
Therefore, the following corollary holds.

\begin{cor}\label{cor:Odd}
Among all the nonsingular billiard trajectories inside a nondegenerate
ellipsoid,
only those with an ellipsoid as caustic can have odd period.
\end{cor}

\subsection{Two conjectures}\label{ssec:TwoConjectures}

\begin{con}
\label{con:FrequencyLocal}
The frequency map is a local diffeomorphism;
i.e., it is nondegenerate:
\[
\det \left(
     \frac{\partial \omega_j}{\partial \lambda_i}(\lambda)
     \right)_{1 \le i,j \le n} \neq 0,
\qquad \forall \lambda \in \Lambda.
\]
\end{con}

This conjecture has several relevant consequences throughout the paper.
Popov and Topalov~\cite{PopovTopalov2009} have shown that the frequency map
is \emph{almost everywhere} nondegenerate when $Q$ is a triaxial ellipsoid
of $\Rset^3$,
although they only consider the components $\Lambda_\sigma$
such that $\sigma_1 = 0$.
The nondegeneracy of the frequency map is important because it is an
essential hypothesis
---although we acknowledge that it can be replaced by some weaker
   R\"ussmann-like nondegeneracy conditions~\cite[\S 2]{Russmann2001}---
in most KAM-like theorems,
which are the standard tool to prove the persistence of Liouville tori
under small smooth perturbations of completely integrable maps.
Therefore, if conjecture~\ref{con:FrequencyLocal} holds,
we can conclude that most of the Liouville tori of the billiard phase space
persist under small smooth perturbations of the ellipsoid.
We shall present evidences for this conjecture in
sections~\ref{sec:2D} and~\ref{sec:3D}.

\begin{con}
\label{con:WindingNumbers}
Winding numbers of nonsingular periodic billiard trajectories
are ordered in a strict decreasing way.
More concretely,
\[
2 \le m_n < \cdots < m_1 < m_0 = {\rm period}.
\]
\end{con}

Inequality $m_n \ge 2$ is immediate, because $c_{2n+1}=a_{n+1}$,
so $m_n$ is even.
Inequalities $m_j \le m_0$ for $j \ge 1$ are also immediate,
because the number of crossings with any fixed hyperplane or the number
of tangential touches with any fixed caustic can not exceed the
number of segments of the periodic billiard trajectory.
The strict inequalities $m_j < m_0$ could be also established
(using the symmetries of the ellipsoid), but we skip the details,
since such a small improvement is not worth the effort.

If both conjectures hold, the components of the frequency map
are also ordered in a strict decreasing way and lie in the range $(0,1/2)$;
that is,
\begin{equation}\label{eq:FrequencyOrder}
0 < \omega_n(\lambda) < \cdots < \omega_1(\lambda) < 1/2,
\qquad \forall \lambda \in \Lambda.
\end{equation}
To prove this, we note that
conjecture~\ref{con:FrequencyLocal} implies that
$\omega(\Lambda)$ is open in $\Rset^n$
and $\omega^{-1}(\Qset^n)$ is dense in $\Lambda$,
whereas conjecture~\ref{con:WindingNumbers} and
relation~(\ref{eq:FrequencyWinding}) imply that
the strict inequalities $0 < \omega_n < \cdots < \omega_1 < 1/2$
hold for rational frequencies.
Therefore,
$0 \le \omega_n(\lambda) \le \cdots \le \omega_1(\lambda) \le 1/2$
for any $\lambda \in \Lambda$,
but these inequalities must be strict because $\omega(\Lambda)$ is open.
We have numerically checked that inequalities~(\ref{eq:FrequencyOrder})
hold for thousands of random choices of
$a_1,\ldots,a_{n+1},\lambda_1,\ldots,\lambda_n$
in ``dimensions'' $n \le 5$.
The details of the experiments for $n=2$ are presented in
section~\ref{sec:3D}.

\subsection{Lower bounds on the periods}
\label{ssec:LowerBounds}

We know from theorem~\ref{thm:Cayley} that the period $m_0$ of any nonsingular
periodic billiard trajectory inside a nondegenerate ellipsoid
$Q \subset \Rset^{n+1}$ verifies that $m_0 \ge n+1$.
This result can be improved in several ways using the ordering of the
winding numbers stated in conjecture~\ref{con:WindingNumbers}.
For instance, the global lower bound $m_0 \ge n+2$ follows directly.
We present below more refined semi-global lower bounds,
holding each one on a different open connected component of the
nonsingular caustic space.
They are obtained by realizing that some winding numbers
must be even in agreement with the first items in remark~\ref{rem:WindingNumbers}.
The lower bound associated to some connected component reaches
the value $2n+2$, which doubles the original lower bound given in
theorem~\ref{thm:Cayley}.

\begin{thm}\label{thm:LowerBounds}
Given any $\sigma = (\sigma_1,\ldots,\sigma_n) \in \{0,1\}^n$,
let $E_\sigma$ be the subset of $\{0,1,\ldots,n\}$ such that:
1) $0 \in E_\sigma \Leftrightarrow \sigma_1 = 1$;
2) $j \in E_\sigma \Leftrightarrow (\sigma_j,\sigma_{j+1}) \neq (1,0)$; and
3) $n \in E_\sigma$.
\begin{enumerate}
\item
If $m_0,\ldots,m_n$ are the winding numbers of a periodic trajectory
with caustic parameter in $\Lambda_\sigma$,
then $m_j$ is even for all $j\in E_\sigma$.
\item
If conjecture~\ref{con:WindingNumbers} holds,
any periodic billiard trajectory inside a nondegenerate ellipsoid
of $\Rset^{n+1}$ whose caustic parameter is in $\Lambda_\sigma$
has period at least
\[
\varkappa(\sigma) :=
\min \left\{ m_0 :
\begin{array}{l}
\mbox{$\exists\; 2 \le m_n < \cdots < m_0$ sequence of integers} \\
\mbox{such that $m_j$ is even for any $j \in E_\sigma$}
\end{array}
\right\}.
\]
\item
Let $\mathbf{1}=(1,\ldots,1)\in \{0,1\}^n$,
$\varsigma=(\ldots,0,1,0,1,0)\in \{0,1\}^n$, and $\sigma \in \{0,1\}^n$.
Then
\[
n+2 = \varkappa(\varsigma) < \varkappa(\sigma) < \varkappa(\mathbf{1}) = 2n+2,
\qquad \forall \; \sigma \neq \mathbf{1},\varsigma.
\]
\end{enumerate}
\end{thm}

\proof
(i)
We recall that $m_j$ must be even in the three first cases listed
in remark~\ref{rem:WindingNumbers}.
This is the key property.
For instance, $m_n$ is always even because $c_{2n+1}=a_{n+1}$.
If $\sigma_1 =  1$, then $\lambda_1 \in (a_1,a_2)$ and $c_1=a_1$, so $m_0$ is even.
If $m_j$ is odd,
then $(c_{2j},c_{2j+1}) = (\lambda_j,\lambda_{j+1})$ and
$\lambda_j,\lambda_{j+1} \in (a_j,a_{j+1})$,
so $(\sigma_j,\sigma_{j+1}) = (1,0)$.
Hence,
we have seen that $(\sigma_j,\sigma_{j+1}) \neq (1,0) \Rightarrow m_j$ is even.

(ii)
This follows directly from the previous item and the definition of $\varkappa(\sigma)$.

(iii)
First, we note that $E_\mathbf{1} = \{0,\ldots,n\}$ and
$E_\varsigma = \{\ldots,n-4,n-2,n\}$.
Therefore,
\(
\varkappa(\mathbf{1}) =
\min \left\{ m_0 :
\mbox{$\exists\; 2 \le m_n < \cdots < m_0$ sequence of even numbers}
\right\} = 2n+2
\), and
\[
\varkappa(\varsigma) =
\min \left\{ m_0 :
\begin{array}{l}
\mbox{$\exists\; 2 \le m_n < \cdots < m_0$ sequence of} \\
\mbox{integers s.t. $m_n,m_{n-2},\ldots$ are even}
\end{array}
\right\} = n+2.
\]

The minimum value of $m_0$ among all integer sequences such that
$2 \le m_n < \cdots < m_0$ is attained at the sequence
$m_j=n+2-j$, $0 \le j \le n$.
Thus, $\varkappa(\sigma) \ge n+2$ for all $\sigma \in \{0,1\}^n$, and
$\varkappa(\sigma) = n+2 \Rightarrow E_\sigma = E_\varsigma
\Rightarrow \sigma = \varsigma$.

On the other hand,
$E_\sigma \subset E_{\sigma'} \Rightarrow \varkappa(\sigma) \le \varkappa(\sigma')$.
Hence, $\varkappa(\sigma) \le \varkappa(\mathbf{1}) = 2n+2$
for all $\sigma \in \{0,1\}^n$.
Finally,
$\varkappa(\sigma) = 2n+2 \Rightarrow E_\sigma = \{0,\ldots,n \}
 \Rightarrow \sigma_j \neq 0 \quad \forall j \Rightarrow \sigma = \mathbf{1}$.
\qed

\begin{remark}
Let $\mathbf{0}=(0,\ldots,0)$ and $\bar{\varsigma}=(\ldots,1,0,1,0,1)$.
Then $E_\mathbf{0} = \{1,\ldots,n\}$ and
$E_{\bar{\varsigma}}=\{\ldots,n-5,n-3,n-1,n\}$,
so $\varkappa(\mathbf{0}) = 2n+1$ and $\varkappa(\bar{\varsigma}) = n+3$.
\end{remark}

\begin{remark}
All the lower bounds $\varkappa(\sigma)$ can be explicitly computed when $n=2$.
In that case,
\[
4 = \varkappa(\varsigma) < 
\varkappa(\mathbf{0}) = 5 = \varkappa(\bar{\varsigma}) <
\varkappa(\mathbf{1}) = 6
\]
where $\varsigma=(1,0)$, $\mathbf{0}=(0,0)$,
$\bar{\varsigma}=(0,1)$, and $\mathbf{1}=(1,1)$.
Therefore,
theorem~\ref{thm:LowerBounds3D} about triaxial ellipsoids of $\Rset^3$
is just a particular case of theorem~\ref{thm:LowerBounds}.
It suffices to realize that $\sigma = \mathbf{0}$, $\sigma=\varsigma$,
$\sigma=\bar{\varsigma}$, and $\sigma=\mathbf{1}$ correspond to the cases
EH1, H1H1, EH2, and H1H2, respectively.
\end{remark}

\begin{remark}
The function $\varkappa:\{0,1\}^n \to \{n+2,\ldots,2n+2\}$ is surjective and has
average
\[
\bar{\varkappa} :=
2^{-n} \sum_{\sigma \in \{0,1\}^n} \varkappa (\sigma) = 3n/2 + 2.
\]
We skip the details, the proof is by induction over $n$.
Thus, these semi-global lower bounds improve the
global lower bound $n+2$ by, in average, almost $50\%$.
\end{remark}

Now, a natural question arises.
Are these semi-global lower bounds optimal?
Optimal does not mean that there exists a $\varkappa(\sigma)$-periodic billiard
trajectory whose caustic parameter is in $\Lambda_\sigma$ inside
\emph{all} nondegenerate ellipsoids, but just inside \emph{some} of them.
And we put another question.
Which are the ellipsoids with such ``minimal'' periodic billiard trajectories?
Both questions become almost trivial for ellipses;
see subsection~\ref{ssec:PeriodicTrajectories2D}.
The case of triaxial ellipsoids of $\Rset^3$ was numerically
answered in the introduction.
The general case remains open, but we conjecture that
all these semi-global lower bounds are optimal.

\section{Billiard inside an ellipse}
\label{sec:2D}

In this section we describe the main properties of
the frequency map when $n=1$,
in which case it is called rotation number and denoted by $\rho$.
Many of these properties are old,
but the observation that the the rotation number is exponentially sharp
at the singular caustic parameter seems to be new.
The known results can be found in the
monographes~\cite{KozlovTreshchev1991,Tabachnikov1995}
and the papers~\cite{Kolodziej1985,ChangFriedberg1988,Tabanov1994,Waalkens_etal1997}.

\subsection{Confocal caustics}

To simplify the exposition, we write the ellipse as
\[
Q = \left\{ (x,y) \in \Rset^2 : \frac{x^2}{a} + \frac{y^2}{b} = 1 \right\},
\qquad a > b > 0,
\]
where we could assume, without loss of generality, that $a=1$;
see remark~\ref{remark:Homogeneous}.
Then any nonsingular billiard trajectory inside $Q$
is tangent to one confocal caustic of the form
\[
Q_\lambda =
\left\{
(x,y) \in \Rset^2 : \frac{x^2}{a-\lambda} + \frac{y^2}{b-\lambda} = 1
\right\},
\]
where the caustic parameter $\lambda$ belongs to the nonsingular
caustic space\footnote{When $\lambda \to b^-$ (resp., $\lambda \to b^+$)
the caustic $Q_\lambda$ flattens into the region of the x-axis enclosed
by (resp., outside) the foci of the ellipse $Q$.
When $\lambda \to a^-$, the caustic flattens into the whole y-axis.}
\begin{equation}\label{eq:CausticSpace2D}
\Lambda = E \cup H,\qquad E =(0,b),\qquad H=(b,a).
\end{equation}
We have chosen those names for the connected components of $\Lambda$
because then $Q_\lambda$ is an ellipse for $\lambda \in E$,
and a hyperbola for $\lambda \in H$.

\subsection{Phase portrait}

We describe now the billiard dynamics inside an ellipse.
This description goes back to Birkhoff~\cite[\S VIII.12]{Birkhoff1927},
so it is rather old and we just list the results.
Concretely, we want to know how the phase space is foliated
by Liouville tori
(invariant curves on which the motion becomes a rigid rotation) and
separatrices
(invariant curves on which the motion tends to some hyperbolic
 periodic trajectories).

\begin{figure}
\iffiguresPDF
\begin{center}
\includegraphics[width=6in]{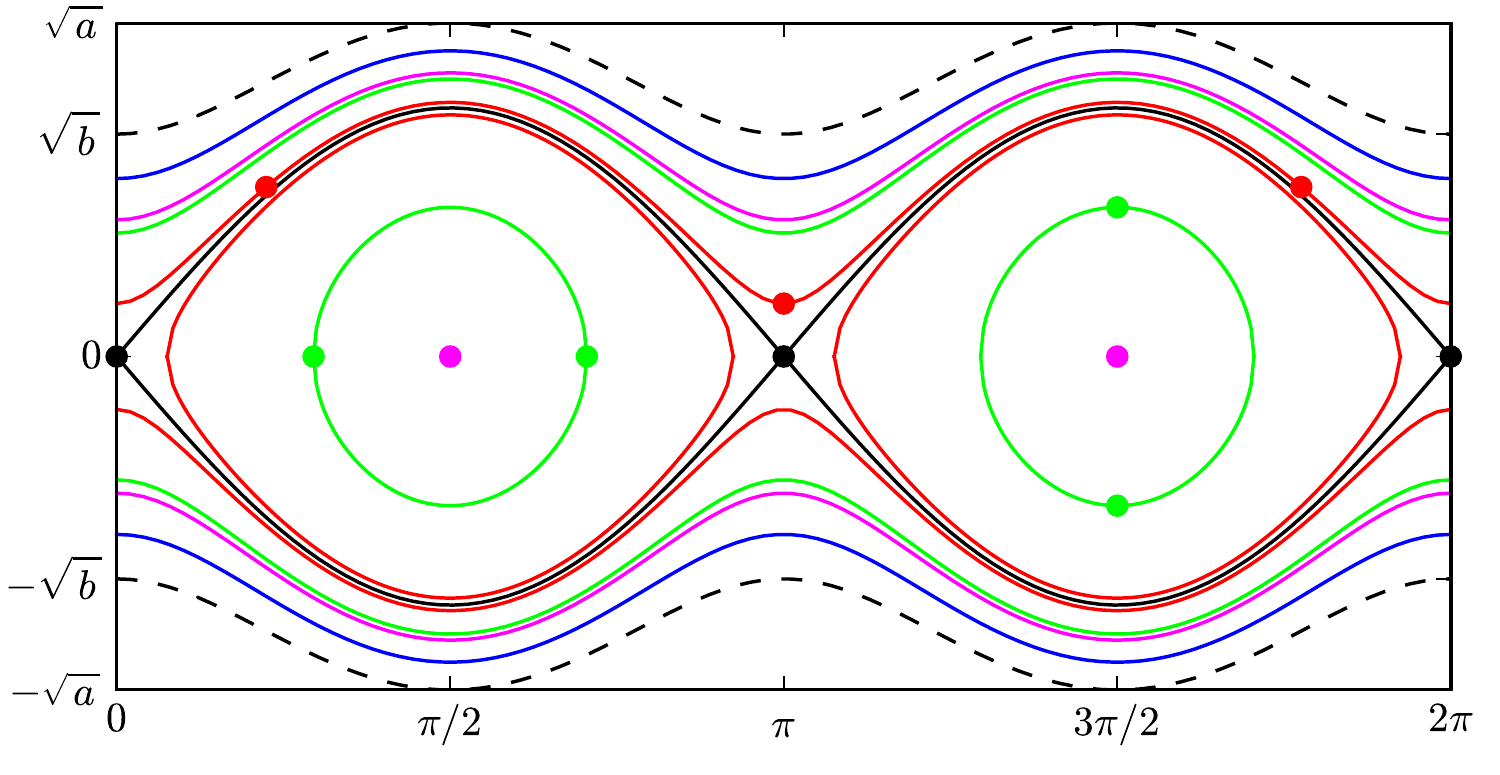}
\end{center}
\else
\vspace{3.275in}
\fi
\caption{Phase portrait of the billiard
map in $(\varphi,r)$ coordinates for $a=1$ and $b=4/9$.
The dashed black lines enclose the phase space~(\ref{eq:PhaseSpace2D}).
The black points are the hyperbolic two-periodic points
corresponding to the oscillation along the major axis of the ellipse.
The black curves are the separatrices of these hyperbolic points.
The magenta points denote the elliptic two-periodic points
corresponding to the oscillation along the minor axis of the ellipse.
The magenta curves are the invariant curves whose rotation number
coincides with the frequency of these elliptic points.
The invariant curves with rotation numbers $1/6$, $1/4$ and $1/3$
are depicted in blue, green and red, respectively.
The red points label a three-periodic trajectory whose caustic is an ellipse.
The green points label a four-periodic trajectory whose caustic is
a hyperbola.}
\label{fig:PhasePortrait}
\end{figure}

Let us put some global coordinates $(\varphi,r)$ over the billiard phase
space $M$ defined in~(\ref{eq:PhaseSpace}), just for visualization purposes.
First, following Birkhoff,
we parameterize the impact points on the ellipse by means of an angular
coordinate $\varphi \in \Tset$.
We take, for instance,
$q = \gamma(\varphi) = (a^{1/2} \cos \varphi,b^{1/2} \sin \varphi)$.
Second, given an outward unitary velocity $p \in \Sset$,
we set $r = \langle \gamma'(\varphi),p \rangle$,
and so $|r| < | \gamma'(\varphi) | = (a \sin^2 \varphi + b \cos^2 \varphi)^{1/2}$.
Then the correspondence $(q,p) \mapsto (\varphi,r)$ allows us to identify
the phase space $M$ with the annulus
\begin{equation}
\label{eq:PhaseSpace2D}
\Aset=
\left\{
(\varphi,r) \in \Tset\times\Rset: r^2 < a \sin^2 \varphi + b \cos^2 \varphi
\right\}.
\end{equation}
In these coordinates,
the caustic parameter becomes $\lambda(\varphi,r) = (a-b)\sin^{2}\varphi +b - r^2$.
The partition of the annulus into invariant level curves
of $\lambda$ is shown in figure~\ref{fig:PhasePortrait}.

Each regular level set contains \emph{two} Liouville curves and represents
the family of tangent lines to a fixed nonsingular caustic $Q_\lambda$.
If $Q_\lambda$ is an ellipse,
each Liouville curve has a one-to-one projection onto the
$\varphi$ coordinate and corresponds to rotations around
$Q_\lambda$ in opposite directions, so they are invariant under $f$.
If $Q_\lambda$ is a hyperbola, then each Liouville
curve corresponds to the impacts on \emph{one} of the two pieces of the
ellipse between the branches of $Q_\lambda$,
so they are exchanged under $f$ and invariant under $f^2$.

The singular level set $\{(\varphi,r) \in \Aset : \lambda(\varphi,r) = b\}$
gives rise to the $\infty$-shaped curve
\[
\lambda^{-1}(b) =
\left\{ (\varphi,r) \in \Aset : r = \pm (a-b)^{1/2} \sin \varphi \right\},
\]
which corresponds to the family of lines through the foci.
This singular level set has rotation number $1/2$;
see~\cite[page 428]{KatokHasselblatt1995}.
The cross points on this singular level represent the
two-periodic trajectory along the major axis of the ellipse,
and the eigenvalues of the differential of the billiard map at them
are positive but different from one:
$\rme^{\pm h}$ with $\cosh^2 h/2 = a/b$ and $h >0$.
On the contrary,
the two-periodic trajectory along the minor axis correspond to the
centers of the regions inside the $\infty$-shaped curve,
and the eigenvalues in that case are conjugate complex of modulus one:
$\rme^{\pm 2\pi \theta \rmi}$ with $\cos^2 \pi \theta = b/a$ and
$0 < \theta < 1/2$.
Therefore, the major axis is a hyperbolic (unstable) two-periodic
trajectory and the minor axis is an elliptic (stable) one.
These are the only two-periodic motions.
The basic results about the stability of two-periodic billiard
trajectories can be found in~\cite{KozlovTreshchev1991,Tabachnikov1995}.

\subsection{Extension and range of the rotation number}

\begin{figure}
\iffiguresPDF
\begin{center}
\includegraphics[width=6in]{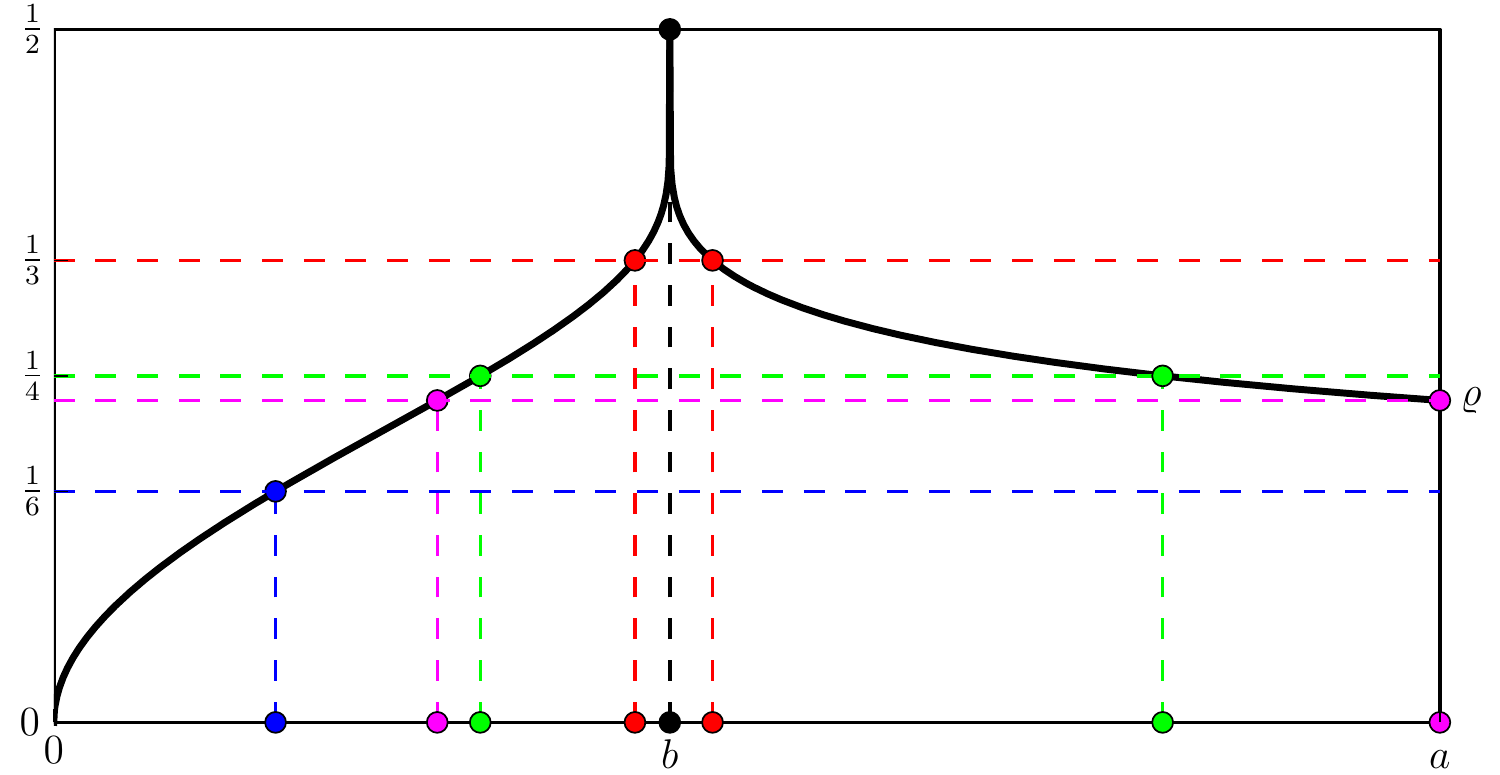}
\end{center}
\else
\vspace{3.275in}
\fi
\caption{The rotation function $\rho(\lambda)$ of the ellipse
for $a=1$ and $b=4/9$.
Colours are taken from figure~\ref{fig:PhasePortrait}.
The parameters $\lambda^0_\pm$ fast approach $b$ as $\rho^0$
tends to $1/2$.}
\label{fig:RotationFunction}
\end{figure}

Let $\rho(\lambda)$ be the \emph{rotation number}
of the billiard trajectories inside the ellipse $Q$ sharing the
nonsingular caustic $Q_\lambda$.
From definition~\ref{defi:FrequencyMap} we get that the function
$\rho:E \cup H \to \Rset$ is given by the quotients of elliptic integrals
\begin{equation}\label{eq:RotationNumber}
\rho(\lambda) = \rho(\lambda;b,a) =
\frac{\int_0^{\min(b,\lambda)} \frac{\rmd s}{\sqrt{(\lambda-s)(b-s)(a-s)}}}
     {2 \int_{\max(b,\lambda)}^a \frac{\rmd s}{\sqrt{(\lambda-s)(b-s)(a-s)}}} =
\frac{\int_\chi^\mu \frac{\rmd t}{\sqrt{t(t-1)(t-\chi)}}}
     {2\int_0^1 \frac{\rmd t}{\sqrt{t(t-1)(t-\chi)}}},
\end{equation}
where the parameters $1 < \chi < \mu$ are given by
$\chi=(a-\underline{m})/(a-\overline{m})$ and $\mu=a/(a-\overline{m})$,
with $\underline{m}=\min(b,\lambda)$ and $\overline{m}=\max(b,\lambda)$.
The second equality follows from the change of variables
$t=(a-s)/(a-\overline{m})$.
The second quotient already appears in~\cite{Chang_etal1993b}.
Other equivalent quotients of elliptic integrals were given
in~\cite{Kolodziej1985,Waalkens_etal1997}.
We have drawn the rotation function $\rho(\lambda)$ in
figure~\ref{fig:RotationFunction},
compare with~\cite[figure 2]{Waalkens_etal1997}.

\begin{pro}\label{pro:2D}
The function $\rho : E \cup H \to \Rset$
given in~(\ref{eq:RotationNumber}) has the following properties.
\begin{enumerate}
\item
It is analytic in $\Lambda = E \cup H$.
\item
It can be continuously extended to the closed interval
$\bar{\Lambda} = \Lambda \cup \partial \Lambda = [0,a]$ with
\[
\rho(0) = 0, \qquad
\rho(b) = 1/2, \qquad
\rho(a) = \varrho,
\]
where the limit value $0 < \varrho < 1/2$ is defined by
$\sin^2 \pi \varrho = b/a$.

\item
Let $\kappa^G$ and $\kappa^S$ be the positive constants given by
\[
\kappa^G =
\left(\sqrt{ab} \int_b^a \frac{\rmd s}{\sqrt{s(s-b)(a-s)}} \right)^{-1},
\qquad
\cosh^2 \kappa^S = a/b.
\]
The asymptotic behavior of $\rho(\lambda)$ at the singular
parameters $\lambda \in \partial \Lambda = \{0,b,a\}$ is:
\begin{enumerate}
\item
$\rho(\lambda) = \kappa^G \lambda^{1/2} + \Or(\lambda^{3/2})$,
as $\lambda \to 0^+$;
\item
$\rho(\lambda) =
 1/2 + \kappa^S/\log |b-\lambda| + \Or\left(1/\log^2 |b-\lambda|\right)$,
as $\lambda \to b$; and
\item
$\rho(\lambda) = \varrho + \Or(a-\lambda)$, as $\lambda \to a^-$.
\end{enumerate}

\item
Given any $\rho^0 \in (\varrho,1/2)$,
let $\lambda^0_-$ be the biggest parameter in $E$
such that $\rho(\lambda^0_-) =\rho^0 $,
and let $\lambda^0_+$ be the smallest parameter in $H$
such that $\rho(\lambda^0_+) = \rho^0$.
Both parameters become exponentially close to the singular
caustic parameter $b$ when $\rho_0$ tends to $1/2$.
In fact,
\[
\lambda^0_\pm =
b \pm 16 (a-b) \rme^{-\kappa^S/(1/2-\rho^0)}
+ \Or\left(\rme^{-2\kappa^S/(1/2-\rho^0)} \right),
\quad \rho^0 \to (1/2)^-.
\]
\end{enumerate}
\end{pro}

\proof
(i)~follows from remark~\ref{remark:Homogeneous}.
The rest of the proof is postponed to~\ref{sap:Proof_2D}.
\qed

\begin{remark}
If conjecture~\ref{con:FrequencyLocal} holds,
then $\rho'(\lambda)$ is positive in $E$ and negative in $H$,
so $\rho(\lambda)$ maps diffeomorphically $E$ onto $(0,1/2)$
and $H$ onto $(\varrho,1/2)$.
In particular, the parameters $\lambda^0_-$ and $\lambda^0_+$ are unique.
The conjecture remains unproven,
but we shall see in proposition~\ref{pro:Increasing} that $\rho(\lambda)$
is increasing in $E$, which suffices to check the unicity of $\lambda^0_-$.
\end{remark}

\begin{remark}
The limit rotation number $\varrho$ is related to the
conjugate complex eigenvalues $\rme^{\pm 2\pi \theta \rmi}$ of the
elliptic two-periodic orbit.
Concretely, $\theta + \varrho= 1/2$.
Besides, $\varrho$ tends to zero when the ellipse flattens
and tends to one half when the ellipse becomes circular.
That is, $\lim_{b/a \to 0^+}\varrho = 0$, and
$\lim_{b/a \to 1^-}\varrho = 1/2$.
\end{remark}

\begin{defi}\label{defi:Extended_rho}
The continuous extension $\rho:[0,a] \to \Rset$ is called the
\emph{(extended) rotation function} of the ellipse $Q$.
\end{defi}

\subsection{Geometric meaning of the rotation number}

Let us assume that the billiard trajectories sharing some
nonsingular caustic $Q_\lambda$ are $m_0$-periodic,
so they describe polygons with $m_0$ sides inscribed in the ellipse $Q$.
Then,
according to theorem~\ref{thm:DR}, equation~(\ref{eq:FrequencyWinding}),
and corollary~\ref{cor:Odd},
it turns out that $\rho(\lambda) = m_1/2m_0$ for some integers
$2 \le m_1 < m_0$ such that $m_1$ is always even whereas $m_0$ can be odd
only when $Q_\lambda$ is an ellipse.
Besides, from the geometric interpretation of the frequency map presented
in section~\ref{sec:FrequencyMap},
we know that:
1) If $Q_\lambda$ is an ellipse,
   the polygons are enclosed between $Q$ and $Q_\lambda$,
   and make $m_1/2$ turns around the origin; and
2) If $Q_\lambda$ is a hyperbola,
   they are contained in the region delimited by $Q$
   and the branches of $Q_\lambda$,
   and cross $m_1$ times the minor axis of the ellipse.

These interpretations can be extended to nonperiodic trajectories.
Concretely,
\[
\rho(\lambda) =
\cases{\lim_{k \to +\infty} n_k/k  \quad & \mbox{if $\lambda \in E$,}\\
       \case{1}{2}\lim_{k\to + \infty} l_k/k \quad & \mbox{if $\lambda \in H$,}}
\]
where $n_k$ (respectively, $l_k$) is the number of turns around the origin
(respectively, crossings of the minor axis) of the first $k$ segments of a
given billiard trajectory with caustic $Q_\lambda$.

\begin{pro}\label{pro:Increasing}
The rotation function $\rho(\lambda)$ is increasing in $E$.
\end{pro}

\proof
Let $\gamma: \Tset \to Q$ be a fixed parameterization of the ellipse $Q$.
Then the billiard dynamics inside $Q$ associated to any caustic $Q_\lambda$,
$\lambda \in E$, induces a circle diffeomorphism $f_\lambda: \Tset \to \Tset$
of rotation number $\rho(\lambda)$.
Let $0 < \lambda_1 < \lambda_2 < b$.
The billiard trajectories sharing the small caustic $Q_{\lambda_2}$
rotate faster than the ones sharing the big caustic $Q_{\lambda_1}$,
so $F_{\lambda_1} < F_{\lambda_2}$ for any two compatible lifts
$F_{\lambda_j}$ of the circle diffeomorphisms $f_{\lambda_j}$.
Then $\rho(\lambda_1) \le \rho(\lambda_2)$;
see~\cite[Proposition~11.1.8]{KatokHasselblatt1995}.
Thus, $\rho(\lambda)$ is nondecreasing and, by analyticity, increasing.
\qed

We have not proved that $\rho(\lambda)$ is decreasing in $H$
because it is not easy to construct an ordered family of
circle diffeomorphisms for caustic hyperbolas.

\subsection{Bifurcations in parameter space}
\label{ssec:Bifurcations2D}

We want to determine all the ellipses $Q = \{ x^2/a + y^2/b = 1\}$,
$0 < b < a$, that have billiard trajectories
with a prescribed rotation number $\rho^0 \in (0,1/2)$ and with
a prescribed type of caustics (ellipses or hyperbolas).
We recall that the rotation function $\rho(\lambda)$
diffeomorphically maps $E$ onto $(0,1/2)$, and $H$ onto $(\varrho,1/2)$.
Therefore, $\rho^0 \in \rho(E)$ for all ellipses $Q$, whereas
\begin{equation}\label{eq:Bifurcations2D}
\rho^0 \in \rho(H) \Leftrightarrow
\varrho < \rho^0 \Leftrightarrow
\sin^2 \pi \varrho < \sin^2 \pi \rho^0 \Leftrightarrow
b < a \sin^2 \pi \rho^0.
\end{equation}
This shows that flat ellipses have more periodic trajectories
than rounded ones.
There exist similar results for triaxial ellipsoids of $\Rset^3$.
See, for instance, propositions~\ref{pro:Criterions} and~\ref{pro:Criterions2}.

\subsection{Examples of periodic trajectories with minimal periods}
\label{ssec:PeriodicTrajectories2D}

\begin{figure}
\iffiguresPDF
\begin{center}
\includegraphics[width=3in]{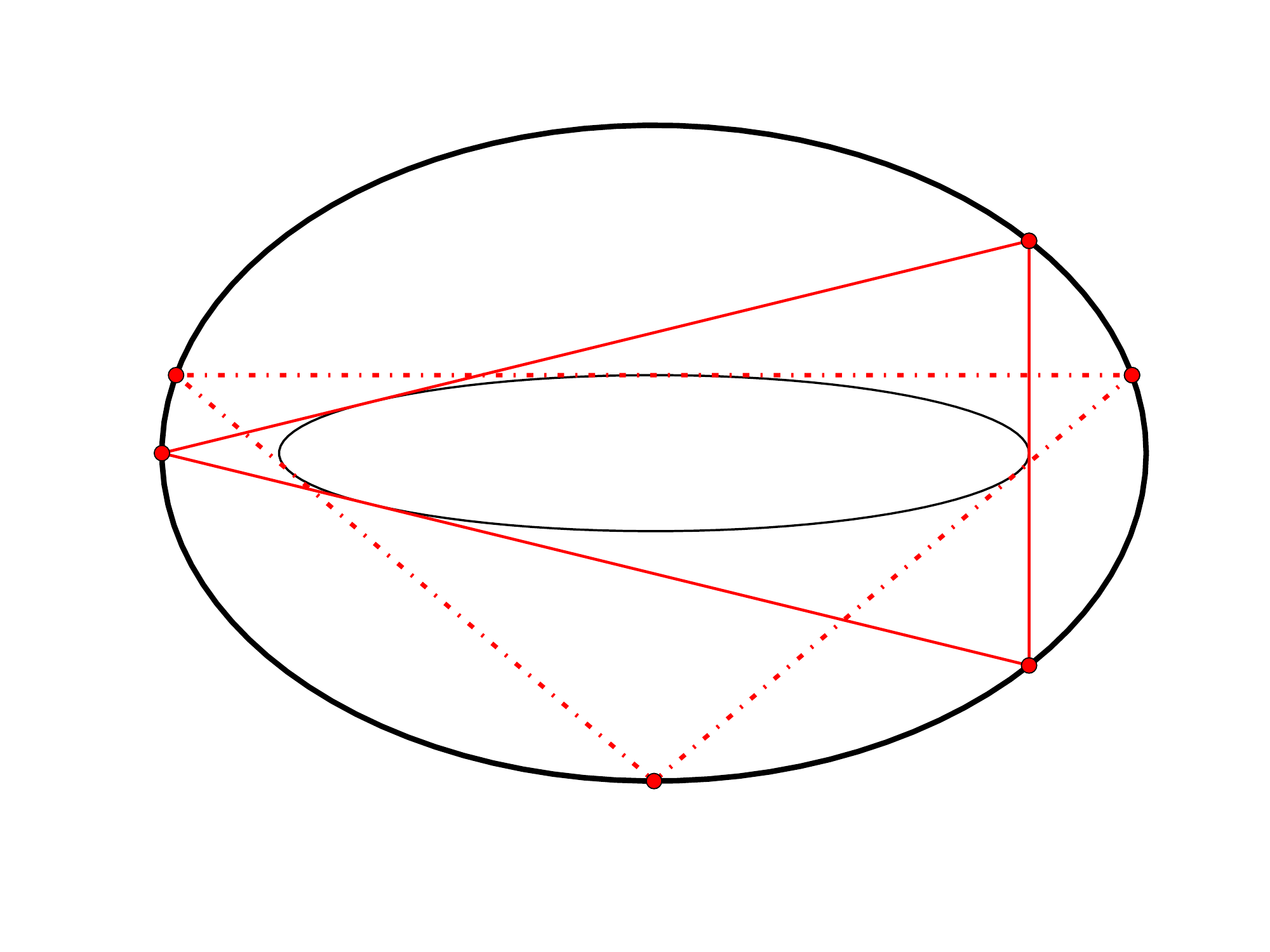}
\includegraphics[width=3in]{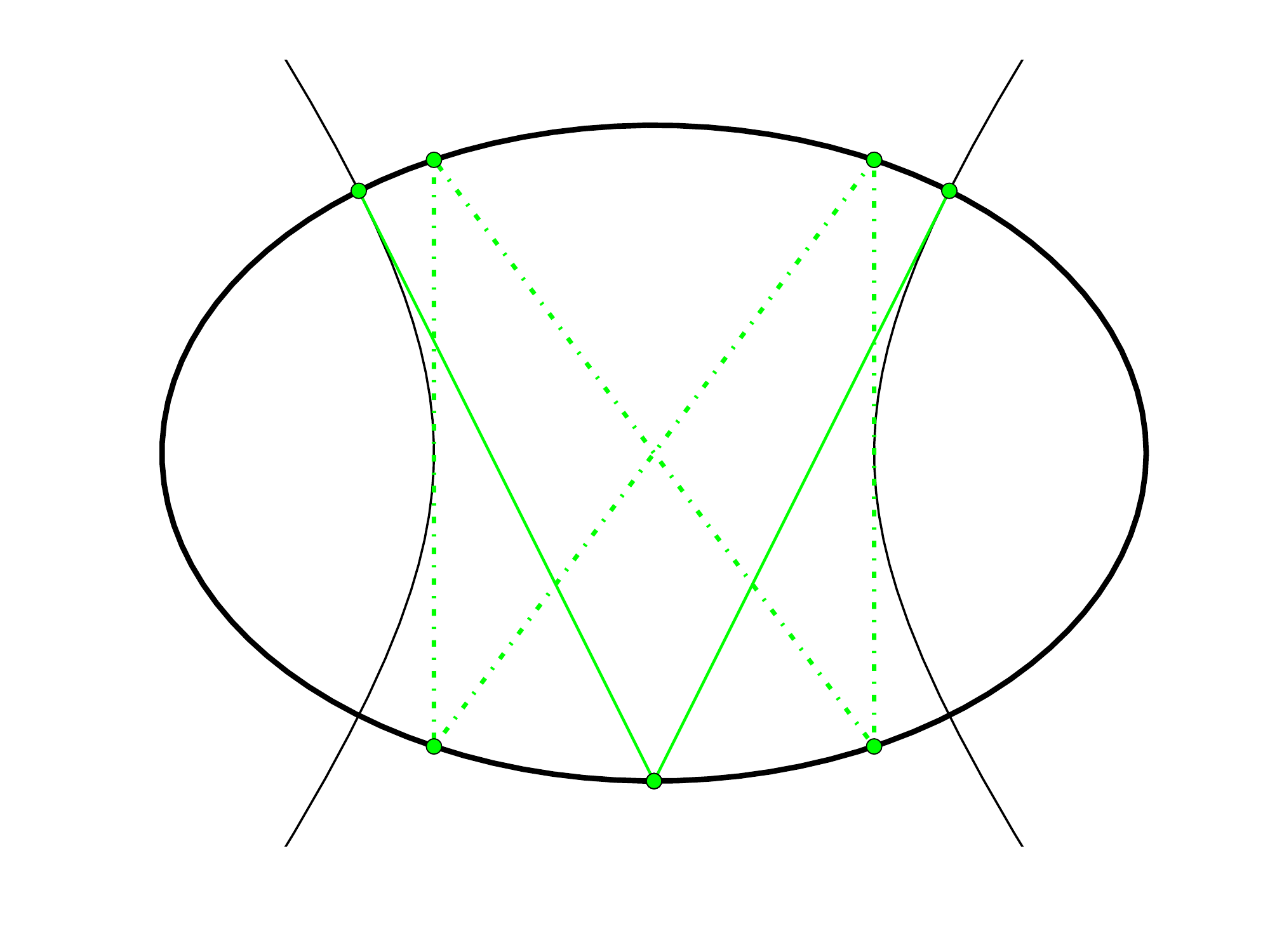}
\end{center}
\else
\vspace{2.47in}
\fi
\caption{Examples of symmetric nonsingular billiard trajectories
with minimal periods for $a=1$ and $b=4/9$.
Left: Period three and the caustic is an ellipse.
Right: Period four and the caustic is a hyperbola.
The continuous lines are reserved for the trajectories
that correspond to the periodic orbits depicted in
figure~\ref{fig:PhasePortrait}.}
\label{fig:PeriodicTrajectories2D}
\end{figure}

The billiard map associated to an ellipse has no fixed points,
its only two-periodic points correspond to the
oscillations along the major or minor axis,
and only the trajectories with an ellipse as caustic can have odd period.
Therefore,
the periodic trajectories with an ellipse as caustic have period
at least three,
whereas the ones with a hyperbola as caustic have period at least four.
These lower bounds are optimal;
see figure~\ref{fig:PeriodicTrajectories2D}.
To be more precise, we set
\begin{equation}\label{eq:MinimalCaustics2D}
\lambda^*_\rmE = \frac{3ab}{a+b+2\sqrt{a^2-ba+b^2}},
\qquad
\lambda^*_\rmH = \frac{ab}{a-b}.
\end{equation}
We note that $\lambda^*_\rmE \in E$ for all $0 < b < a$,
and $\lambda^*_\rmH \in H$ for all $0 < b < a/2$.
The trajectories with caustic $Q_{\lambda^*_E}$ are three-periodic,
the ones with caustic $Q_{\lambda^*_H}$ are four-periodic.
The proof is an elementary exercise in Euclidean geometry.
We leave it to the reader.
Finally,
we deduce from the geometric interpretation of the rotation number given
before that $\rho(\lambda^\ast_\rmE) = 1/3$ and $\rho(\lambda^*_\rmH) = 1/4$.
This second identity explains the restriction $b<a/2$;
see~(\ref{eq:Bifurcations2D}).

\section{Billiard inside a triaxial ellipsoid of $\Rset^3$}
\label{sec:3D}

The previous section sets the basis of this one.
Roughly speaking, we want to follow the same steps
---extension of the frequency map and description of its range---
in order to find the same results
---bifurcations in the parameter space and minimal periodic trajectories.
But the study of ellipsoids is harder,
which has two unavoidable consequences.
First, statements and proofs of the analytical results are more cumbersome.
Second, some results remain unproven, so we shall present numerical experiments
and semi-analytical arguments as support.

\subsection{Confocal caustics}

The caustics of a billiard inside a triaxial ellipsoid are described
in several places.
The representation of the caustic space shown in figure~\ref{fig:CausticSpace}
can also be found in~\cite{Knorrer1985,Waalkens_etal1999,DragovicRadnovic2009}.

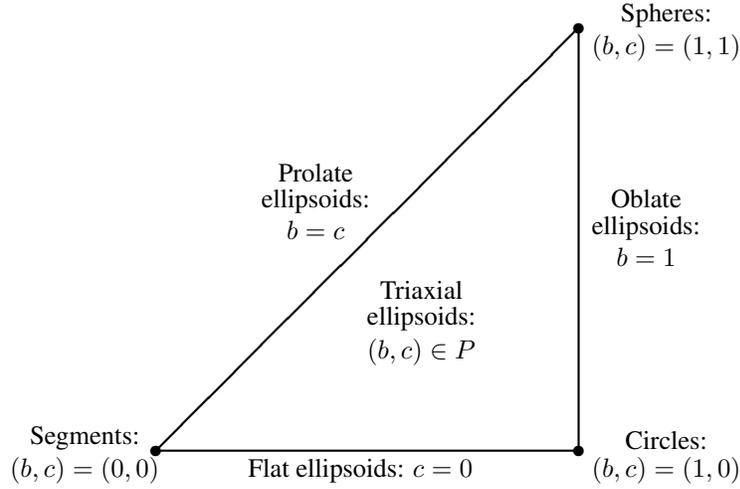
\begin{figure}
\begin{center}
\begin{picture}(200,200)(0,0)
\footnotesize

\put(20,20){\circle*{4}}
\put(180,20){\circle*{4}}
\put(180,180){\circle*{4}}

\thicklines
\put(20,20){\line(1,0){160}}
\put(180,20){\line(0,1){160}}
\put(20,20){\line(1,1){160}}

\put(185,170){\shortstack{Spheres: \\ $(b,c)=(1,1)$}}
\put(60,100){\shortstack{Prolate \\ ellipsoids: \\ $b=c$ }}
\put(185,90){\shortstack{Oblate \\ ellipsoids: \\ $b=1$ }}
\put(100,55){\shortstack{Triaxial \\ ellipsoids: \\ $(b,c) \in P$ }}
\put(55,10){Flat ellipsoids: $c=0$ }
\put(-35,10){\shortstack{Segments: \\ $(b,c)=(0,0)$ }}
\put(185,10){\shortstack{Circles: \\ $(b,c)=(1,0)$}}
\end{picture}
\end{center}
\caption{The triangular parameter space $P$.}
\label{fig:ParameterSpace}
\end{figure}

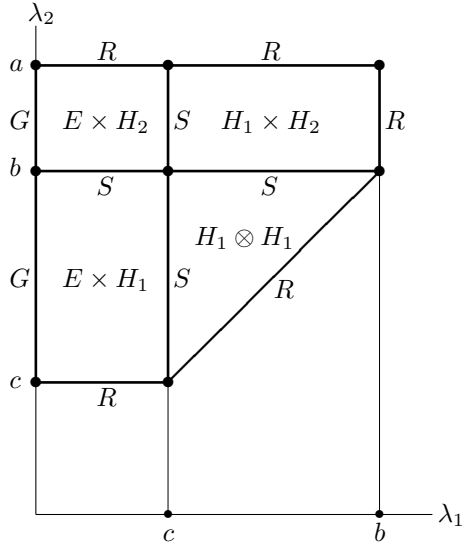
\begin{figure}
\begin{center}
\begin{picture}(150,200)(0,0)
\footnotesize
\put(50,0){\circle*{3}}
\put(130,0){\circle*{3}}
\put(0,50){\circle*{4}}
\put(0,130){\circle*{4}}
\put(0,170){\circle*{4}}
\put(50,50){\circle*{4}}
\put(50,130){\circle*{4}}
\put(50,170){\circle*{4}}
\put(130,130){\circle*{4}}
\put(130,170){\circle*{4}}

\put(-10,48){$c$}
\put(-10,128){$b$}
\put(-10,168){$a$}
\put(-3,187){$\lambda_2$}
\put(48,-10){$c$}
\put(128,-10){$b$}
\put(152,-3){$\lambda_1$}
\put(-10,86){$G$}
\put(-10,146){$G$}
\put(23,41){$R$}
\put(23,172){$R$}
\put(85,172){$R$}
\put(132,146){$R$}
\put(90,82){$R$}
\put(23,121){$S$}
\put(85,121){$S$}
\put(52,86){$S$}
\put(52,146){$S$}
\put(10,86){$E \times H_1$}
\put(10,146){$E \times H_2$}
\put(70,146){$H_1 \times H_2$}
\put(60,102){$H_1 \otimes H_1$}

\thinlines
\put(0,0){\line(1,0){150}}
\put(0,0){\line(0,1){185}}

\linethickness{1pt}
\put(0,50){\line(1,0){50}}
\put(0,130){\line(1,0){130}}
\put(0,170){\line(1,0){130}}
\put(0,50){\line(0,1){120}}
\put(50,50){\line(0,1){120}}
\put(130,130){\line(0,1){40}}
\thicklines
\put(50,50){\line(1,1){80}}

\linethickness{0.1pt}
\put(50,0){\line(0,1){50}}
\put(130,0){\line(0,1){130}}
\end{picture}
\end{center}
\caption{The nonsingular caustic space
$\Lambda = (E \times H_1) \cup (E  \times H_2) \cup (H_1 \otimes H_1)
           \cup (H_1 \times H_2)$ with its border
$\partial \Lambda = G \cup R \cup S \cup \Lambda^0$.}
\label{fig:CausticSpace}
\end{figure}

We write the triaxial ellipsoid as
\[
Q =
\left\{
(x,y,z) \in \Rset^3 :
\frac{x^2}{a} + \frac{y^2}{b} + \frac{z^2}{c} = 1
\right\},
\qquad a > b > c > 0.
\]
We could assume, again without loss of generality, that $a=1$.
Then the parameter space of triaxial ellipsoids in $\Rset^3$
can be represented as the triangle
\begin{equation}\label{eq:ParameterSpace}
P = \left\{ (b,c)\in\Rset^2 : 0 < c < b <1 \right\},
\end{equation}
whose edges represent ellipsoids with a symmetry of revolution
(oblate and prolate ones) or flat ellipsoids,
as illustrated in figure~\ref{fig:ParameterSpace}.
We shall write the statements of the main results for
arbitrary values of $a$, but we shall take $a=1$ in the pictures.

From theorem~\ref{thm:JacobiChasles}, we know that any nonsingular billiard trajectory
inside the ellipsoid $Q$ is tangent to \emph{two} distinct nonsingular
caustics of the confocal family
\[
Q_\lambda =
\left\{
(x,y,z) \in \Rset^3 :
\frac{x^2}{a-\lambda} + \frac{y^2}{b-\lambda} + \frac{z^2}{c-\lambda}
= 1
\right\}.
\]
The caustic $Q_\lambda$ is an ellipsoid for $\lambda \in E$,
a one-sheet hyperboloid when $\lambda \in H_1$,
and a two-sheet hyperboloid if $\lambda \in H_2$, where
\[
E = (0,c),\qquad H_1 = (c,b),\qquad H_2 = (b,a).
\]

In order to have a clearer picture of how these caustics change,
let us explain the situation when $\lambda$ approaches the
singular values $c$, $b$, or $a$.
First,
when $\lambda\to c^-$ (respectively, $\lambda\to c^+$),
the caustic $Q_\lambda$ flattens into the region of the coordinate plane
$\pi_\rmz = \{ z = 0 \}$ enclosed by (respectively, outside)
the \emph{focal ellipse}
\begin{equation}\label{eq:FocalEllipse}
Q^\rmz_c =
\left\{
(x,y,0)\in\Rset^3 :
\frac{x^2}{a-c} + \frac{y^2}{b-c} = 1
\right\}.
\end{equation}
Second, when $\lambda\to b^-$ (resp., $\lambda\to b^+$),
the caustic $Q_\lambda$ flattens into the region of the coordinate plane
$\pi_\rmy = \{ y=0 \}$ between (resp., outside) the branches of the
\emph{focal hyperbola}
\[
Q^\rmy_b =
\left\{
(x,0,z)\in\Rset^{3} :
\frac{x^2}{a-b} - \frac{z^2}{b-c} = 1
\right\}.
\]
Third, the caustic flattens into the whole coordinate plane
$\pi_\rmx = \{ x = 0 \}$ when $\lambda\to a^-$.

We recall that not all combinations of caustics can take place.
For instance, both caustics can not be ellipsoids.
The four possible combinations are denoted by EH1, H1H1, EH2, and H1H2.
Hence, the caustic parameter $\lambda=(\lambda_1,\lambda_2)$
belongs to the nonsingular caustic space
\begin{equation}\label{eq:CausticSpace3D}
\Lambda =
(E \times H_1) \cup (H_1 \otimes H_1) \cup
(E  \times H_2) \cup (H_1 \times H_2),
\end{equation}
where $H_1 \otimes H_1 =
\{ \lambda \in H_1 \times H_1 : \lambda_1 < \lambda_2 \}$.
For instance, $\lambda \in E \times H_1$ for trajectories of type EH1,
which means that $Q_{\lambda_1}$ is an ellipsoid and $Q_{\lambda_2}$
is a one-sheet hyperboloid.

\subsection{The extension of the frequency map}

To begin with,
we extend the frequency map $\omega: \Lambda \to \Rset^2$ to the
borders of the four components of the caustic space~(\ref{eq:CausticSpace3D}),
in the same way that the rotation number was extended to the endpoints
of the two caustic intervals~(\ref{eq:CausticSpace2D}).
The extension depends strongly on the ``piece''
of the border under consideration.
Hence, we need some notations for such ``pieces''.

The set $\Lambda$ is the union of three open rectangles and one open
isosceles rectangular triangle. In total $\Lambda$ has eleven edges and
eight vertexes. We consider the partitions
\[
\partial \Lambda = \Lambda^0 \cup \Lambda^1,\qquad
\Lambda^1 = G \cup R \cup S,
\]
where $\Lambda^1$ is the set of edges,
$\Lambda^0$ is the set of vertexes, and
$S$, $G$ and $R$ are the sets formed by the four inner edges,
the two left edges, and the remaining five edges, respectively.
See figure~\ref{fig:CausticSpace}.
We shall see that the frequency map is quite singular
(in fact, exponentially sharp) at the four edges in $S$,
quite regular at the five edges in $R$,
and it is somehow related to the geodesic flow on the ellipsoid $Q$
at the two edges in $G$.
That motivates the notation.

Next, we shall check that the frequency map of the triaxial ellipsoid $Q$
can be continuously extended to the borders of the caustic space in
such a way that its values on the edges and vertexes can be expressed
in terms of exactly six functions of one variable that ``glue'' well.
Three of them are the extended rotation functions
associated to the three ellipses obtained by sectioning $Q$ with
the coordinate planes $\pi_\rmx$, $\pi_\rmy$, and $\pi_\rmz$.
That is,
they are the functions
$\rho_\rmx : [0,b] \to \Rset$, $\rho_\rmy : [0,a] \to \Rset$,
and $\rho_\rmz : [0,a] \to \Rset$ defined as
\[
\rho_\rmx(\lambda) = \rho(\lambda;c,b),\qquad
\rho_\rmy(\lambda) = \rho(\lambda;c,a), \qquad
\rho_\rmz(\lambda) = \rho(\lambda;b,a),
\]
using the notation in~(\ref{eq:RotationNumber}).
The other three functions are defined in terms of the former ones as follows.
Let $\underline{m}=\min(\lambda,c)$ and $\overline{m}=\max(\lambda,c)$.
Let $T_\rmx(s) = (\lambda-s)(c-s)(b-s)$,
$T_\rmy(s) = (c-s)(\lambda-s)(a-s)$,
and $T_\rmz(s) = (\underline{m}-s)(b-s)(a-s)$.
Then we consider the functions $\nu_\rmx : [0,b] \to \Rset$, $\nu_\rmy: [b,a] \to \Rset$,
and $\nu_\rmz: [0,b] \to \Rset$ defined by the identities
\begin{eqnarray*}
\int_0^{\underline{m}} \frac{\rmd s}{(a-s)\sqrt{T_\rmx(s)}} & -
2\rho_\rmx(\lambda) \int_{\overline{m}}^b \frac{\rmd s}{(a-s)\sqrt{T_\rmx(s)}} & +
\frac{2\pi \nu_\rmx(\lambda)}{\sqrt{-T_\rmx(a)}} = 0, \\
\int_0^c \frac{\rmd s}{(b-s)\sqrt{T_\rmy(s)}} & +
2\rho_\rmy(\lambda)  \int_{\lambda}^a \frac{\rmd s}{(s-b)\sqrt{T_\rmy(s)}} & -
\frac{2\pi \nu_\rmy(\lambda)}{\sqrt{-T_\rmy(b)}} = 0, \\
\int_0^{\underline{m}} \frac{\rmd s}{(\overline{m}-s)\sqrt{T_\rmz(s)}} & +
2\rho_\rmz(\underline{m}) \int_b^a \frac{\rmd s}{(s-\overline{m})\sqrt{T_\rmz(s)}} & -
\frac{2\pi \nu_\rmz(\lambda)}{\sqrt{-T_\rmz(\overline{m})}} = 0.
\end{eqnarray*}

\begin{lem}\label{lem:Nu}
The functions $\nu_\rmx$, $\nu_\rmy$, and $\nu_\rmz$ have the following properties.
\begin{enumerate}
\item
They are analytic in $E \cup H_1$, $H_2$, and $E\cup H_1$, respectively.
\item
They can be continuously extended to $[0,b]$, $[b,a]$, and $[0,b]$,
respectively.
\item
Their asymptotic behavior at the endpoints
$\lambda \in \partial E \cup \partial H_1 \cup \partial H_2 = \{0,c,b,a\}$ are:
\begin{enumerate}
\item
$\nu_\rmx(\lambda) = \Or(\lambda^{1/2})$, as $\lambda \to 0^+$;
\item
$\nu_\rmz(\lambda) = \Or(\lambda^{1/2})$, as $\lambda \to 0^+$;
\item
$\nu_\rmx(\lambda) = \rho_\rmz(a) + \Or\left(1/\log|c-\lambda|\right)$,
as $\lambda \to c$;
\item
$\nu_\rmz(\lambda) = 1/2 + \Or\left(|\lambda-c|^{1/2}\right)$,
as $\lambda \to c$;
\item
$\nu_\rmx(\lambda) = \rho_\rmy(a) + \Or(b-\lambda)$, as $\lambda \to b^-$;
\item
$\nu_\rmz(\lambda) = \rho_\rmz(c) + \Or\left((b-\lambda)^{1/2}\right)$,
as $\lambda \to b^-$;
\item
$\nu_\rmy(\lambda) = \rho_\rmz(c) + \Or\left((\lambda-b)^{1/2}\right)$,
as $\lambda\to b^+$; and
\item
$\nu_\rmy(\lambda) = \rho_\rmx(b) + \Or(a-\lambda)$, as $\lambda\to a^-$.
\end{enumerate}

\end{enumerate}
\end{lem}

\proof
We know that the function $\rho_\rmx(\lambda) = \rho(\lambda;c,b)$ is analytic
in $\lambda$, $c$, and $b$, as long as $0 < c < b$ and $\lambda \in E \cup H_1$.
Besides, the integrand $(a-s)^{-1}(T_\rmx(s))^{-1/2}$ is analytic with respect
to the variable of integration $s$ in the intervals of integration
$(0,\underline{m})$ and $(\overline{m},b)$,
and with respect to the parameters $\lambda$, $c$, $b$, and $a$,
as long as $0 < c < b < a$ and $\lambda \in E \cup H_1$.
Hence, the function $\nu_\rmx(\lambda)=\nu_\rmx(\lambda;c,b,a)$ is analytic in
its four variables, as long as $0 < c < b < a$ and $\lambda \in E \cup H_1$.
The analyticity of $\nu_\rmy$ and $\nu_\rmz$ follows
from similar arguments.

The study of the asymptotic behavior of the functions $\nu_\rmx$, $\nu_\rmy$,
and $\nu_\rmz$ has been deferred to~\ref{sap:Nux}, \ref{sap:Nuy},
and~\ref{sap:Nuz}, respectively.
\qed

\begin{remark}\label{rem:NuMonotone}
We have numerically observed that $\nu_\rmx$ and $\nu_\rmz$ are
increasing in $E$ and decreasing in $H_1$,
whereas $\nu_\rmy$ is increasing in $H_2$,
but we have not been able to prove it.
\end{remark}

\begin{thm}\label{thm:3D}
The frequency map $\omega : \Lambda \to \Rset^2$ has the
following properties.
\begin{enumerate}
\item
It is analytic in $\Lambda$.
\item
It can be continuously extended to the border $\partial \Lambda$,
and the extension has the form
\begin{eqnarray*}
\omega(0,\lambda_2) & = (0,0) \quad & \mbox{for } c \le \lambda_2 \le b,  \\
\omega(\lambda_1,b) & = ( \rho_\rmy(\lambda_1) , \rho_\rmy(\lambda_1) ) \quad
&\mbox{for } 0 \le \lambda_1 \le b, \\
\omega(c,\lambda_2) & = ( 1/2, \rho_\rmz(\lambda_2) ) \quad
&\mbox{for } c \le \lambda_2 \le a, \\
\omega(\lambda_1,a) & = ( \rho_\rmx(\lambda_1) , \nu_\rmx(\lambda_1) ) \quad
&\mbox{for } 0 \le \lambda_1 \le b, \\
\omega(b,\lambda_2) & = ( \nu_\rmy(\lambda_2) , \rho_\rmy(\lambda_2) ) \quad
&\mbox{for } b \le \lambda_2 \le a, \\
\omega(\lambda_1,c) & = ( \nu_\rmz(\lambda_1) , \rho_\rmz(\lambda_1) ) \quad
&\mbox{for } 0 \le \lambda_1 \le c, \\
\omega(\lambda_1,\lambda_1) & = ( \nu_\rmz(\lambda_1) , \rho_\rmz(c) ) \quad
&\mbox{for } c \le \lambda_1 \le b.
\end{eqnarray*}
\item
Its asymptotic behavior at the eleven edges in
$\Lambda^1 = G \cup S \cup R$ is:
\begin{enumerate}
\item
$\omega(\lambda_1,\lambda_2) =
\kappa^G(\lambda_2) \lambda_1^{1/2} + \Or(\lambda_1^{3/2})$,
as $\lambda_1 \to 0^+$;
\item
$\omega(\lambda_1,\lambda_2) - \omega(c,\lambda_2) \asymp
\kappa^S(c,\lambda_2) /\log |c-\lambda_1|$,
as $\lambda_1 \to c$;
\item
$\omega(\lambda_1,\lambda_2) - \omega(\lambda_1,b) \asymp
\kappa^S(\lambda_1,b) /\log |b-\lambda_2|$,
as $\lambda_2 \to b$; and
\item
$\omega(\lambda) - \omega(\lambda^R) = \Or(\lambda-\lambda^R)$,
as $\lambda \to \lambda^R \in R$;
\end{enumerate}
for some analytic functions $\kappa^G : H_1 \cup H_2 \to \Rset_+^2$ and
$\kappa^S: S \to \Rset^2$.
\item
Its asymptotic behavior at the eight vertexes in $\Lambda^0$ is:
\begin{enumerate}
\item
$\omega(\lambda_1,\lambda_2) = \Or(\lambda_1^{1/2})$,
as $(\lambda_1,\lambda_2) \to (0^+,c^+)$;
\item
$\omega(\lambda_1,\lambda_2) =
 \Or(\lambda_1^{1/2})$,
as $(\lambda_1,\lambda_2) \to (0^+,b)$;
\item
$\omega(\lambda_1,\lambda_2) = \Or(\lambda_1^{1/2})$,
as $(\lambda_1,\lambda_2) \to (0^+,a^-)$;
\item
$\omega(\lambda_1,\lambda_2) =
 (1/2,\rho_\rmz(c)) + \Or(1/\log |c-\lambda_1|,\lambda_2-c)$,
as $(\lambda_1,\lambda_2) \to (c,c^+)$;
\item
$\omega(\lambda_1,\lambda_2) =
 (1/2,1/2) + \Or(1/\log |c-\lambda_1|,1/\log |b-\lambda_2|)$,
as $(\lambda_1,\lambda_2) \to (c,b)$;
\item
$\omega(\lambda_1,\lambda_2) =
 (1/2,\rho_\rmz(a)) + \Or(1/\log |c-\lambda_1|,a-\lambda_2)$,
as $(\lambda_1,\lambda_2) \to (c,a^-)$;
\item
$\omega(\lambda_1,\lambda_2) =
 (\rho_\rmy(b),\rho_\rmy(b)) + \Or(b-\lambda_1,1/\log |b-\lambda_2|)$,
as $(\lambda_1,\lambda_2) \to (b^-,b)$; and
\item
$\omega(\lambda_1,\lambda_2) =
 (\rho_\rmx(b),\rho_\rmy(a)) + \Or(b-\lambda_1,a-\lambda_2)$,
as $(\lambda_1,\lambda_2) \to (b^-, a^-)$.
\end{enumerate}
\end{enumerate}
\end{thm}

\begin{figure}
\begin{center}
\begin{picture}(400,120)(0,0)
\footnotesize

\put(0,97){EH1:}
\thinlines
\put(50,100){\line(1,0){350}}
\linethickness{2pt}
\put(90,100){\line(1,0){10}}
\put(160,100){\line(1,0){80}}
\put(270,100){\line(1,0){80}}
\put(90,100){\circle*{5}}
\put(160,100){\circle*{5}}
\put(270,100){\circle*{5}}
\put(350,100){\circle*{5}}
\put(88,90){$0$}
\put(158,90){$c$}
\put(268,90){$b$}
\put(348,90){$a$}
\put(86,106){$c_0$}
\put(156,106){$c_2$}
\put(266,106){$c_4$}
\put(346,106){$c_5$}
\put(100,100){\circle*{5}}
\put(240,100){\circle*{5}}
\put(96,90){$\lambda_1$}
\put(236,90){$\lambda_2$}
\put(96,106){$c_1$}
\put(236,106){$c_3$}

\put(0,67){H1H1:}
\thinlines
\put(50,70){\line(1,0){350}}
\linethickness{2pt}
\put(90,70){\line(1,0){70}}
\put(190,70){\line(1,0){10}}
\put(270,70){\line(1,0){80}}
\put(90,70){\circle*{5}}
\put(160,70){\circle*{5}}
\put(270,70){\circle*{5}}
\put(350,70){\circle*{5}}
\put(88,60){$0$}
\put(158,60){$c$}
\put(268,60){$b$}
\put(348,60){$a$}
\put(86,76){$c_0$}
\put(156,76){$c_1$}
\put(266,76){$c_4$}
\put(346,76){$c_5$}
\put(190,70){\circle*{5}}
\put(200,70){\circle*{5}}
\put(186,60){$\lambda_1$}
\put(196,60){$\lambda_2$}
\put(186,76){$c_2$}
\put(196,76){$c_3$}

\put(0,37){EH2:}
\thinlines
\put(50,40){\line(1,0){350}}
\linethickness{2pt}
\put(90,40){\line(1,0){60}}
\put(160,40){\line(1,0){110}}
\put(280,40){\line(1,0){70}}
\put(90,40){\circle*{5}}
\put(160,40){\circle*{5}}
\put(270,40){\circle*{5}}
\put(350,40){\circle*{5}}
\put(88,30){$0$}
\put(158,30){$c$}
\put(268,30){$b$}
\put(348,30){$a$}
\put(86,46){$c_0$}
\put(156,46){$c_2$}
\put(266,46){$c_3$}
\put(346,46){$c_5$}
\put(150,40){\circle*{5}}
\put(280,40){\circle*{5}}
\put(146,30){$\lambda_1$}
\put(276,30){$\lambda_2$}
\put(146,46){$c_1$}
\put(276,46){$c_4$}

\put(0,7){H1H2:}
\thinlines
\put(50,10){\line(1,0){350}}
\linethickness{2pt}
\put(90,10){\line(1,0){70}}
\put(260,10){\line(1,0){10}}
\put(340,10){\line(1,0){10}}
\put(90,10){\circle*{5}}
\put(160,10){\circle*{5}}
\put(270,10){\circle*{5}}
\put(350,10){\circle*{5}}
\put(88,0){$0$}
\put(158,0){$c$}
\put(268,0){$b$}
\put(348,0){$a$}
\put(86,16){$c_0$}
\put(156,16){$c_1$}
\put(266,16){$c_3$}
\put(346,16){$c_5$}
\put(260,10){\circle*{5}}
\put(340,10){\circle*{5}}
\put(256,0){$\lambda_1$}
\put(336,0){$\lambda_2$}
\put(256,16){$c_2$}
\put(336,16){$c_4$}

\end{picture}
\end{center}
\caption{The four possible configurations of the ordered sequence
$0 < c_1 < \cdots < c_5$.
Thick lines denote intervals of integration.
Each one of the displayed configurations illustrates some collapse:
geodesic flow limit (and type EH1),
simple regular collapse (and type H1H1),
double singular collapse (and type EH2), and
double regular collapse (and type H1H2).}
\label{fig:C}
\end{figure}
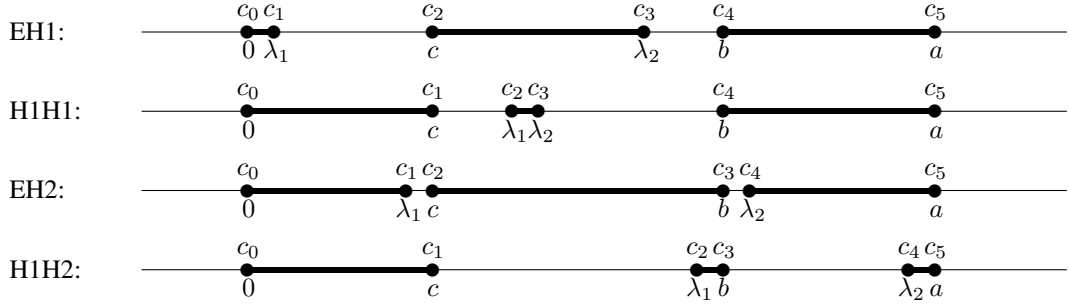

\proof
Once fixed the parameters $a > b > c > 0$ of the ellipsoid and the
couple of caustic parameters $\lambda_1$ and $\lambda_2$, we set
\[
\{ c_1,\ldots,c_5 \} = \{ a,b,c \} \cup \{ \lambda_1,\lambda_2\},
\qquad c_0 := 0 < c_1 < \cdots < c_5.
\]
Four configurations are possible; see figure~\ref{fig:C}.
We said in remark~\ref{remark:Homogeneous} that
the frequency is analytic in $c_1,\ldots,c_5$ provided that
$0 < c_1 < \cdots < c_5$.
In particular, this implies that the frequency is analytic in the
caustic parameter provided it belongs to $\Lambda$.

The frequency map is expressed in terms of six hyperelliptic integrals
over the intervals $(0,c_1)$, $(c_2,c_3)$, and $(c_4,c_5)$
---represented in thick lines in figure~\ref{fig:C}.
See definition~\ref{defi:FrequencyMap}.
We face its asymptotic behavior at the border
$\partial \Lambda = \Lambda^0 \cup \Lambda^1$,
which requires the study of the asymptotic behavior of the six
hyperelliptic integrals when some interval defined
by the ordered sequence $0 < c_1 < \cdots < c_5$ collapses to a point.
Therefore, there are exactly five simple collapses.
The collapse of the first interval is called
\emph{geodesic flow limit}: $c_1 \to 0^+$,
the collapse of the second or fourth intervals is called \emph{singular}:
$c_2 - c_1 \to 0^+$ or $c_4 - c_3 \to 0^+$, and
the collapse of the third or fifth intervals is called \emph{regular}:
$c_3 - c_2 \to 0^+$ or $c_5 - c_4 \to 0^+$.
Thus, regular collapses imply that the interval of integration of a
couple of hyperelliptic integrals collapses to a point;
whereas singular collapses imply the connection of two consecutive intervals of
integration.
See figure~\ref{fig:C}.
It is immediate to check that this terminology agrees with
the partition $\Lambda^1 = G \cup R \cup S$,
whereas double collapses ---that is, two simultaneous simple collapses---
correspond to the eight vertexes in $\Lambda^0$.

The asymptotic behavior of the frequency map at the eleven edges in
$\Lambda^1 = G \cup R \cup S$ is deduced from several results
disseminated through~\ref{ap:HyperellipticIntegrals}.
In short, some technical lemmas are listed in~\ref{sap:TechnicalLemmas},
some notations are introduced in~\ref{sap:Frequency},
the geodesic flow limit is studied in~\ref{sap:GeodesicFlowLimit},
simple regular collapses are analyzed in~\ref{sap:SimpleRegularCollapse},
and simple singular collapses are computed in~\ref{sap:SimpleSingularCollapse}.
For instance, one can trace the definition of the functions $\nu_\rmx$,
$\nu_\rmy$, and $\nu_\rmz$ to equation~(\ref{eq:Rcondition}).
The reader is encouraged to consult the appendix.
Here, we just note that the appendix deals with the general
high-dimensional setup, since the computations do not become
substantially more involved when the dimension is increased.

The computations regarding the vertexes
have also been relegated to~\ref{ap:HyperellipticIntegrals},
although for the sake of brevity we have written out only
the computations for two vertexes.
Vertex $\lambda=(c,b)$ in~\ref{sap:TotalSingularCollapse}
---which corresponds to the unique double singular collapse---,
and vertex $\lambda=(b,a)$ in~\ref{sap:TotalRegularCollapse}
---which correspond to the unique double regular collapse.
The study of the remaining six vertexes does not require
additional ideas.
For instance, the three vertexes related to the geodesic flow limit
can be simultaneously dealt with simply by using lemma~\ref{lem:G},
which ensures that the hyperelliptic integrals over
$(c_0,c_1)=(0,\lambda_1)$ are $\Or(\lambda_1^{1/2})$ as $\lambda_1 \to 0^+$.

Finally, we realize that the extended frequency map
$\omega:\bar{\Lambda} \to \Rset^2$ is continuous because the extensions
``glue'' well at the eight vertexes; see lemma~\ref{lem:Nu}.
For instance, let us consider the vertex $(b,b)$.
We obtain from the three statements of the theorem regarding this vertex that
\[
\omega(b,b) =
(\rho_\rmy(b),\rho_\rmy(b)) =
(\nu_\rmy(b),\rho_\rmy(b)) =
(\nu_\rmz(b),\rho_\rmz(c)),
\]
which is consistent:
$\nu_\rmy(b) = \nu_\rmz(b) =\rho_\rmz(c) =
 \rho(c;b,a) = \rho(b;c,a) = \rho_\rmy(b)$.
\qed

\begin{defi}\label{defi:Extended_omega}
The continuous extension $\omega:\bar{\Lambda} \to \Rset^2$ is called
the \emph{(extended) frequency map} of the ellipsoid $Q$.
\end{defi}

The origin of the terminology ``geodesic flow limit'' can be explained
as follows.
The phase space of the geodesic flow on an triaxial ellipsoid
$Q \subset \Rset^3$ was completely described by Kn\"orrer~\cite{Knorrer1980}.
Any nonsingular geodesic on $Q$ oscillates between two symmetric
curvature lines obtained by intersecting $Q$ with some hyperboloid
$Q_\lambda$, $\lambda \in H_1 \cup H_2$.
The rotation number of those oscillations is the quotient
\[
\rho^G(\lambda) =
\frac{\int_c^{\min(b,\lambda)} \frac{s \rmd s}{\sqrt{T^G(s)}}}
     {\int_{\max(b,\lambda)}^a \frac{s \rmd s}{\sqrt{T^G(s)}}},
\qquad T^G(s) = -s(\lambda-s)(c-s)(b-s)(a-s),
\]
see~\cite[\S 4.1]{DragovicRadnovic2009}.
This rotation number $\rho^G(\lambda)$ can be continuously
extended to the closed interval $[c,a]$ with $\rho^G(b)=1$.
On the other hand,
the geodesic flow on the ellipsoid $Q$ with caustic lines
$Q \cap Q_{\lambda_2}$ can be obtained as a limit of the
billiard dynamics inside $Q$ when its first caustic
$Q_{\lambda_1}$ approaches $Q$; that is,
when $\lambda_1 \to 0^+$, so that $(\lambda_1,\lambda_2) \to G$.
Therefore, it is natural to look for a relation between the function
$\kappa^G = (\kappa^G_1, \kappa^G_2) : H_1 \cup H_2 \to \Rset_+^2$ and
the rotation number $\rho^G : H_1 \cup H_2 \to \Rset_+$.

\begin{lem}
$\rho^G = \kappa^G_2/\kappa^G_1$.
Thus,
$\omega_2(\lambda_1,\lambda_2)/\omega_1(\lambda_1,\lambda_2) =
 \rho^G(\lambda_2) + \Or(\lambda_1)$, as $\lambda_1 \to 0^+$.
\end{lem}

\proof
In~\ref{ap:G} we will check that $\kappa^G$ is the unique solution
of the linear system
\[
2
\left(
\begin{array}{rr} K^G_{01} & -K^G_{02} \\ K^G_{11} & -K^G_{12} \end{array}
\right)
\left(
\begin{array}{c} \kappa^G_1 \\ \kappa^G_2 \end{array}
\right) =
\left(
\begin{array}{c} K^G_{00} \\ 0 \end{array}
\right),
\]
where $K^G_{ij} = \int_{c_{2j}}^{c_{2j+1}} (T^G(s))^{-1/2} s^i \rmd s$,
$K^G_{00} = 2 (abc\lambda)^{-1/2}$, and
$\{c_2,c_3,c_4,c_5\} = \{a,b,c,\lambda\}$ with $c_2 < c_3 < c_4 < c_5$.
Therefore, since $\lambda \in H_1 \cup H_2$, it turns out that
$c_2=c$, $c_3=\min(b,\lambda)$, $c_4=\max(b,\lambda)$, and $c_5=a$.
Finally, $\kappa^G_2/\kappa^G_1 = K^G_{11}/K^G_{12} = \rho^G$.
\qed

\subsection{On the Jacobian of the frequency map}

\begin{figure}
\iffiguresPDF
\includegraphics[width=6.18in,height=6.18in]{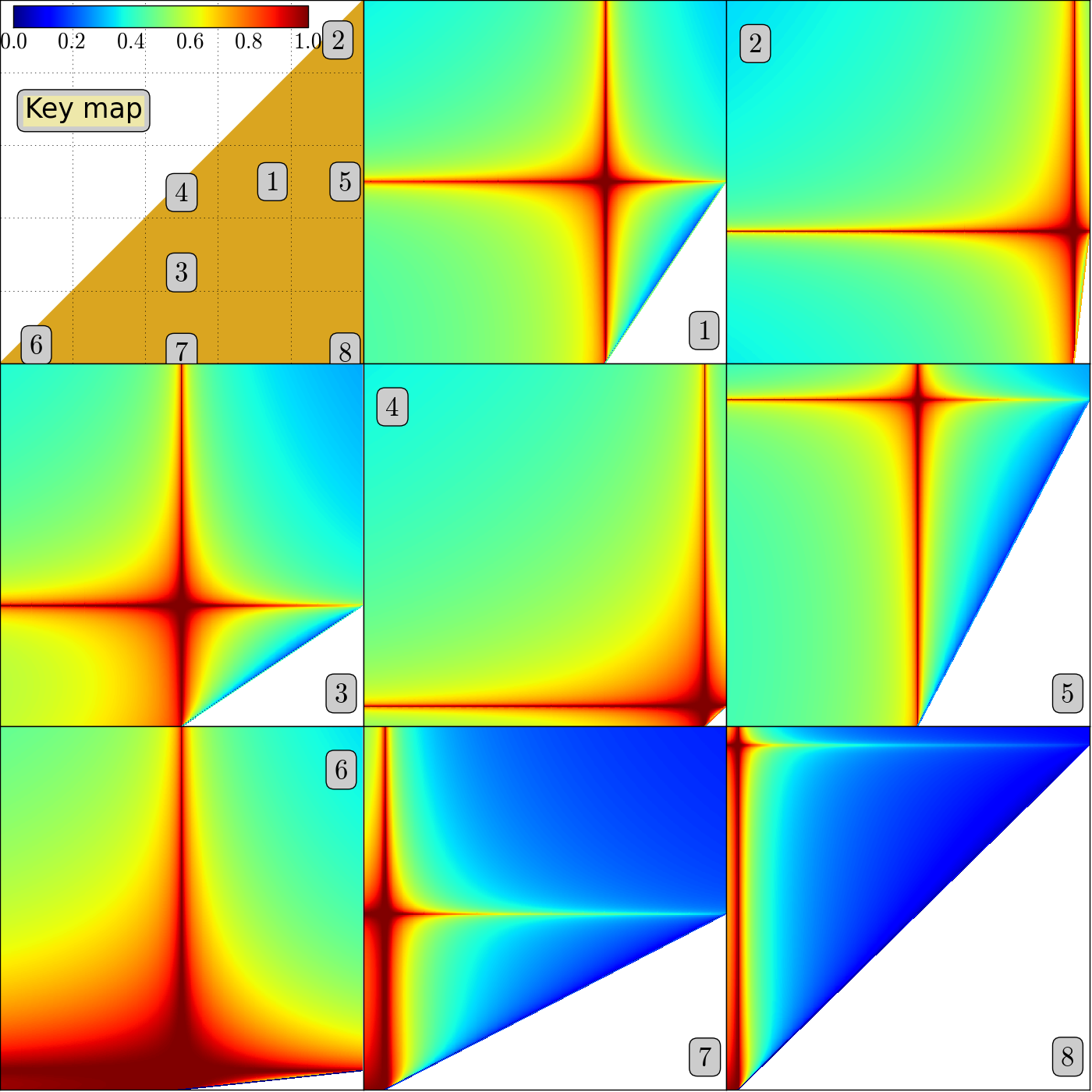}
\else
\vspace{6.26in}
\fi
\caption{The normalized Jacobian $J_\ast:\Lambda \to [0,1]$
of the frequency map for eight different ellipsoids.}
\label{fig:Jacobians}
\end{figure}

We present the numerical experiments about conjecture~\ref{con:FrequencyLocal}
stated in section~\ref{sec:FrequencyMap}.
We have computed the Jacobian of the frequency map
\[
J:\Lambda \to \Rset,\qquad
J(\lambda) :=
\det \left(
\frac{\partial \omega_j}{\partial \lambda_i}(\lambda)
     \right)_{i,j = 1,2}
\]
for several ellipsoids, in order to check that it never vanishes.
Its visualization close to the four inner edges labelled with
the letter $S$ in figure~\ref{fig:CausticSpace} has a technical difficulty.
To understand this fact,
one can look at the graph of the rotation number $\rho(\lambda)$ shown
in figure~\ref{fig:RotationFunction}.
The derivative $\rho'(\lambda)$ explodes at $\lambda=b$,
which would make difficult its visual representation.
The problem is worse in the spatial case,
because the frequency map has the same kind of ``inverse logarithm''
singularity at the four inner edges instead of at a single point.

We overcome the visualization problem by representing the normalized
Jacobian
\[
J_*:\Lambda \to [0,1],\qquad
J_*(\lambda) = (1-\exp(-|J(\lambda)|))^{1/4}.
\]
The exponential function is intended to cancel the exponentially
sharp behavior of the Jacobian at the inner edges.
The exponent $1/4$ has been chosen by trial and error to obtain more informative plots.
The normalized Jacobian ranges over the interval $[0,1]$.
We note that $J_\ast = 0 \Leftrightarrow J = 0$ and
$J_\ast = 1 \Leftrightarrow |J| = \infty$.
The results are shown in figure~\ref{fig:Jacobians}.
In the upper left corner,
we have displayed the parameter space $P$ introduced
in~(\ref{eq:ParameterSpace}), and sketched in
figure~\ref{fig:ParameterSpace}.
We study the eight ellipsoids that correspond to the eight points
in $P$ labelled from 1) to 8).
In particular, we have chosen at least one sample of each ``kind''
of ellipsoid:
1) standard, 2) almost spheric, 3) standard, 4) almost prolate,
5) almost oblate, 6) close to a segment, 7) close to a flat solid ellipse,
and 8) close to a flat circle.
The color palette is a classical one:
cold colors represent low values, hot colors represent high values.
The neighborhood of the inner edges is always a ``hot'' region;
that is, the Jacobian is always big on that region.
On the contrary, the Jacobian tends to zero close to the hypotenuse of
the $H_1\otimes H_1$ region. This can be seen from a symmetry reasoning.
Furthermore,
the Jacobian never vanishes, not even in the cases~7) and~8),
which correspond to almost flat ellipsoids.

\subsection{The range of the frequency map}\label{ssec:Range3D}

\begin{figure}
\iffiguresPDF
\begin{center}
\includegraphics[width=6.18in]{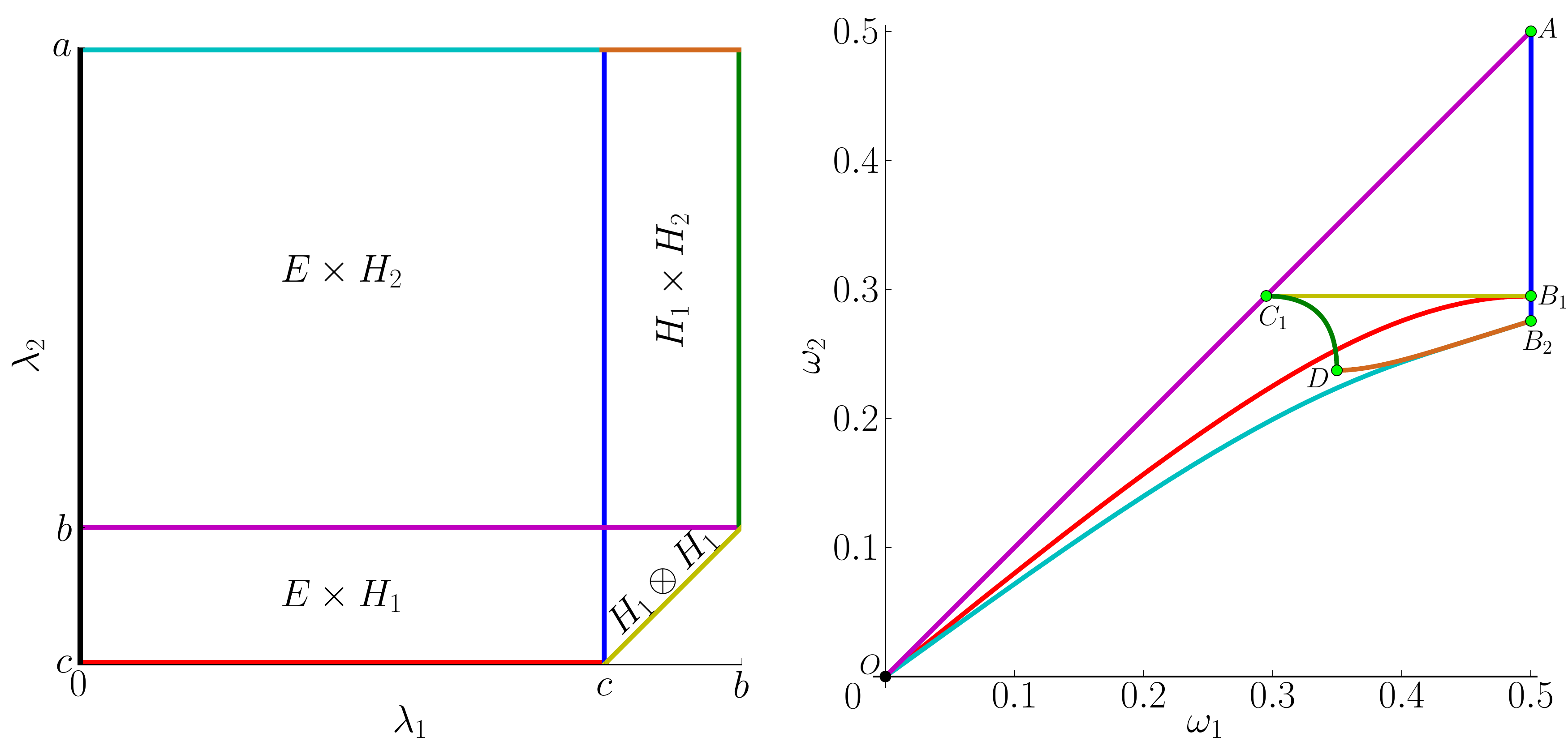}
\end{center}
\else
\fi
\caption{The extended frequency map $\omega:\bar{\Lambda} \to \bar{\Omega}$
\emph{on} the edges of the caustic space for $a=1$, $b=0.58$, and $c=0.46$.
Left: Caustic space. Right: Frequency space.}
\label{fig:Borders}
\end{figure}

We recall that if the two conjectures stated in
subsection~\ref{ssec:TwoConjectures} hold,
then the components of the frequency map are ordered as stated
in~(\ref{eq:FrequencyOrder}).
Thus, the range of the frequency map should be a subset of the frequency space
\[
\Omega =
\left\{
(\omega_1,\omega_2) \in \Rset^2 : 0 < \omega_2 < \omega_1 < 1/2
\right\}.
\]
We visualize in figure~\ref{fig:Borders} how each edge of the caustic space
is mapped onto the frequency space.
All the depicted curves have been numerically computed from exact formulae
given in theorem~\ref{thm:3D}.
We have represented the caustic space $\Lambda$ at the left side,
and the frequency space $\Omega$ at the right side.
Each colored segment in the caustic space is mapped onto the curve
of the same color in the frequency space.
The black segment in $\Lambda$ ---which represents the geodesic flow limit---
is mapped onto the origin $O=(0,0)$.
The point $(c,b)$ is mapped onto $A = (1/2,1/2)$.
The images of the magenta and blue segments are folded at this point $A$.
Henceforth,
$[AB]$ stands for the segment with endpoints $A$ and $B$,
and $\triangle[ABC]$ stands for the interior of the triangle
with vertexes $A$, $B$, $C$.
We see that $\omega(E \times H_1)$ is enclosed by the magenta segment $[OA]$,
the blue segment $[AB_1]$, and a red smooth curve from $B_1$ to $O$;
$\omega(H_1 \otimes H_1) = \triangle[AB_1C_1]$;
$\omega(E \times H_2)$ is enclosed by the magenta segment $[OA]$,
the blue segment $[AB_2]$, and a cyan smooth curve from $B_2$ to $O$; and
$\omega(H_1 \times H_2)$ is enclosed by the magenta segment $[C_1A]$,
the blue segment $[AB_2]$, a brown smooth curve from $B_2$ to $D$,
and a green smooth  curve from $D$ to $C_1$.

The points $A$, $B_1$, $B_2$, $C_1$, and $D$ can be
explicitly expressed in terms of the parameters of the ellipsoid $Q$.
Let $0 < \varrho_\rmx,\varrho_\rmy,\varrho_\rmz,\varrho_\ast < 1/2$
be the quantities defined by
\begin{equation}\label{eq:varrho}
\sin^2 \pi \varrho_\rmx = c/b, \quad
\sin^2 \pi \varrho_\rmy = c/a, \quad
\sin^2 \pi \varrho_\rmz = b/a, \quad
\varrho_\ast = \rho(c;b,a),
\end{equation}
where $\rho(\lambda;b,a)$ is the rotation number~(\ref{eq:RotationNumber}).
From the formulae contained in theorem~\ref{thm:3D},
we get that $A = (1/2,1/2)$, $B_1 = (1/2,\varrho_\ast)$,
$B_2 = (1/2,\varrho_\rmz)$, $C_1 = (\varrho_\ast,\varrho_\ast)$,
and $D=(\varrho_\rmx,\varrho_\rmy)$.
We note that $D \in \Omega$, since $\varrho_\rmy < \varrho_\rmx$.
In fact,
$\varrho_\rmy < \varrho_\rmz$ and $\varrho_\rmy < \varrho_\ast < \varrho_\rmx$,
although we do not have a rigorous proof of the inequalities involving $\varrho_\ast$.
The four quantities defined in~(\ref{eq:varrho})
can be interpreted in terms of the restriction of the billiard dynamics
to suitable planar sections of the original ellipsoid.
For instance,
$\varrho_\ast$ is the rotation number of the trajectories
contained in the section by the plane $\pi_\rmz$
whose caustic is the focal ellipse~(\ref{eq:FocalEllipse}).

We give now some numerical estimates on the
size and the shape of the four ranges.

\begin{num}\label{num:Ranges}
Let $0 < \varrho_\rmx,\varrho_\rmy,\varrho_\rmz,\varrho_\ast < 1/2$
be the quantities defined in~(\ref{eq:varrho}).
Let $\varrho = \max(\varrho_\ast,\varrho_\rmz)$.
Let $O=(0,0)$, $A = (1/2,1/2)$, $B_1 = (1/2,\varrho_\ast)$,
$B_2 = (1/2,\varrho_\rmz)$, $B=(1/2,\varrho)$,
$C_1 = (\varrho_\ast,\varrho_\ast)$, $C_2 = (\varrho_\rmz,\varrho_\rmz)$,
$C=(\varrho,\varrho)$, and $D=(\varrho_\rmx,\varrho_\rmy)$.
Then:
\begin{enumerate}
\item
$\triangle[A B_j C_j] \subsetneq
 \omega(E \times H_j) \subsetneq
 \triangle[AB_jO]$
for $j=1,2$;
\item
$\omega(H_1 \otimes H_1) = \triangle[AB_1C_1]$; and
\item
$\triangle[A B C] \subsetneq
 \omega(H_1 \times H_2) \subsetneq
 \omega(E \times H_2) \cap
 \{ (\omega_1,\omega_2) \in \Omega :
    \omega_1 > \varrho_\ast,\; \omega_2 > \varrho_\rmy \}$.
\end{enumerate}
\end{num}

Next, we enlighten some practical consequences of these estimates.
To begin with, let us present four simple criteria to decide
if the ellipsoid has billiard trajectories of frequency
$\omega^0 = (\omega_1^0, \omega_2^0) \in \Omega$ and
of caustic type EH1, H1H1, EH2, or H1H2.
Compare with the criterion for the existence of billiard trajectories
inside an ellipse with rotation number $\rho^0 \in (0,1/2)$ and
a caustic hyperbola given in~(\ref{eq:Bifurcations2D}).

\begin{pro}\label{pro:Criterions}
If numerical result~\ref{num:Ranges} holds,
then the following criteria can be applied.
\begin{enumerate}
\item
If $\omega^0_2 > \rho(c;b,a)$,
then $\omega^0 \in \omega(E \times H_1)$.
If $\omega^0_2/2\omega^0_1 \le \rho(c;b,a)$,
then $\omega^0 \not \in \omega(E \times H_1)$.
\item
$\omega^0 \in \omega(H_1 \otimes H_1)$ if and only if $\omega^0_2 > \rho(c;b,a)$.
\item
If $b < a \sin^2 \pi\omega^0_2$, then
$\omega^0 \in \omega(E \times H_2)$.
If $b \ge a \sin^2(\pi\omega^0_2/2\omega^0_1)$, then
$\omega^0 \not \in \omega(E \times H_2)$.
\item
If $\omega^0_2 > \rho(c;b,a)$ and $b < a \sin^2 \pi\omega^0_2$,
then $\omega^0 \in \omega(H_1 \times H_2)$.
If $\omega^0_1 \le \rho(c;b,a)$,
or $c \ge a \sin^2 \pi \omega^0_2$,
or $b \ge a \sin^2(\pi\omega^0_2/2\omega^0_1)$,
then $\omega^0 \not \in \omega(H_1 \times H_2)$.
\end{enumerate}
Hence, there exist billiard trajectories of the four caustic types
when $\omega^0_2$ is big enough:
$\omega^0_2 > \rho(c;b,a)$ and $\sin^2 \pi\omega^0_2 > b/a$.
On the contrary,
there does not exist any of such trajectories when
$\omega^0_2/\omega^0_1$ is small enough:
$\omega^0_2/\omega^0_1 \le 2 \rho(c;b,a)$ and
$\sin^2(\pi\omega^0_2/2\omega^0_1) \le b/a$.
\end{pro}

\proof
The first and third criteria follow from
$\triangle[A B_j C_j] \subset
 \omega(E \times H_j) \subset
 \triangle[AB_jO]$, $j=1,2$.
The second one follows from the identity
$\omega(H_1 \otimes H_1) = \triangle[AB_1C_1]$.
The last one follows from the the last item
of numerical result~\ref{num:Ranges}.
\qed

We can also understand how the range of the frequency map depends
on the shape of the ellipsoid.
It suffices to see how the quantities
$\varrho_\rmx$, $\varrho_\rmy$, $\varrho_\rmz$, and $\varrho_\ast$
depend on the parameters $0 < c < b < a$.
On the one hand,
if the ellipsoid flattens
---that is, if $c$ decreases, but $a$ and $b$ remain fixed---,
then $\varrho_\ast$ decreases,
so $\omega(E \times H_1)$ and $\omega(H_1 \otimes H_1)$ expand.
Indeed, both ranges tend to cover the whole space $\Omega$
for flat ellipsoids: $c \to 0^+$,
whereas they collapse to the empty set for prolate ellipsoids: $c \to b^-$.
On the other hand,
if the ellipsoid becomes more oblate
---that is, if $b$ increases, but $a$ and $c$ remain fixed---,
then $\varrho_\rmz$ increases,
so $\omega(E \times H_2)$ contracts.
Indeed, $\omega(E \times H_2)$ tends to cover $\Omega$ for
``segments'': $b\to 0^+$,
but collapses to the empty set for oblate ellipsoids: $b \to a^-$.
The behavior of $\omega(H_1 \times H_2)$ is more complicated,
because its vertex $D = (\varrho_\rmx,\varrho_\rmy)$ can be at
any point of the frequency space $\Omega$; see~(\ref{eq:varrho}).
Anyway, if the ellipsoid becomes spheric
---that is, $c$ and $b$ approach $a$---,
then $\varrho_\rmz$ and $\varrho_\ast$ tend to one half,
so $B_j$ tends to $A$ and the four ranges collapse to the empty set.
This means that \emph{the more spheric is an ellipsoid,
the poorer are its four types of nonsingular billiard dynamics.}
Some of the criteria stated in propositions~\ref{pro:Criterions} 
and~\ref{pro:Criterions2} quantify this general principle.

\begin{figure}
\iffiguresPDF
\includegraphics[width=6.18in,height=6.18in]{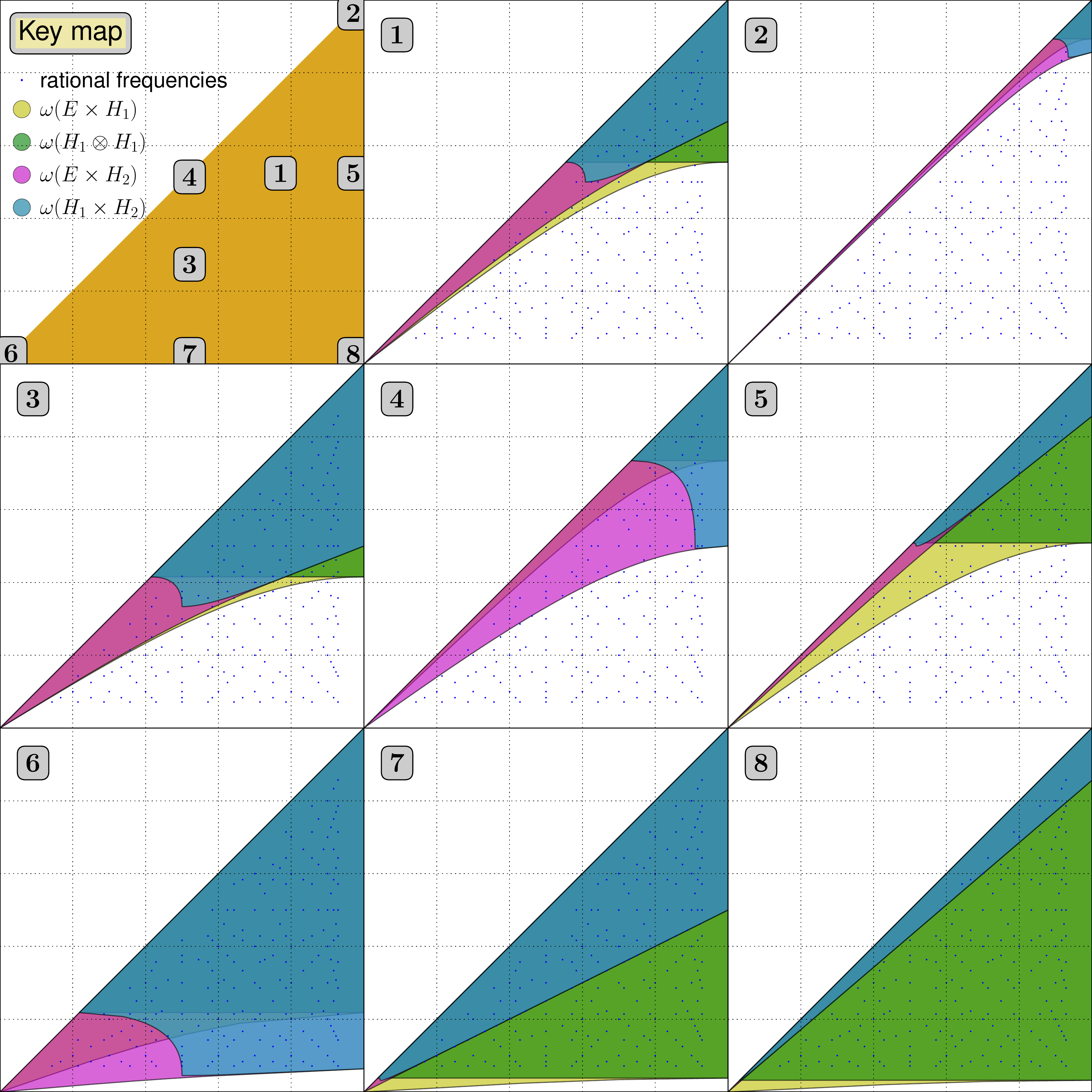}
\else
\vspace{6.26in}
\fi
\caption{Ranges of the frequency map for eight different ellipsoids.}
\label{fig:Ranges}
\end{figure}

The ranges of the frequency map for eight different ellipsoids
are shown in figure~\ref{fig:Ranges}.
In the left upper picture, we have again marked the
ellipsoids as points in the parameter space~(\ref{eq:ParameterSpace}).
The image sets $\omega(E \times H_1)$, $\omega(H_1 \otimes H_1)$,
$\omega(E \times H_2)$, and $\omega(H_1 \times H_2)$ are depicted
in yellow, green, magenta, and blue, respectively.
The transparency allows to visualize simultaneously all four sets.
We can check all their properties stated in numerical result~\ref{num:Ranges},
together with the ones regarding their dependence on
the shape of the ellipsoids.
Blue dots correspond to rational frequencies with small common
denominators.

In the previous paragraphs the range of the frequency map has been described
by mixing analytic formulae and numerical computations,
but some properties can be justified.
The following proposition is an example.

\begin{pro}\label{pro:Local2Global}
If conjecture~\ref{con:FrequencyLocal} holds, then
$\omega(H_1 \otimes H_1) = \triangle[AB_1C_1] \subsetneq \omega(E \times H_1)$,
and $\omega: H_1 \otimes H_1 \to \triangle[AB_1C_1]$ is
a global diffeomorphism.
Besides,
$\omega^0 \in \omega(H_1 \otimes H_1)$ if and only if
$\omega^0_2 > \varrho_\ast := \rho(c;b,a)$.
Finally, $\lim_{c \to 0^+} \triangle[AB_1C_1] = \Omega$ and
$\lim_{c \to b^-} \triangle[AB_1C_1] = \emptyset$.
\end{pro}

\proof
If $U = H_1 \otimes H_1$,
then $X = \partial U$ is the triangle with vertexes
$\tilde{A}=(c,b)$, $\tilde{B}_1=(c,c)$, $\tilde{C}_1=(b,b)$.
Using the formulae for the extended frequency map established
in theorem~\ref{thm:3D},
we get that $\omega([\tilde{A}\tilde{B}_1]) =[A B_1]$,
$\omega([\tilde{B}_1\tilde{C}_1]) =[B_1 C_1]$, and
$\omega([\tilde{C}_1\tilde{A}]) =[C_1 A]$.
Thus, $Y = \omega(X)$ is the triangle with vertexes $A$, $B_1$, $C_1$.
In particular, $X$ and $Y$ are Jordan curves,
so the frequency map $\omega:U \to \Rset^2$ verifies
the hypotheses of lemma~\ref{lem:LocalGlobal} in~\ref{ap:Topological}.
Hence, $\omega(U) = \triangle[AB_1C_1]$
and $\omega:U \to \omega(U)$ is a global diffeomorphism.

In order to prove the strict inclusion
$\triangle[AB_1C_1]\subsetneq \omega(E \times H_1)$,
it suffices to see that the red curve from $O$ to $B_1$ in the
right picture of figure~\ref{fig:Borders} is strictly below the
yellow segment $[B_1C_1]$.
And this is equivalent to prove the inequality
\[
\rho_\rmz(\lambda_1) < \rho_\rmz(c), \qquad \forall \lambda_1 \in (0,c),
\]
due to the formulae for the extended frequency map contained in
theorem~\ref{thm:3D}.
This inequality was proved in proposition~\ref{pro:Increasing}.
Finally, we note that
$\lim_{c \to 0^+} \varrho_\ast = 0$ and
$\lim_{c \to b^-} \varrho_\ast = 1/2$;
see the second item of proposition~\ref{pro:2D}.
\qed

\subsection{Geometric meaning of the frequency map}

Let $m_0,m_1,m_2$ be the winding numbers of a periodic
billiard trajectory of type EH1.
Then $m_0$ is the period.
Besides, according to remark~\ref{rem:WindingNumbers},
$m_1$ and $m_2/2$ are the number of times along one period that the trajectory
crosses the coordinate plane $\pi_\rmz  = \{ z=0 \}$ and
the number of times along one period that it rotates around
the coordinate axis $a_\rmz = \{ x = y = 0 \}$, respectively.
Therefore,
the components of the frequency map have the following geometric meaning:
$\omega_1 = m_1/2m_0$ is the number of oscillations around $\pi_\rmz$
per period,
whereas $\omega_2 = m_2/2m_0$ is the number of rotations around $a_\rmz$
per period.
Thus, it is quite natural to say that
$\omega_1$ is the \emph{$\rmz$-oscillation number} and
$\omega_2$ is the \emph{$\rmz$-rotation number} of the trajectory.

As in the planar case,
these interpretations are extended to quasiperiodic trajectories.
If $\lambda = (\lambda_1,\lambda_2) \in E \times H_1$,
then $Q_{\lambda_1}$ is an ellipsoid,
$Q_{\lambda_2}$ is a one-sheet hyperboloid,
and
\[
\omega(\lambda) = \lim_{k \to +\infty} (n_k,l_k)/k,
\]
where $n_k$ (respectively, $l_k$) is the number of oscillations
around $\pi_\rmz$ (respectively, number of rotations around $a_\rmz$)
of the first $k$ segments of \emph{any} given trajectory with caustics
$Q_{\lambda_1}$ and $Q_{\lambda_2}$.

\begin{table}
\begin{small}
\begin{center}
\begin{tabular}{c|c|c|c|c}
Type & $m_1$ & $m_2$ & $\omega_1$ & $\omega_2$ \\
\hline
EH1 & Crossings of $\pi_\rmz$ & Half-turns around $a_\rmz$ & $\rmz$-oscillation & $\rmz$-rotation \\
EH2 & Half-turns around $a_\rmx$ & Crossings of $\pi_\rmx$ & $\rmx$-rotation & $\rmx$-oscillation \\
H1H1 & Touches of $Q_{\lambda_j}$ & Half-turns around $a_\rmz$ & (H1-oscillation)/2 & $\rmz$-rotation \\
H1H2 & Crossings of $\pi_\rmy$ & Crossings of $\pi_\rmx$ & $\rmy$-oscillation & $\rmx$-oscillation
\end{tabular}
\end{center}
\end{small}
\caption{Geometric meaning of the winding numbers and the frequency vector
 when $Q \subset \Rset^3$.}
\label{tab:GeometricMeaning}
\end{table}

The billiard trajectories of other types can be analyzed
following similar arguments.
The results are listed in table~\ref{tab:GeometricMeaning}
and can be checked by visual inspection;
see figure~\ref{fig:PeriodicTrajectories3D}.

Finally, we stress a point already commented in remark~\ref{rem:Factor_2}.
If the trajectory is of type H1H1
---that is, if both caustics are one-sheet hyperboloids---,
then the winding number $m_1$ is the number of (alternate)
tangential touches with the caustics, so $\omega_1 = m_1/2m_0$ is half
the number of oscillations between the one-sheet hyperboloids per period.
In that situation,
we call $2\omega_1$ the \emph{H1-oscillation number} of the trajectory.
In particular, it can happen that $m_0 \omega \not \in \Zset^2$.
For instance, if the winding numbers are $m_0=4$, $m_1=3$, and $m_2=2$,
the period is four, but $\omega=(3/8,1/4)$.

\subsection{Bifurcations in parameter space}
\label{ssec:Bifurcations3D}

We want to determine the ellipsoids that have billiard trajectories with
a prescribed frequency and with a prescribed caustic type.
We recall that each ellipsoid is represented by a point in
$P = \{(b,c) \in \Rset^2 : 0 < c < b < 1\}$, because $a=1$.
Let $P^0_1$, $P^0_2$, $P^0_3$, and $P^0_4$
be the four regions of $P$ that correspond to ellipsoids with
billiard trajectories of frequency $\omega^0$ and caustic type
EH1, H1H1, EH2, and H1H2, respectively.
Their shapes are described below.

\begin{num}\label{num:BifurcationsG}
Once fixed any frequency vector
$\omega^0 = (\omega_1^0, \omega_2^0) \in \Omega$,
let $b^0_1 = b^0_2 = 1$, $b^0_3 = b^0_4 = \sin^2(\pi\omega^0_2/2\omega^0_1)$,
$c^0_1 = c^0_3 = \beta^0_4 = \sin^2 \pi \omega^0_2/\sin^2 \pi \omega^0_1$,
and $c^0_2 = c^0_4 = \sin^2 \pi \omega^0_2$.
Then
\[
P^0_j =
\left \{ (b,c) \in P : b <  b^0_j, \quad c <  g^0_j(b) \right \},
\quad  1 \le j \le 4,
\]
for some continuous functions $g^0_j: [0,b^0_j] \to \Rset$ such that
\begin{enumerate}
\item
$g^0_1$ is concave increasing in $[0,1]$, $0 < g^0_1(b) < b$ for all $b\in(0,1)$,
and $g^0_1(1) = c^0_1$;
\item
$g^0_2$ is concave increasing in $[0,1]$,
$0 < g^0_2(b) < g^0_1(b)$ for all $b\in(0,1)$,
and $g^0_2(1) = c^0_2$;
\item
$g^0_3$ is the identity in $[0,c^0_3]$,
concave decreasing in $[c^0_3,b^0_ 3]$,
and $g^0_3(b^0_ 3) = 0$; and
\item
$g^0_4$ is increasing in $[0,\beta^0_4]$,
concave decreasing in $[\beta^0_4,b^0_4]$,
$c^0_4 b/\beta^0_4 < g^0_4(b) < b$ for all $b\in(0,\beta^0_4)$,
$0 < g^0_4(b) < g^0_3(b)$ for all $b\in(\beta^0_4,b^0_4)$,
$g^0_4(\beta^0_4)=c^0_4$, and $g^0_4(b^0_4) = 0$.
\end{enumerate}
\end{num}

\begin{remark}
Numerical result~\ref{num:BifurcationsP} follows from
numerical result~\ref{num:BifurcationsG}
just by choosing suitable rational frequency vectors:
$\omega^0=(2/5,1/5)$ in the cases EH1 and EH2,
$\omega^0=(3/8,1/4)$ in the case H1H1,
and $\omega^0=(1/3,1/6)$ in the case H1H2.
We stress that
inequality $g^\ast_1 < g^\ast_2$ in numerical result~\ref{num:BifurcationsP}
and inequality $g^0_2 < g^0_1$ in numerical result~\ref{num:BifurcationsG}
are not contradictory,
because the first one refers to two different frequency vectors:
$(2/5,1/5)$ and $(3/8,1/4)$.
\end{remark}

\begin{remark}
We have numerically checked that $g^0_4$
is not concave in $[0,\beta^0_4]$.
\end{remark}

\begin{remark}
Inclusions $P^0_2 \subset P^0_1$ and $P^0_4 \subset P^0_3$
---and so, inequalities $g^0_2(b) < g^0_1(b)$ and $g^0_4(b) < g^0_3(b)$---
are in direct agreement with inclusions
$\omega(H_1 \otimes H_1) \subset \omega(E \times H_1)$ and
$\omega(H_1 \times H_2) \subset \omega(E \times H_2)$ mentioned in
numerical result~\ref{num:Ranges}.
\end{remark}

\begin{figure}
\iffiguresPDF
\begin{center}
\includegraphics[width=5.5in]{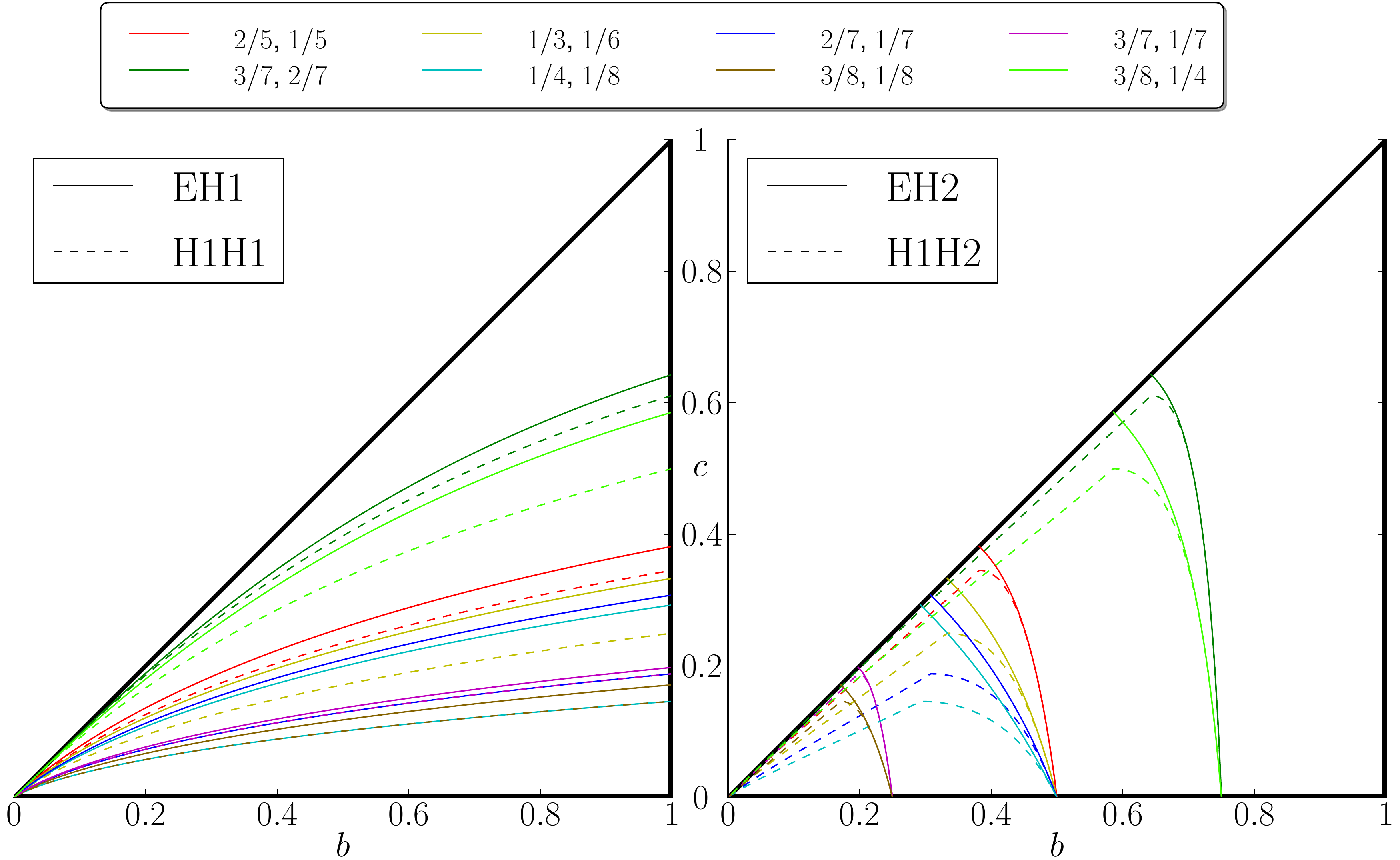}
\includegraphics[width=5.5in]{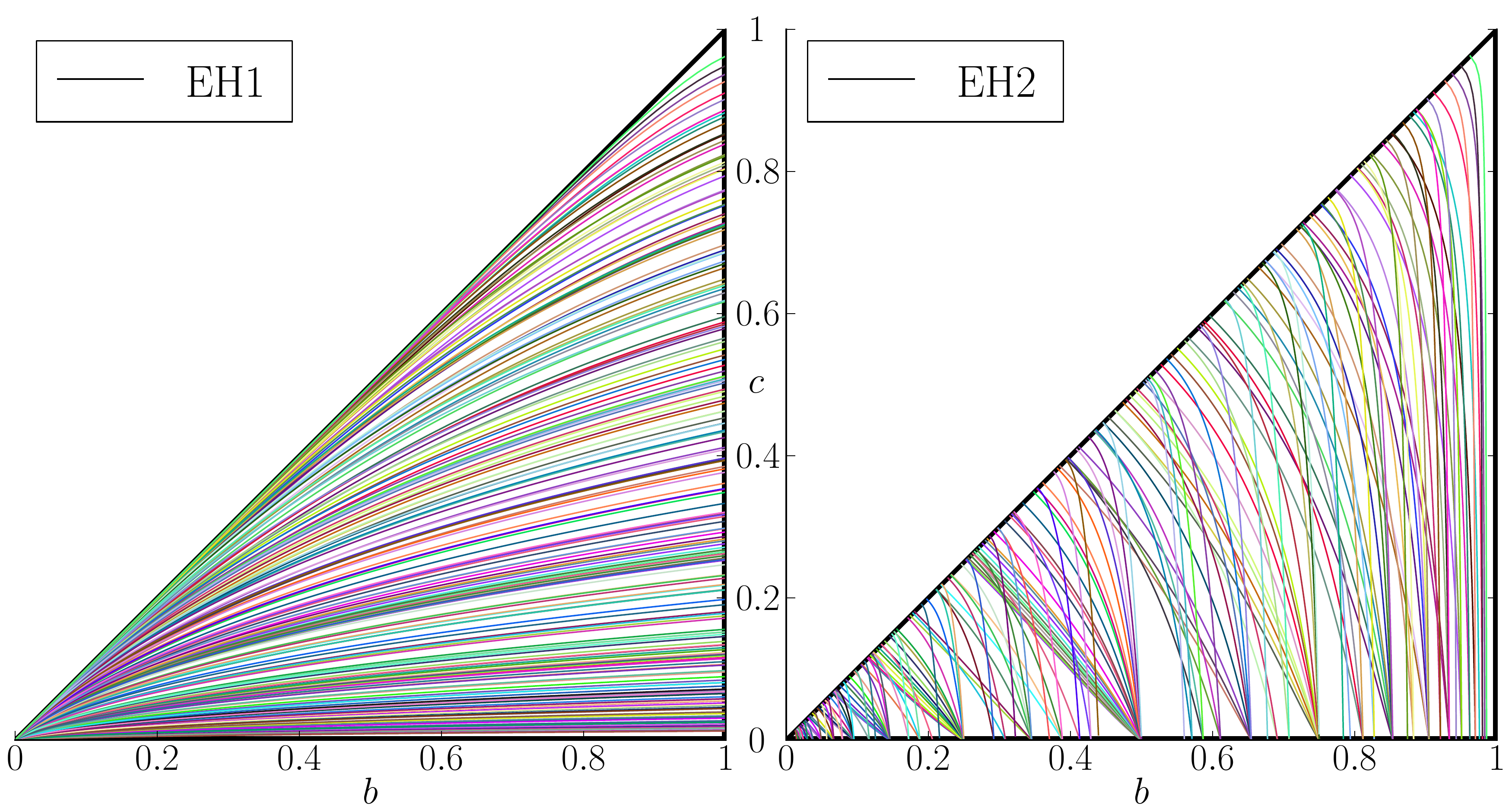}
\end{center}
\else
\vspace{6.4in}
\fi
\caption{Some bifurcation curves in the parameter space $P$.}
\label{fig:Bifurcations}
\end{figure}

Some bifurcations curves corresponding to the graphs of the functions
$g^0_j : [0,b^0_j] \to \Rset$ are presented in figure~\ref{fig:Bifurcations}.
On top of this figure we consider the eight rational frequencies
with the smallest denominators.
The inclusions $P^0_2 \subset P^0_1$ and $P^0_4 \subset P^0_3$ can be easily
visualized, since all dashed curves are below their continuous pairs.
On the bottom,
we depict the bifurcation curves associated to the rational frequencies
marked with blue dots in figure~\ref{fig:Ranges}.
We have needed a multiple precision arithmetic to compute the
bifurcation curves close to some of their endpoints,
since the involved root-finding problems become quite singular
at them.
The programs have been written using the PARI system~\cite{PARI}.

Next, we describe four more criteria to decide
if an ellipsoid has billiard trajectories of a given frequency.
They are similar to the four ones established in
proposition~\ref{pro:Criterions}.

\begin{pro}\label{pro:Criterions2}
If numerical result~\ref{num:BifurcationsG} holds,
then the following criteria can be applied.
\begin{enumerate}
\item
If $c < c^0_1 b$, then $\omega^0 \in \omega(E \times H_1)$.
If $c \ge c^0_1 a$, then $\omega^0 \not \in \omega(E \times H_1)$.
\item
If $c < c^0_2 b$, then $\omega^0 \in \omega(H_1 \otimes H_1)$.
If $c \ge c^0_2 a$, then $\omega^0 \not \in \omega(H_1 \otimes H_1)$.
\item
If $(b^0_3 - c^0_3)c <  c^0_3 (b^0_ 3 a - b)$,
then $\omega^0 \in \omega(E \times H_2)$.
If $b \ge b^0_3 a$, then $\omega^0 \not \in \omega(E \times H_2)$.
\item
If $\beta^0_4 c < \min \big(c^0_4 b, \beta^0_4 b^0_4 a + (c^0_4- b^0_4)b \big)$,
then $\omega^0 \in \omega(H_1 \times H_2)$.
If $b \ge b^0_4 a$ or $c \ge c^0_4 a$,
then $\omega^0 \not \in \omega(H_1 \times H_2)$.
\end{enumerate}
\end{pro}

\proof
From numerical result~\ref{num:BifurcationsG},
we get that $T^0_j := \triangle[O \Gamma^0_j \Delta^0_j] \subset P^0_j$,
where $O = (0,0)$, $\Gamma^0_1 = \Gamma^0_2 = (1,0)$,
$\Gamma^0_3 = (b^0_3,0)$, $\Gamma^0_4 = (b^0_4,0)$, $\Delta^0_1 = (1,c^0_1)$,
$\Delta^0_2 = (1,c^0_2)$, $\Delta^0_3 = (c^0_3,c^0_3)$,
and $\Delta^0_4 = (\beta^0_4,c^0_4)$.
It is straightforward to check that a point $(b,c) \in P$ belongs to the
triangles $T^0_1$, $T^0_2$, $T^0_3$, and $T^0_4$ if and only if
$c < c^0_1 b$, $c < c^0_2 b$, $(b^0_3 - c^0_3)c <  c^0_3 (b^0_ 3 - b)$, and
$\beta^0_4 c < \min \big(c^0_4 b, \beta^0_4 b^0_4 + (c^0_4- b^0_4)b \big)$,
respectively.
This proves the first part of each criterion for $a=1$.
To prove the general case, it suffices to take into account that
its formulae are homogeneous in the parameters $a,b,c$.

The second parts follow from similar arguments.
For instance, $P^0_j \cap \{ c \ge c^0_j\} = \emptyset$ for $j=1,2$,
since $g^0_j(b)$ are increasing in $[0,1]$ and $g^0_j(1) = c^0_j$.
\qed

\begin{remark}
Proposition~\ref{pro:Criterions} has been obtained by fixing
the ellipsoid and looking at the frequency space.
On the contrary,
proposition~\ref{pro:Criterions2} has been derived by fixing
the frequency vector and looking at the parameter space.
Of course, both approaches are equivalent,
but their criteria are slightly different.
The second ones are computationally simpler,
because they do not involve any elliptic integral.
\end{remark}

Although the description of the regions $P^0_j$
has a strong numerical component,
some results can be proved.
The following proposition is an example.

\begin{pro}\label{pro:BifurcationsH1H1}
Let $\omega^0 = (\omega_1^0, \omega_2^0) \in \Omega$
be a fixed frequency vector.
If conjecture~\ref{con:FrequencyLocal} holds,
then $P^0_2 \subset P^0_1$, $P^0_2$ only depends on $\omega^0_2$, and
$P^0_2 = \left \{ (b,c) \in P : 0 <  b <  1, \; c <  g^0_2(b) \right \}$
for some increasing analytic function $g^0_2: (0,1) \to \Rset$ such that
$0 < g^0_2(b) < b$ for all $b\in (0,1)$.
Besides, $\lim_{b \to 1^-} g^0_2(b) = \sin^2 \pi \omega^0_2$.
Finally,
\[
g^0_2(b) =
\cases{ 3b/\big(1+b+2\sqrt{1-b+b^2}\big) & \mbox{for $\omega^0_2 = 1/3$,} \\
        b/(1+b) & \mbox{for $\omega^0_2 = 1/4$.}}
\]
\end{pro}

\proof
If conjecture~\ref{con:FrequencyLocal} holds,
then $\omega(H_1 \otimes H_1) = \triangle(AB_1C_1)\subset \omega(E \times H_1)$;
see proposition~\ref{pro:Local2Global}.
Therefore, $P^0_2 \subset P^0_1$, because
\[
(b,c) \in P^0_2 \Rightarrow
\omega^0 \in \omega(H_1 \otimes H_1) \subset \omega(E \times H_1) \Rightarrow
(b,c) \in P^0_1.
\]
On the other hand,
since $A=(1/2,1/2)$, $B_1 =(1/2,\varrho_\ast)$ and $C_1=(\varrho_\ast,\varrho_\ast)$,
we deduce that
\[
(b,c) \in P^0_2 \Leftrightarrow
\omega^0 \in \omega(H_1 \otimes H_1) = \triangle(AB_1C_1)\Leftrightarrow
\omega^0_2 > \varrho_\ast := \rho(c;b,a).
\]
We know that the rotation function $\rho(\lambda) = \rho(\lambda;b,a)$
is increasing in $(0,b)$, $\rho(0) = 0$, and $\rho(b) = 1/2$.
Hence, the function $g^0_2:(0,1) \to \Rset$, $0 < g^0_2(b) <b$,
is implicitly defined by
\begin{equation}\label{eq:Implicit}
\rho(g^0_2(b);b,a) = \omega^0_2,
\end{equation}
where $\omega^0_2 \in (0,1/2)$ and $a=1$ are fixed parameters.
Analyticity of $g^0_2$ follows from the Implicit Function Theorem,
since conjecture~\ref{con:FrequencyLocal} also implies that
$\partial_1 \rho(\lambda;b,a) \neq 0$ for all $\lambda \in (0,b) \cup (b,a)$.
Indeed, this derivative is positive in $(0,b)$ and negative in $(b,a)$,
because $\rho(\cdot;b,a)$ is increasing in $(0,b)$ and decreasing in $(b,a)$.
Besides,
we know that $\partial_2 \rho(\lambda;b,a) = \partial_1 \rho(b;\lambda,a)$
from the symmetry $\rho(\lambda;b,a) = \rho(b;\lambda,a)$;
see~(\ref{eq:RotationNumber}).
Hence,
by differentiating equation~(\ref{eq:Implicit}) with respect to $b$
and setting $c= g^0_2(b) \in (0,b)$, we get that
\[
\big(g^0_2\big)'(b) =
- \partial_2 \rho(c;b,a) / \partial_1 \rho(c;b,a) =
- \partial_1 \rho(b;c,a) / \partial_1 \rho(c;b,a) > 0.
\]
Using proposition~\ref{pro:2D}, we know that
$\lim_{b \to a^-} \sin^2 \pi \rho(c;b,a) =
 \lim_{b \to a^-} \sin^2 \pi \rho(b;c,a) = c/a$.
Thus,
we deduce $\lim_{b \to 1^-} g^0_2(b) = \sin^2 \pi \omega^0_2$,
since $a=1$.

Finally,
we must find the values $c \in (0,b)$ such that $\rho(c;b,a)$
is equal to $1/3$ or $1/4$.
That is,
we must find the values of $c\in (0,b)$ such that the billiard
trajectories inside the ellipse $\{x^2/a + y^2/b = 1\}$
with caustic $\{ x^2/(a-c) + y^2/(b-c) = 1\}$ have period
three or four.
This is an old result that goes back to Cayley~\cite{Cayley1854}.
For instance, $c=ab/(a+b)$ in the four-periodic case.
The value for the three-periodic case was given in
equation~(\ref{eq:MinimalCaustics2D}).
\qed

The fact that $P^0_2$ only depends on $\omega^0_2$ can be visualized
on top left in figure~\ref{fig:Bifurcations}.
The two dashed curves with $\omega^0_2 = 1/8$ coincide,
as well as the two ones with $\omega^0_2 = 1/7$.

\subsection{On the ubiquity of almost singular trajectories}
\label{ssec:Singular3D}

\begin{figure}
\iffiguresPDF
\begin{center}
\includegraphics[width=6.18in]{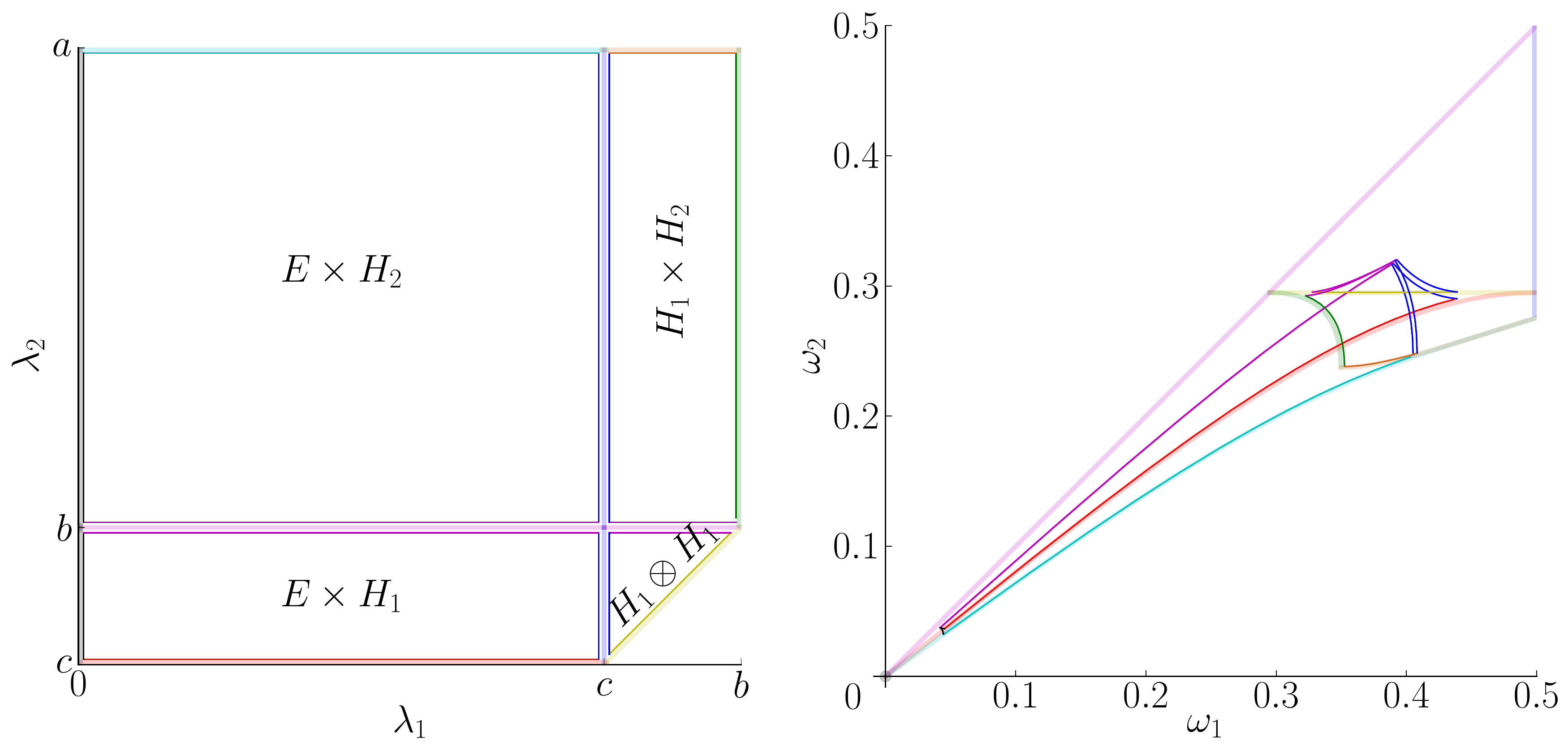}
\end{center}
\else
\fi
\caption{The extended frequency map $\omega:\bar{\Lambda} \to \bar{\Omega}$
\emph{close to} the edges of the caustic space for $a=1$, $b=0.58$,
and $c=0.46$. (Compare with figure~\ref{fig:Borders}.)}
\label{fig:ApproximatedBorders}
\end{figure}

In figure~\ref{fig:ApproximatedBorders} we have superposed
the edges and borders (drawn in light colors)
already displayed in figure~\ref{fig:Borders},
and some new segments and curves (drawn in heavy colors).
In the caustic space, these new segments are close to the original edges.
To be precise, the distance between them and the edges is equal to
$c/100 = 4.6 \cdot 10^{-3}$.
Nevertheless, the images of the black, magenta, and blue ones are
far from their corresponding borders in the frequency space.
This phenomenon seems stronger on the magenta and blue borders.
It has to do with the fact that, as stated in theorem~\ref{thm:3D},
the frequency map has an \emph{inverse logarithm} singularity at the
blue and magenta edges of the caustic space.
Therefore, one must be \emph{exponentially close} to them,
just to be close to their images.
On the other hand,
the frequency map has a \emph{squared root} singularity at the
black edges of the caustic space.
Thus, one must be \emph{quadratically close} to them,
just to be close to the origin in the frequency space.

We deduce from this phenomenon that
\emph{billiard trajectories with some almost singular caustic are ubiquitous.}
Let us describe a quantitative sample of this principle using
figure~\ref{fig:ApproximatedBorders}.
Let $T$ be the triangle delimited by the yellow, blue and magenta
thin segments that are close to the edges of $H_1 \otimes H_1$.
It turns out that the area of $\omega(H_1 \otimes H_1)$ is approximately
16 times the area of $\omega(T)$.
Hence, if we look for billiard trajectories of type H1H1 inside $Q$
with a random frequency in $\omega(H_1 \otimes H_1)$,
their caustic parameter $\lambda=(\lambda_1,\lambda_2)$ shall verify
$\min (|\lambda_1 -c|,|\lambda_2 - b|) < 4.6 \cdot 10^{-3}$
with a probability approximately equal to $94\%$.
It suffices to note that $15/16 = 0.9375$.

\subsection{Examples of periodic trajectories with minimal periods}
\label{ssec:PeriodicTrajectories3D}

\begin{figure}
\iffiguresPDF
\begin{center}
\includegraphics{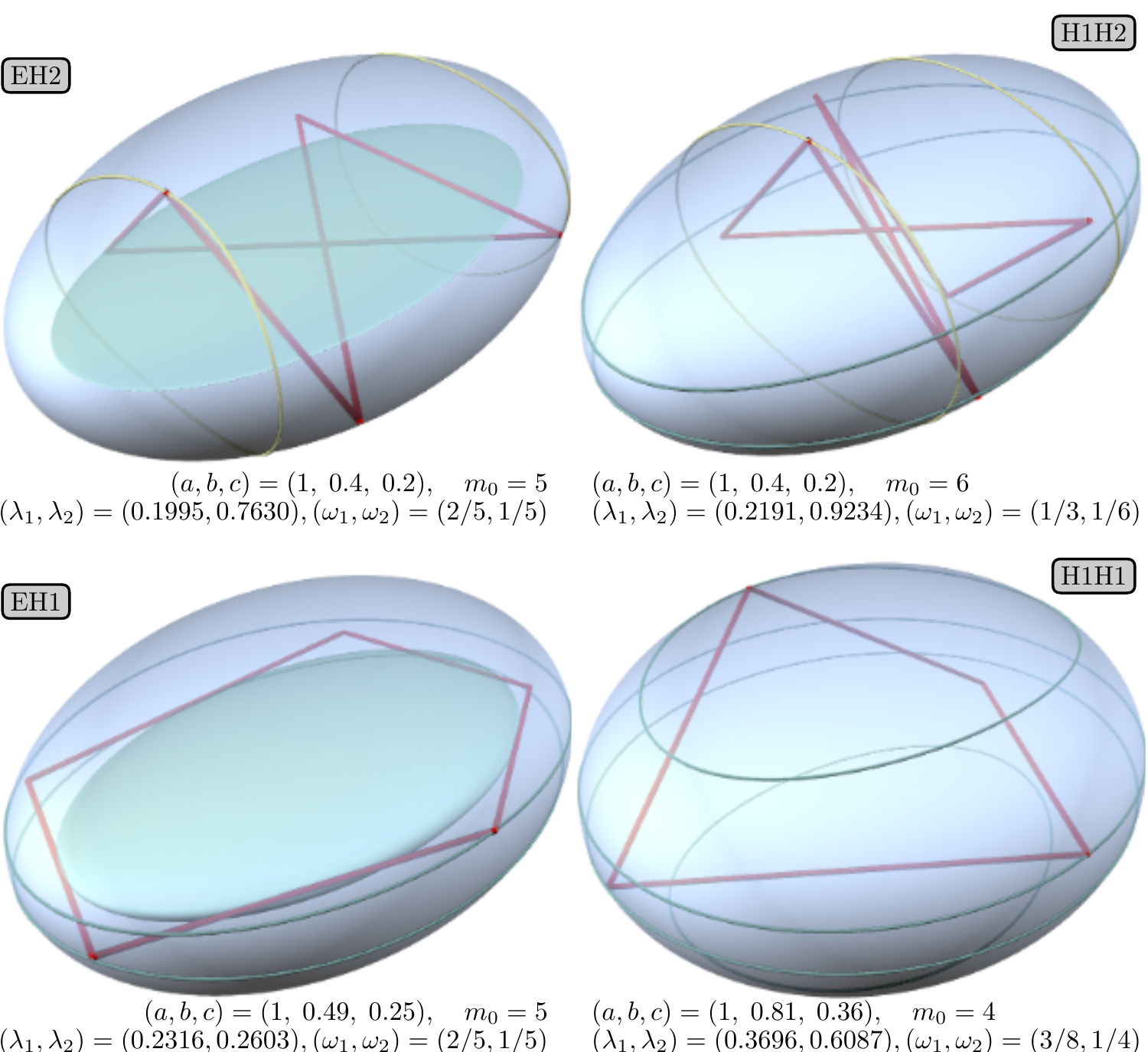}
\end{center}
\else
\vspace{6.4in}
\fi
\caption{Examples of symmetric nonsingular billiard trajectories with
minimal periods.
Lines in red represent the particle's trajectory.
Lines in green and yellow correspond to the intersections of the original
 ellipsoid with the caustic 1-sheet and 2-sheet hyperboloids, respectively.
In the cases EH1 and EH2, the caustic ellipsoid is also depicted.}
\label{fig:PeriodicTrajectories3D}
\end{figure}

We have numerically computed some symmetric periodic trajectories to check
that the lower bounds stated in theorem~\ref{thm:LowerBounds3D} are optimal;
see figure~\ref{fig:PeriodicTrajectories3D}.
All these trajectories are almost singular.
Concretely, $c - \lambda_1 \simeq 5 \cdot 10^{-4}$ in the case EH2;
$\lambda_1 - c \simeq 2 \cdot 10^{-2}$ in the case H1H2;
$c-\lambda_1 \simeq 2 \cdot 10^{-2}$ and
$\lambda_2 - c \simeq 10^{-2}$ in the case EH1;
and $\lambda_1 - c \simeq 10^{-2}$ in the case H1H1.
Of course, we did not look for almost singular trajectories,
but we got them anyway.

Considering the values given in figure~\ref{fig:PeriodicTrajectories3D},
and bearing in mind table~\ref{tab:GeometricMeaning},
we have $(m_1, m_2) = (4, 2)$ for the EH2 trajectory,
so it performs two turns around the coordinate axis $a_\rmx$ and
crosses twice the coordinate plane $\pi_\rmx$.
As well, $(m_1, m_2) = (4, 2)$ for the EH1 trajectory,
meaning four crossings with $\pi_\rmz$ and just one turn around $a_\rmz$.
Again, we have $(m_1, m_2) = (4, 2)$ for the H1H2 trajectory,
meaning four crossings with $\pi_\rmy$ and two crossings with $\pi_\rmx$.
Finally, $(m_1, m_2) = (3, 2)$ for the H1H1 trajectory,
which corresponds to three tangential touches with each of the
caustics and a single turn around $a_\rmz$.
Each of those geometric interpretations has been verified on
the corresponding trajectory.

\section{Billiard inside a nondegenerate ellipsoid of $\Rset^{n+1}$}
\label{sec:nD}

We describe briefly the high-dimensional version of some of
the analytical results already shown in the spatial case.
We denote again the nondegenerate ellipsoid as in~(\ref{eq:Q})
and the nonsingular caustic space as in~(\ref{eq:NonsingularSet}).

By analogy with the spatial case, we consider three disjoint partitions:
\[
\partial \Lambda = \cup_{k=0}^{n-1} \Lambda^k,\qquad
\Lambda^{n-1} = G \cup R \cup S,\qquad
S = \cup_{j=1}^n S_j.
\]
With regard to the first one,
$\Lambda^k$ is the $k$-dimensional border of $\Lambda$.
That is, $\Lambda^0$ is the set of vertexes,
$\Lambda^1$ is the set of edges, $\Lambda^2$ is the set of faces,
and so on.
The second one mimics the distinction among geodesic flow limits,
simple regular collapses, and simple singular collapses
already seen in the previous section.
For instance, $G=\{ \lambda \in \Lambda^{n-1} : \lambda_1 = 0 \}$.
The asymptotic behavior of the frequency map in each one of these
three situations is expected to be dramatically different;
see the theorem below.
The last partition labels the component of the caustic parameter
that becomes singular:
$S_j = \{ \lambda \in \Lambda^{n-1} : \lambda_j = a_j \}$.
Besides, given any caustic parameter $\lambda \in \Lambda$
we shall denote by $\lambda^{S_j} \in S_j$ the caustic parameter obtained
from $\lambda$  by substituting its $j$-th component with $a_j$.
Finally, we introduce the $(n-1)$-dimensional set
\[
G_* =
\{(\lambda_2,\ldots,\lambda_n) \in \Rset^{n-1} :
  (0,\lambda_2,\ldots,\lambda_n) \in G \},
\]
which turns out to be the nonsingular caustic space for the
geodesic flow on the ellipsoid.
We note that $S_1 = \{c\} \times (H_1 \cup H_2)$,
$S_2 = (E \cup H_1) \times \{b\}$, and $G_\ast = H_1 \cup H_2$
with the notations used in the previous section for triaxial
ellipsoids of $\Rset^3$.

\begin{thm}\label{thm:nD}
The frequency map $\omega : \Lambda \to \Rset^n$ has the
following properties.
\begin{enumerate}
\item
It is analytic in $\Lambda$.
\item
It can be continuously extended to the border $\partial \Lambda$,
the extended map being as follows:
\begin{enumerate}
\item
It vanishes at $\bar{G}$;
\item
One of its components can be explicitly written as a function of the rest
 at $\bar{R}$;
\item
Its first component is equal to $1/2$ at $\bar{S}_1$;
\item
Its $l$-th component is equal to the $(l-1)$-th component at $\bar{S}_l$
for $2 \le l \le n$; and
\item
Its ``free'' components are an $(n-1)$-dimensional frequency of
the billiard inside the section of the original ellipsoid by a suitable
coordinate hyperplane at $\bar{R} \cup \bar{S}$.
\end{enumerate}
Besides, the restriction of the continuous extended map to any of
the $k$-dimensional connected components of $\Lambda^k$, $1 \le k \le n-1$,
is analytic.
\item
Its asymptotic behavior at $\Lambda^{n-1} = G \cup S \cup R$ is:
\begin{enumerate}
\item
$\omega(\lambda) =
 \kappa^G(\lambda_2,\ldots,\lambda_n) \lambda_1^{1/2} + \Or(\lambda_1^{3/2})$,
as $\lambda_1 \to 0^+$;
\item
$\omega(\lambda) - \omega(\lambda^{S_j}) \asymp \kappa^S(\lambda^{S_j}) /\log |a_j-\lambda_j|$,
as $\lambda_j \to a_j$;
\item
$\omega(\lambda) - \omega(\lambda^R) = \Or(\lambda-\lambda^R)$,
as $\lambda \to \lambda^R \in R$;
\end{enumerate}
for some analytic functions $\kappa^G : G_* \to \Rset_+^n$
and $\kappa^S: S \to \Rset^n$.
\end{enumerate}
\end{thm}

\proof
It follows from the same arguments and computations that in the spatial case.
The arguments are not repeated.
The computations with hyperelliptic integrals have been relegated
to~\ref{ap:HyperellipticIntegrals}.
\qed

We recall that, once fixed the parameters $a_1,\ldots,a_{n+1}$ of the ellipsoid and
the caustic parameters $\lambda_1,\ldots,\lambda_n$,
we write the $2n+1$ positive numbers
\[
\{c_1,\ldots,c_{2n+1}\} =
\{a_1,\ldots,a_{n+1}\} \cup \{\lambda_1,\ldots,\lambda_n\}
\]
in an ordered way: $c_0:= 0 < c_1 < \cdots < c_{2n+1}$.
Then the frequency $\omega(\lambda)$ is defined in terms of some
hyperelliptic integrals over the intervals $(c_{2j},c_{2j+1})$.
If two consecutive elements of $\{c_0,\ldots,c_{2n+1}\}$ collide,
then $\omega(\lambda)$ is, a priori, not well-defined.
Thus, it is natural to ask:
How does $\omega(\lambda)$ behave at these collisions?

In the previous theorem we have solved this question at the set
$\Lambda^{n-1} = G \cup R \cup S$,
which covers just the geodesic flow limit: $c_1 \to 0^+$,
the $n$ simple regular collapses: $c_{2l+1},c_{2l} \to c^*$ for some $l$,
and the $n$ simple singular collapses: $c_{2l-1},c_{2l} \to c^*$ for some $l$.
But there are many more (multiple) collapses, from double ones to total ones.
Double collapses correspond to the set $\Lambda^{n-2}$.
Total collapses have multiplicity $n$,
so they correspond to set of vertexes $\Lambda^0$.

We believe that it does not make sense to describe the asymptotic
behavior of the frequency map at all of them,
since the behavior in each case must be the expected one.
In order to convince the reader of the validity of this claim,
we end the paper with a couple of extreme cases.

As a first example,
let us consider the vertex $\hat{\lambda}=(a_1,\ldots,a_n) \in \Lambda^0$.
It represents the unique total singular collapse,
because it is the unique common vertex of the $2^n$ open connected
components of the caustic space:
\[
\bigcap_{\sigma \in \{0,1\}^n} \bar{\Lambda}_\sigma =
\{ \hat{\lambda} \} = \bigcap_{j=1}^n \bar{S}_j.
\]
Using that the point $\hat{\lambda}$ belongs to all the closures $\bar{S}_j$,
from theorem~\ref{thm:nD} we get that $\omega(\hat{\lambda}) = (1/2,\ldots,1/2)$.
Which is the asymptotic behavior of $\omega$ at this vertex?
In~\ref{sap:TotalSingularCollapse} it is proved that
\[
\omega(\lambda) = \omega(\hat{\lambda}) +
\Or(1/\log |a_1-\lambda_1|,\ldots,1/\log |a_n-\lambda_n|\big), \qquad
\lambda \to \hat{\lambda}.
\]
This behavior is singular in the $n$ caustic coordinates, as expected.

On the contrary,
the vertex $\tilde{\lambda}=(a_2,\ldots,a_{n+1}) \in \Lambda^0$ represents the
unique total regular collapse,
so we predict a regular behavior in the $n$ caustic coordinates.
In~\ref{sap:TotalRegularCollapse} we show that
$\omega(\tilde{\lambda}) = (\widetilde{\omega}_1,\ldots,\widetilde{\omega}_n)$,
where the limit frequencies $0 < \widetilde{\omega}_j < 1/2$ are defined
as $\sin^2 \pi \widetilde{\omega}_j = a_1/a_{j+1}$,
and the asymptotic behavior is
\[
\omega(\lambda) = \omega(\tilde{\lambda}) + \Or (\tilde{\lambda}-\lambda), \qquad
\lambda \to \tilde{\lambda}^-.
\]
Once more, the frequency map has the expected behavior.

\section{Conclusion and further questions}

We studied periodic trajectories of billiards inside
nondegenerate ellipsoids of $\Rset^{n+1}$.
First,
we trivially extended the definition of the frequency map
$\omega$ to any dimension,
we presented two conjectures about $\omega$ based on numerical computations,
and we deduced from the second one some lower bounds on the periods.
Next, we proved that $\omega$ can be continuously extended to any
singular value of the caustic parameters,
although it is exponentially sharp at the ``inner'' singular caustic
parameters.
Finally, we focused on ellipses and triaxial ellipsoids,
where we found examples of trajectories
whose periods coincide with the previous lower bounds.
We also computed several bifurcation curves.
Despite these results, many unsolved questions remain.
We indicate just four.

The most obvious challenge is to tackle any of the conjectures,
although it does not look easy.
We have already devoted some efforts without success.
We believe that the proof of any of these conjectures requires
either a deep use of algebraic geometry or to rewrite the
frequency map as the gradient of a ``Hamiltonian'';
see~\cite[\S 4]{WaalkensDullin2002}.

Another interesting question is to describe completely the phase space
of billiards inside ellipsoids in $\Rset^{n+1}$ for $n\ge 2$.
A rich hierarchy of invariant objects appears in these billiards:
Liouville maximal tori, low-dimensional tori, normally hyperbolic manifolds
whose stable and unstable manifolds are doubled, et cetera.
For instance, the stable and unstable invariant manifolds of the two-periodic
hyperbolic trajectory corresponding to an oscillation along the major axis
of the ellipsoid were fully described in~\cite{Delshams_etal2001}.

Third,
we plan to give a complete classification of the symmetric periodic
trajectories inside generic ellipsoids~\cite{CasasRamirez}.
To present the problem,
let us consider the symmetric periodic trajectories inside
an ellipse displayed in figure~\ref{fig:PeriodicTrajectories2D}.
On the one hand, the three-periodic trajectory drawn in a continuous red line
has an impact point on (and is symmetric with respect to) the $x$-axis.
On the other hand, the four-periodic trajectory drawn in a dashed green line
has a couple of segments passing through (and is symmetric with respect to)
the origin.
It is immediate to realize that there do not exist neither a trajectory
with a hyperbola as caustic like the first one, neither a trajectory
with an ellipse as caustic like the second one.
The problem consists of describing all possible kinds of symmetric
periodic trajectories once fixed the type of the $n$ caustics for
ellipsoids in $\Rset^{n+1}$.
Once these trajectories were well understood,
we could study their persistence under small symmetric
perturbations of the ellipsoid,
and the break-up of the Liouville tori on which they live.
Similar results have already been found in other billiard frameworks:
homoclinic trajectories inside ellipsoids of $\Rset^{n+1}$
with a unique major axis~\cite{Bolotin_etal2004},
and periodic trajectories inside circumferences of the
plane~\cite{RamirezRos2006}.

Finally, we look for simple formulae to express the caustic parameters
$\lambda_1,\ldots,\lambda_n$ that give rise to periodic trajectories of
small periods in terms of the parameters $a_1,\ldots,a_{n+1}$ of the
ellipsoid.
As a by-product of those formulae,
one can find algebraic expressions for the functions $g^\ast_j(b)$
that appear in numerical result~\ref{num:BifurcationsP}.
This is a work in progress~\cite{RamirezRos1}.

\ack
P. S. Casas was supported in part by MCyT-FEDER Grant MTM2006-00478 (Spain).
R. Ram{\'\i}rez-Ros was supported in part by MICINN-FEDER Grant MTM2009-06973 (Spain)
and CUR-DIUE Grant 2009SGR859 (Catalonia).
Useful conversations with Jaume Amor\'os, \`Alex Haro, Yuri Fedorov,
and Carles Sim\'o are gratefully acknowledged.

\appendix

\section{Computations with hyperelliptic integrals}
\label{ap:HyperellipticIntegrals}

\subsection{Technical lemmas}
\label{sap:TechnicalLemmas}

\begin{lem}\label{lem:0}
Let $f_\epsilon \in C^0([\alpha,\beta])$ be a family of functions
such that $f_\epsilon = f_0 + \Or(\epsilon)$ in the $C^0$-topology.
Then
\[
I_\epsilon =
\int_\alpha^\beta \frac{f_\epsilon(s) \rmd s}{\sqrt{(s-\alpha)(\beta-s)}} =
\int_\alpha^\beta \frac{f_0(s) \rmd s}{\sqrt{(s-\alpha)(\beta-s)}} +
\Or(\epsilon).
\]
\end{lem}
\proof
$|I_\epsilon - I_0| \le
 |f_\epsilon - f_0|_{C^0([\alpha,\beta])}
\int_\alpha^\beta ((s-\alpha)(\beta-s))^{-1/2} \rmd s =
\pi |f_\epsilon - f_0|_{C^0([\alpha,\beta])} = \Or(\epsilon)$.
\qed

\begin{lem}\label{lem:R}
Let $f \in C^1([m,M])$ with $m < \alpha < \beta < M$
and $\epsilon = \beta - \alpha$.
Then
\[
\int_\alpha^\beta \frac{f(s) \rmd s}{\sqrt{(s-\alpha)(\beta-s)}} =
\pi f(\alpha) + \Or(\epsilon) = \pi f(\beta) + \Or(\epsilon),\qquad
\epsilon \to 0^+.
\]
\end{lem}
\proof
Using the Mean Value Theorem for integrals,
we get that there exists some $s_0 \in [\alpha,\beta]$ such that
the integral is equal to
$f(s_0) \int_\alpha^\beta ((s-\alpha)(\beta-s))^{-1/2} \rmd s = \pi f(s_0)$.
\qed

\begin{lem}\label{lem:G}
Let $f\in C^1([0,M])$ with $0 < \epsilon < M$. Then
\[
I_\epsilon = \int_0^\epsilon \frac{f(s) \rmd s}{\sqrt{\epsilon-s}} =
2 f(0) \epsilon^{1/2} + \Or(\epsilon^{3/2}),\qquad
\epsilon \to 0^+.
\]
\end{lem}
\proof
$I_\epsilon = \left[-2 (\epsilon-s)^{1/2}f(s) \right]_{s=0}^{s=\epsilon}
+ 2 \int_0^\epsilon (\epsilon-s)^{1/2}f'(s)\rmd s =
2 f(0) \epsilon^{1/2} + \Or(\epsilon^{3/2})$.
\qed

\begin{lem}\label{lem:S}
Let $f \in C^1([\alpha,\beta])$.
Set
$\eta = f(\alpha) \log(4\beta-4\alpha) +
 \int_\alpha^\beta (s-\alpha)^{-1}(f(s) - f(\alpha)) \rmd s$,
$\xi = \int_\alpha^\beta (s-\alpha)^{-3/2} (f(s) - f(\alpha) ) \rmd s$,
$\mu = f(\beta) \log(4\beta-4\alpha) +
           \int_\alpha^\beta (\beta-s)^{-1} (f(s)-f(\beta) ) \rmd s$, and
$\psi = \int_\alpha^\beta (\beta-s)^{-3/2} (f(s) - f(\beta)) \rmd s$.
Then
\begin{eqnarray*}
\int_\alpha^\beta \frac{f(s) \rmd s}{\sqrt{(s+\epsilon-\alpha)(s-\alpha)}} =
- f(\alpha) \log \epsilon + \eta + \Or(\epsilon \log \epsilon),\qquad
\epsilon \to 0^+, \\
\int_\alpha^\beta \frac{f(s) \rmd s}{(s+\epsilon-\alpha)\sqrt{s-\alpha}} =
\pi f(\alpha) \epsilon^{-1/2} + \xi + \Or(\epsilon^{1/2}),\qquad
\epsilon \to 0^+,\\
\int_\alpha^\beta \frac{f(s) \rmd s}{\sqrt{(\beta+\epsilon-s)(\beta -s)}} =
- f(\beta) \log \epsilon + \mu + \Or(\epsilon \log \epsilon),\qquad
\epsilon \to 0^+, \\
\int_\alpha^\beta \frac{f(s) \rmd s}{(\beta+\epsilon-s)\sqrt{\beta -s}} =
\pi f(\beta) \epsilon^{-1/2} + \psi + \Or(\epsilon^{1/2}),\qquad
\epsilon \to 0^+.
\end{eqnarray*}
The first (respectively, last) two estimates also hold when $f$ has
a singularity at $s=\beta$ (respectively, at $s=\alpha$),
provided $f \in L^1([\alpha,\beta])$.
\end{lem}
\proof
We split the first integral as
$I_\epsilon = \tilde{\eta} + \hat{I}_\epsilon - \tilde{I}_\epsilon$,
where $\tilde{\eta} = \int_\alpha^\beta (s-\alpha)^{-1}(f(s) - f(\alpha)) \rmd s$
is a constant, and
\[
\fl
\hat{I}_\epsilon =
\int_\alpha^\beta \frac{f(\alpha) \rmd s}{\sqrt{(s+\epsilon-\alpha)(s-\alpha)}},\qquad
\tilde{I}_\epsilon =
\int_\alpha^\beta \frac{f(s) - f(\alpha)}{s-\alpha}
                  \left(1-\sqrt{\frac{s-\alpha}{s+\epsilon-\alpha}}\right)\rmd s.
\]
By performing the change $x^2=s-\alpha$ in the integral $\hat{I}_\epsilon$,
we get that
\[
\fl
\hat{I}_\epsilon =
2 \int_0^{\sqrt{\beta-\alpha}} \frac{f(\alpha) \rmd x}{\sqrt{x^2 + \epsilon}} =
2 f(\alpha) \left[\log\left(x+\sqrt{x^2+\epsilon}\right) \right]_{x=0}^{x=\sqrt{\beta-\alpha}}
 = - f(\alpha) \log \epsilon + \hat{\eta} + \Or(\epsilon),
\]
where $\hat{\eta} = f(\alpha) \log(4\beta-4\alpha)$ is another constant.
Thus,
to get the first formula with constant $\eta = \hat{\eta} + \tilde{\eta}$
it suffices to see that $\tilde{I}_\epsilon = \Or(\epsilon \log \epsilon)$.

Once fixed some $\gamma \in (\alpha,\beta)$,
we decompose the integral $\tilde{I}_\epsilon$ as the sum
$\tilde{J}_\epsilon + \tilde{K}_\epsilon$, where
$\tilde{J}_\epsilon = \int_\alpha^\gamma \tilde{f}(s) r_\epsilon(s) \rmd s$,
$\tilde{K}_\epsilon = \int_\gamma^\beta \tilde{f}(s) r_\epsilon(s) \rmd s$,
and
\[
\tilde{f}(s) = \frac{f(s) - f(\alpha)}{s-\alpha},\qquad
r_\epsilon(s) = 1-\sqrt{\frac{s-\alpha}{s+\epsilon-\alpha}}.
\]
First, we consider the interval $[\alpha,\gamma]$.
Then
$|\tilde{f}|_\infty = \max \{ |\tilde{f}(s)| : \alpha \le s \le \gamma \} < \infty$
and $r_\epsilon(s)$ is positive in $[\alpha,\gamma]$.
Set $\delta=\gamma-\alpha$.
Using again the change $x^2=s-\alpha$, we see that
\begin{eqnarray*}
\fl
|\tilde{f}|_\infty^{-1} |\tilde{J}_\epsilon | & \le &
\int_\alpha^\gamma r_\epsilon(s) \rmd s =
\delta - \int_\alpha^\gamma \sqrt{\frac{s-\alpha}{s+\epsilon-\alpha}} \rmd s =
\delta - 2 \int_0^{\sqrt{\delta}} \frac{x^2 \rmd x}{\sqrt{x^2+\epsilon}} \\
\fl
& = &
\delta -
\left[ x \sqrt{x^2 + \epsilon} + \epsilon \log\left(x+\sqrt{x^2+\epsilon}\right)
 \right]_{x=0}^{x=\sqrt{\delta}}
= -\frac{\epsilon}{2} \log \epsilon + \Or(\epsilon) = \Or(\epsilon \log \epsilon).
\end{eqnarray*}
Concerning the other interval,
we note that $r_\epsilon(s)$ is positive and decreasing in $[\gamma,\beta]$.
Hence,
$\max \{ |r_\epsilon(s)| : \gamma \le s \le \beta \} = r_\epsilon(\gamma)$,
and so
$|\tilde{K}_\epsilon| \le
 r_\epsilon(\gamma) \int_\gamma^\beta |\tilde{f}(s)| \rmd s =
 \Or(\epsilon)$.
This ends the proof of the first formula.

We split the second integral as
$L_\epsilon = \xi + \hat{L}_\epsilon - \tilde{L}_\epsilon$,
where $\xi$ is the constant given in the statement of the lemma, and
\[
\fl
\hat{L}_\epsilon =
\int_\alpha^\beta
\frac{f(\alpha) \rmd s}{(s+\epsilon-\alpha)\sqrt{s-\alpha}},\qquad
\tilde{L}_\epsilon =
\int_\alpha^\beta \frac{f(s) - f(\alpha)}{(s-\alpha)^{3/2}}
                  \left(1-\frac{s-\alpha}{s+\epsilon-\alpha}\right)\rmd s.
\]
By performing the change $x = s-\alpha$ in the integral $\hat{L}_\epsilon$,
we get that
\[
\fl
\hat{L}_\epsilon =
\int_0^{\beta-\alpha} \frac{f(\alpha) \rmd x}{(x+\epsilon)\sqrt{x}} =
2f(\alpha)\epsilon^{-1/2}
\left[\atan \left(\epsilon^{-1/2}x^{1/2}\right) \right]_{x=0}^{x=\beta-\alpha} =
\pi f(\alpha) \epsilon^{-1/2} + \Or(\epsilon^{1/2}).
\]
Thus,
to get the second formula it suffices to see that
$\tilde{L}_\epsilon = \Or(\epsilon^{1/2})$,
which follows from similar computations than the ones above.

The last formulae are obtained by performing the change
of variables $s-\alpha = \beta -t$ in the former ones. 
\qed

\begin{cor}\label{cor:S}
Let $f \in C^1([m,M])$ with $m < \alpha_- < \alpha_+ < \beta_- < \beta_+ < M$, and
\[
I = I(\alpha_-,\alpha_+,\beta_-,\beta_+) =
\int_{\alpha_+}^{\beta_-}
\frac{f(s) \rmd s}{\sqrt{(s - \alpha_-)(s - \alpha_+)(\beta_- -s)(\beta_+ -s)}}.
\]
Let $\alpha_*$ and $\beta_*$ be two reals such that $m < \alpha_* < \beta_* < M$.
Let $\epsilon = (\epsilon_1,\epsilon_2) \in \Rset_+^2$ with
$\epsilon_1 = \alpha_+ - \alpha_-$ and  $\epsilon_2 = \beta_+ - \beta_-$.
Then there exists a constant $\zeta \in \Rset$ such that
\[
\fl
I =
-\frac{f(\alpha_*) (1+ \Or(\epsilon_2) ) \log \epsilon_1 +
f(\beta_*) (1+ \Or(\epsilon_1) ) \log \epsilon_2}{\beta_*-\alpha_*} + \zeta +
 \Or(\epsilon_1 \log \epsilon_1,\epsilon_2 \log \epsilon_2),
\]
as $\alpha_\pm \to \alpha_*$ and $\beta_\pm \to \beta_*$,
so that $\epsilon = (\epsilon_1,\epsilon_2) \to (0^+,0^+)$.
\end{cor}

\proof
It follows by applying the first and third estimates of the previous lemma
to the integrals $\int_{\alpha_+}^\gamma$ and $\int_\gamma^{\beta_-}$
for some point $\gamma \in (\alpha_+,\beta_-)$,
although before we must fix the lower limit of the first integral with
the change $x-\alpha_*=s-\alpha_+$,
and the upper limit of the second integral with the change
$x-\beta_*=s-\beta_-$.
\qed

\begin{lem}\label{lem:LinearSystems}
Let $\mathbf{K}_\epsilon \omega_\epsilon = \tau_\epsilon$ be a family of
square linear systems defined for $\epsilon > 0$.
\begin{enumerate}
\item
If the limits $\mathbf{K} = \lim_{\epsilon \to 0^+} \mathbf{K}_\epsilon$ and
$\tau = \lim_{\epsilon \to 0^+} \tau_\epsilon$ exist,
and $\mathbf{K}$ is nonsingular, then
\[
\omega_\epsilon = \omega +
\Or(|\mathbf{K}_\epsilon - \mathbf{K}|, | \tau_\epsilon - \tau | ),
\qquad \epsilon \to 0^+,
\]
where $\omega=\mathbf{K}^{-1} \tau$ is the unique solution of the nonsingular
limit system $\mathbf{K} \omega = \tau$.

\item
If, in addition, the matrix $\mathbf{K}_\epsilon$ and the vector $\tau_\epsilon$ are
differentiable at $\epsilon=0$,
then the solution $\omega_\epsilon$ also is differentiable at $\epsilon=0$.
To be more precise, if
\[
\mathbf{K}_\epsilon = \mathbf{K} + \epsilon \mathbf{L} + \order(\epsilon),\qquad
\tau_\epsilon = \tau + \epsilon \zeta + \order(\epsilon),\qquad
\epsilon \to 0^+,
\]
for some square matrix $\mathbf{L}$ and some vector $\zeta$, then
\[
\omega_\epsilon = \omega + \epsilon \kappa + \order(\epsilon),\qquad
\epsilon \to 0^+,
\]
where $\omega=\mathbf{K}^{-1} \tau$ and
$\kappa = \mathbf{K}^{-1}(\zeta- \mathbf{L}\omega)$.
\end{enumerate}
\end{lem}
\proof
Both results follow directly from classical error bounds
in numerical linear algebra.
See, for instance,~\cite[\S 2.7]{GolubVanLoan1996}.
\qed

\subsection{Properties of the rotation number}
\label{sap:Proof_2D}

Let us write the rotation number as the quotient
$\rho(\lambda) = \Delta(\lambda)/2K(\lambda)$,
where
\[
\Delta(\lambda) = \int_0^{\min(b,\lambda)} \frac{\rmd s}{\sqrt{T(s)}},
\qquad
K(\lambda) = \int_{\max(b,\lambda)}^a \frac{\rmd s}{\sqrt{T(s)}},
\]
and $T(s) = (\lambda-s)(b-s)(a-s)$.

The study of the limit $\lambda \to 0^+$ is easy.
From lemma~\ref{lem:G}, we get the estimate
$\Delta(\lambda) = 2(ab)^{-1/2} \lambda^{1/2} + \Or(\lambda^{3/2})$,
as $\lambda \to 0^+$,
whereas from lemma~\ref{lem:0} we get that
\[
K(\lambda) =
\int_b^a \frac{\rmd s}{\sqrt{s(s-b)(a-s)}} + \Or(\lambda),\qquad
\lambda \to 0^+.
\]
By combining both estimates we get that
$\rho(\lambda)  = \kappa^G \lambda^{1/2} + \Or(\lambda^{3/2})$,
so $\lim_{\lambda \to 0^+} \rho(\lambda) = 0$.

Next, we consider $\lambda \to b^-$.
After some computations based on lemma~\ref{lem:S},
we get
\begin{eqnarray*}
K(b-\epsilon) & = &
- c^{-1/2} \log \epsilon + \eta + \Or(\epsilon \log \epsilon),
\qquad \epsilon \to 0^+, \\
\Delta(b-\epsilon) & = &
-c^{-1/2} \log \epsilon + \mu + \Or(\epsilon \log \epsilon),
\qquad \epsilon \to 0^+,
\end{eqnarray*}
where $c = a-b$, $\eta = \hat{\eta} + \tilde{\eta}$,
$\mu=\hat{\mu}+\tilde{\mu}$, with
$\hat{\eta} = c^{-1/2} \log 4c$, $\hat{\mu} = c^{-1/2} \log 4b$, and
\begin{eqnarray*}
\fl
\tilde{\eta} & = &
 \int_b^a \left( \frac{1}{\sqrt{a-s}} - \frac{1}{\sqrt{a-b}} \right) \frac{\rmd s}{s-b} =
 \frac{2}{\sqrt{c}} \int_0^{\sqrt{c}} \frac{\rmd x}{x+\sqrt{c}} =
 \frac{\log 4}{\sqrt{c}}, \\
\fl
\tilde{\mu}  & = &
 \int_0^b \left( \frac{1}{\sqrt{a-s}} - \frac{1}{\sqrt{a-b}} \right) \frac{\rmd s}{b-s} =
 -\frac{2}{\sqrt{c}} \int_{\sqrt{c}}^{\sqrt{a}} \frac{\rmd x}{x+\sqrt{c}} =
 \frac{1}{\sqrt{c}} \log \frac{4c}{(\sqrt{a}+\sqrt{c})^2}.
\end{eqnarray*}
We have used the change $x^2 = a-s$ in both integrals.
Let $\eta_* = c^{1/2} \eta = \log 16c$ and
$\mu_* = c^{1/2} \mu = \log 16bc(a^{1/2}+c^{1/2})^{-2}$.
Then we have the estimate
\begin{equation}\label{eq:Slog}
\fl
2\rho(b-\epsilon) =
\frac{\Delta(b-\epsilon)}{K(b-\epsilon)} =
\frac{1 - c^{1/2} \mu \log^{-1} \epsilon + \Or(\epsilon)}
     {1 - c^{1/2} \eta \log^{-1} \epsilon + \Or(\epsilon)} =
\frac{1 - \mu_* \log^{-1} \epsilon}{1 - \eta_* \log^{-1} \epsilon}
 + \Or(\epsilon),
\end{equation}
as $\epsilon \to 0^+$.
Thus,
$\kappa^S = (\eta_* - \mu_*)/2 =
 \log \left( (a/b)^{1/2}+(c/b)^{1/2} \right) =
 \log \left(d + (d^2-1)^{1/2} \right) = \acosh d$, where $d=(a/b)^{1/2}$.
This implies that $\cosh^2 \kappa^S = a/b$.
Besides, estimate~(\ref{eq:Slog}) is the key to prove that
the caustic parameter $\lambda_-^0$ is exponentially close to $b$.
Once fixed $\rho^0 \lesssim 1/2$,
let $\lambda^0_-\in E$ be the unique caustic parameter such that
$\rho(\lambda_-^0 ) = \rho^0$, $0 < \epsilon = b-\lambda^0_- \ll 1$,
and $\delta = \log^{-1} \epsilon$.
By finding $\delta^{-1} = \log \epsilon$ in estimate~(\ref{eq:Slog}),
we get
\[
\log \epsilon = 1/\delta =
\eta_* + \frac{\mu_* - \eta_*}{1 - 2\rho^0} + \Or(\epsilon).
\]
Using that $\kappa^S = (\eta_* - \mu_*)/2$ and $\eta_* = \log 16c$,
we check that $\lambda^0_- = b - \epsilon$, with
\[
\epsilon =
\rme^{1/\delta} =
\rme^{\eta_* - 2\kappa^S/(1-2\rho^0) + \Or(\epsilon)} =
16c \rme^{-\kappa^S/(1/2-\rho^0)} +
\Or\left (\rme^{-2\kappa^S/(1/2-\rho^0)} \right),
\]
as $\rho^0 \to (1/2)^-$.
The limit $\lambda \to b^+$ is completely analogous.
We omit the computations.

With respect to the limit $\lambda \to a^-$,
we note that $T(s) = (b-s)(a-s)^2 + \Or(a-\lambda)$ uniformly
in the interval $[0,b]$.
Hence,
\[
\fl
\Delta(\lambda) =
\int_0^b \frac{\rmd s}{(a-s)\sqrt{b-s}} + \Or(a-\lambda) =
\frac{2}{\sqrt{a-b}} \atan \sqrt{\frac{b}{a-b}}+\Or(a-\lambda),\quad
\lambda \to a^-.
\]
Besides, from lemma~\ref{lem:R} we get the estimate
$K(\lambda) = \pi (a-b)^{-1/2} + \Or(a-\lambda)$,
as $\lambda \to a^-$.
Therefore, $\rho(\lambda) = \varrho + \Or(a-\lambda)$,
as $\lambda \to a^-$,
where the limit value $\varrho \in (0,1/2)$
is defined by $\tan^2 \pi \varrho = b/(a-b)$.
That is, $\sin^2 \pi \varrho = b/a$.

\subsection{Another characterization of the frequency}
\label{sap:Frequency}

We associate a ``frequency'' $\omega=\varpi(c) \in \Rset^n$ to any
vector $c=(c_1,\ldots,c_{2n+1}) \in \Rset^{2n+1}$ such that
$c_0 := 0 < c_1 < \cdots < c_{2n+1}$ in the following way.
First, we consider:
\begin{itemize}
\item
The polynomial
$T(s) = \prod_{i=1}^{2n+1} (c_i -s) \in \Rset_{2n+1}[s]$,
which is positive in the $n+1$ intervals of the form $(c_{2j},c_{2j+1})$;
\item
The $n+1$ linear functionals
$P(s) \mapsto \mathcal{K}_j [P(s)] =
 \int_{c_{2j}}^{c_{2j+1}} (T(s))^{-1/2}P(s) \rmd s$;
\item
The $n+1$ column vectors
$K_j = (\mathcal{K}_j[1],\mathcal{K}_j[s],\ldots,\mathcal{K}_j[s^{n-1}])^t
 \in \Rset^n$;
\item
The $n\times n$ nonsingular matrix $\mathbf{K} = (-K_1,\ldots,(-1)^n K_n)$; and
\item
The linear functionals
$\mathcal{K}[P(s);\omega] =
\mathcal{K}_0[P(s)] +
 2 \sum_{j=1}^n (-1)^j \omega_j \mathcal{K}_j[P(s)]$, for $\omega \in \Rset^n$.
\end{itemize}
The hypothesis $c_1 > 0$ is not essential to get a nonsingular
matrix $\mathbf{K}$,
but it suffices to assume the strict inequalities $c_1 < \cdots < c_{2n+1}$;
see~\cite[\S III.3]{Griffiths1989}.

\begin{lem}
There exists an unique $\omega \in \Rset^n$ such that
\begin{equation}\label{eq:Frequency2}
\mathcal{K}[P(s);\omega] = 0,\qquad
\forall P(s) \in \Rset_{n-1}[s],
\end{equation}
or equivalently, such that $K_0 + 2 \mathbf{K}\omega = 0$,
which is the matricial form of the linear system given in~(\ref{eq:Frequency}).
\end{lem}

\proof
By taking the basis $\{1,s,\ldots,s^{n-1}\}$ of $\Rset_{n-1}[s]$,
we see that condition~(\ref{eq:Frequency2}) is equivalent to
the linear system $K_0+ 2 \mathbf{K}\omega = 0$.
\qed

Therefore, condition~(\ref{eq:Frequency2}) is an equivalent
characterization of the frequency.
From now on, $\omega=\varpi(c)$ stands for the frequency computed
through the previous steps.

\subsection{Geodesic flow limit: $c_1 \to 0^+$}\label{ap:G}
\label{sap:GeodesicFlowLimit}

Let $K^G_{00} = 2 (\prod_{i=2}^n c_i)^{-1/2}$,
$K_0^G = (K^G_{00},0,\ldots,0) \in \Rset^n$, and
$T^G(s) = -s \prod_{i=2}^{2n+1} (c_i -s)$.
Let $\mathbf{K}^G$ be the $n\times n$ nonsingular matrix associated to
the vector $c^G = (0,c_2,\ldots,c_{2n+1})$.
Let $\kappa^G \in \Rset^n$ be the unique solution of the
linear system $K^G_0 + 2 \mathbf{K}^G \kappa^G = 0$.
Then
\begin{equation}\label{eq:Glimit}
\omega = \kappa^G c_1^{1/2} + \Or(c_1^{3/2}),\qquad c_1 \to 0^+.
\end{equation}

The proof is short.
First, we note that $T = T^G + \Or(c_1)$ uniformly in $[0,c_{2n+1}]$.
Thus, using lemma~\ref{lem:0}, we get that
$\mathbf{K} = \mathbf{K}^G + \Or(c_1)$ as $c_1 \to 0^+$.
And using lemma~\ref{lem:G} we see that
$K_0 = K_0^G c^{1/2}_1 + \Or(c^{3/2}_1)$ as $c_1 \to 0^+$.
Therefore,
the linear systems $K^G_0 + 2 \mathbf{K}^G \kappa^G = 0$ and
$c_1^{-1/2} K_0 + 2 \mathbf{K} (c_1^{-1/2}\omega) = 0$
are $\Or(c_1)$-close, being $\mathbf{K}^G$ nonsingular,
so~(\ref{eq:Glimit}) follows from the first item in lemma~\ref{lem:LinearSystems}.

\subsection{Simple regular collapse:
            $c_{2l+1},c_{2l} \to c^*$ for some $l=1,\ldots,n$}
\label{sap:SimpleRegularCollapse}

Set $c^R=(c_1,\ldots,c_{2l-1},c_{2l+2},\ldots,c_{2n+1}) \in \Rset^{2n-1}$.
Let $T^R(s) = \prod_{i \neq 2l,2l+1} (c_i - s)$ be the polynomial
associated to $c^R$.
Let $\mathcal{K}^R_j$ and $\mathcal{K}^R$ be the functionals associated to $c^R$.
Let $\omega^R = (\omega^R_1,\ldots,\omega^R_n) \in\Rset^n$, where
$\omega^R_{\neq l} := (\omega^R_1,\ldots,\omega^R_{l-1},\omega^R_{l+1},\ldots,\omega^R_n) =
 \varpi(c^R) \in \Rset^{n-1}$ is the frequency associated to $c^R$,
and $\omega^R_l \in \Rset$ is determined by
\begin{equation}\label{eq:Rcondition}
\fl
\int_0^{c_1} \frac{\rmd s}{|c^*-s|\sqrt{T^R(s)}} +
2 \sum_{j\neq l} \int_{c_{2j}}^{c_{2j+1}} \frac{(-1)^j \omega^R_j \rmd s}{|c^*-s|\sqrt{T^R(s)}}
+ \frac{(-1)^l 2 \pi \omega^R_l}{\sqrt{-T^R(c^*)}} = 0.
\end{equation}
Let $\epsilon = c_{2l+1}-c_{2l}$.
Then
\begin{equation}\label{eq:Rlimit}
\omega = \omega^R + \Or(\epsilon),\qquad
c_{2l+1},c_{2l} \to c^*.
\end{equation}

In order to prove this claim,
we first observe that characterization~(\ref{eq:Frequency2}) is equivalent
to the system of $n$ linear equations
\begin{equation}\label{eq:R1}
\left\{
\begin{array}{l}
\mathcal{K}[(c^*-s)s^i;\omega] = 0 \mbox{ for $i=0,\ldots,n-2$} \\
\mathcal{K}[1;\omega] = 0
\end{array}
\right. ,
\end{equation}
because $\{1,c^*-s,\ldots,(c^*-s)s^{n-2}\}$ is a basis of $\Rset_{n-1}[s]$.
Now, using lemmas~\ref{lem:0} and~\ref{lem:R}, we deduce the estimates
\[
\fl
\mathcal{K}_j[(c^* -s)s^i] =
 \textstyle{\int_{c_{2j}}^{c_{2j+1}}
 \frac{(c^*-s) s^i \rmd s}{|c^*-s| \sqrt{T^S(s)} + \Or(\epsilon) }} =
\cases{ \mathcal{K}^R_j[s^i] + \Or(\epsilon)     & \mbox{if $j < l$}, \\
        \Or(\epsilon)                            & \mbox{if $j = l$}, \\
       -\mathcal{K}^R_{j-1}[s^i] + \Or(\epsilon) & \mbox{if $j > l$};}
\]
and
\[
\mathcal{K}_j[1] =
 \cases{\pi\left(-T^R(c^*)\right)^{-1/2} + \Or(\epsilon) & \mbox{if $j = l$}, \\
        \int_{c_{2j}}^{c_{2j+1}} \frac{\rmd s}{|c^*-s|\sqrt{T^R(s)}}
           + \Or(\epsilon) & \mbox{otherwise}.}
\]

Therefore, the linear system~(\ref{eq:R1}) is $\Or(\epsilon)$-close to the
nonsingular linear system
\[
\left\{
\begin{array}{l}
\mathcal{K}^R[s^i;\omega^R_{\neq l}] = 0 \mbox{ for $i=0,\ldots,n-2$} \\
\mbox{condition~(\ref{eq:Rcondition})}
\end{array}
\right. ,
\]
and the asymptotic formula~(\ref{eq:Rlimit}) follows from
the first item in lemma~\ref{lem:LinearSystems}.

\subsection{Simple singular collapse:
            $c_{2l-1},c_{2l} \to c^*$ for some $l=1,\ldots,n$}
\label{sap:SimpleSingularCollapse}

Set $c^S=(c_1,\ldots,c_{2l-2},c_{2l+1},\ldots,c_{2n+1}) \in \Rset^{2n-1}$.
Let $T^S(s) = \prod_{i \neq 2l-1,2l} (c_i - s)$ be the polynomial
associated to $c^S$.
Let $\mathcal{K}^S_j$ and $\mathcal{K}^S$ be the functionals associated to $c^S$.
Let $\omega^S = (\omega^S_1,\ldots,\omega^S_n) \in\Rset^n$, where
$\omega^S_{\neq l} := (\omega^S_1,\ldots,\omega^S_{l-1},\omega^S_{l+1},\ldots,\omega^S_n) =
 \varpi(c^S) \in \Rset^{n-1}$
and
\[
\omega^S_l = \cases{1/2 & \mbox{if $l=1$} \\ \omega^S_{l-1} & \mbox{otherwise}}.
\]
Let $\epsilon = c_{2l}-c_{2l-1} > 0$ and $\delta= |\log \epsilon|^{-1} > 0$.
Then there exists $\kappa^S \in \Rset^n$ such that
\begin{equation}\label{eq:Slimit}
\omega = \omega^S + \delta \kappa^S + \order(\delta),\qquad
c_{2l-1},c_{2l} \to c^*.
\end{equation}
To prove this claim, we set $d = \sqrt{T^S(c^*)} > 0$.
We know that characterization~(\ref{eq:Frequency2}) is equivalent
to the system of $n$ linear equations
\begin{equation}\label{eq:S}
\left\{
\begin{array}{l}
\mathcal{K}[\delta d;\omega] = 0 \\
\mathcal{K}[(c^*-s)s^i;\omega] = 0 \mbox{ for $i=0,\ldots,n-2$}
\end{array}
\right. ,
\end{equation}
because $\{\delta d,c^*-s,\ldots,(c^*-s)s^{n-2}\}$
is a basis of $\Rset_{n-1}[s]$.
Now, using lemmas~\ref{lem:0} and~\ref{lem:S},
we deduce the following asymptotic estimates.
On the one hand,
there exist some constants $\zeta_0,\zeta_1,\ldots,\zeta_n \in \Rset$ such that
\[
\mathcal{K}_j[\delta d] = \delta d \mathcal{K}_j[1] =
\cases{1 + \zeta_j \delta + \Or(\epsilon)   & \mbox{if $j = l-1,l$}, \\
       \zeta_j \delta + \Or(\epsilon\delta) & \mbox{otherwise}.}
\]
On the other hand,
\[
\fl
\mathcal{K}_j[(c^* -s)s^i] =
 \textstyle{\int_{c_{2j}}^{c_{2j+1}}
 \frac{(c^*-s) s^i \rmd s}{|c^*-s| \sqrt{T^S(s)} + \Or(\epsilon) }} =
\cases{ \mathcal{K}^S_j[s^i] + \Or(\epsilon)      & \mbox{if $j < l-1$}, \\
        \textstyle{\int_{c_{2l-2}}^{c^*} \frac{s^i \rmd s}{\sqrt{T^S(s)}}}
           + \Or(\epsilon) & \mbox{if $j=l-1$}, \\
       -\textstyle{\int_{c^*}^{c_{2l+1}} \frac{s^i \rmd s}{\sqrt{T^S(s)}}}
           + \Or(\epsilon) & \mbox{if $j=l$}, \\
       -\mathcal{K}^S_{j-1}[s^i] + \Or(\epsilon) & \mbox{if $j > l$}.}
\]
In particular,
$\mathcal{K}_{l-1}[(c^* -s)s^i] - \mathcal{K}_l [(c^*-s)s^i] =
\int_{c_{2l-2}}^{c_{2l+1}} \frac{s^i \rmd s}{\sqrt{T^S(s)}} + \Or(\epsilon) =
 \mathcal{K}^S_{l-1}[s^i] + \Or(\epsilon)$.

We assume now that $l\neq 1$. The case $l=1$ is studied later on.
Since $\epsilon \ll \delta$,
the linear system~(\ref{eq:S}) is $\Or(\delta)$-close to the
nonsingular linear system
\[
\left\{
\begin{array}{l}
2 (-1)^{l-1} (\omega^S_{l-1} - \omega^S_l) = 0 \\
\mathcal{K}^S[s^i;\omega^S_{\neq l}] =
2(-1)^l (\omega^S_l - \omega^S_{l-1})
\textstyle{\int_{c^*}^{c_{2l+1}} \frac{s^i \rmd s}{\sqrt{T^S(s)}}}
\mbox{ for $i=0,\ldots,n-2$}
\end{array}
\right. ,
\]
which in its turn is equivalent to the linear system
\begin{equation}\label{eq:S0}
\left\{
\begin{array}{l}
\omega^S_l = \omega^S_{l-1} \\
\mathcal{K}^S[s^i;\omega^S_{\neq l}] = 0 \mbox{ for $i=0,\ldots,n-2$}
\end{array}
\right.
\end{equation}
whose unique solution is $\omega^S_{\neq l} = \varpi(c^S)$ and $\omega^S_l = \omega^S_{l-1}$.

Thus, the asymptotic formula $\omega = \omega^S + \Or(\delta)$ follows from
the first item in lemma~\ref{lem:LinearSystems}.
In fact, this result can be improved using the second item
in lemma~\ref{lem:LinearSystems}.
It suffices to note that the linear system~(\ref{eq:S}) is not only
$\Or(\delta)$-equivalent to~(\ref{eq:S0}),
but it is differentiable at $\delta = 0$.
Hence, (\ref{eq:Slimit}) holds for some vector $\kappa^S$ that could be
explicitly computed in terms of the limit system
and the constants $\zeta_0,\ldots,\zeta_n$.

If $l=1$, the linear system~(\ref{eq:S}) is $\Or(\delta)$-equivalent
to the nonsingular linear system
\[
\left\{
\begin{array}{l}
\omega^S_1 = 1/2 \\
\mathcal{K}^S[s^i;\omega^S_{\neq 1}] = 0 \mbox{ for $i=0,\ldots,n-2$}
\end{array}
\right. ,
\]
and the proof ends with just the same arguments that for $l \neq 1$.
We omit the details.

\subsection{Total regular collapse:
            $c_{2l+1},c_{2l} \to c_l^*$ for all $l=1,\ldots,n$}
\label{sap:TotalRegularCollapse}

Let us study the case of $n$ simultaneous collapses,
all of them regular.
That is, once fixed a vector $c^* =(c^*_1,\ldots,c^*_n) \in \Rset^n$
such that $0 < c_1 < c^*_1 < \cdots < c^*_n$,
we study the asymptotic behavior of the frequency $\omega=\varpi(c)$
when $c_{2l+1},c_{2l} \to c_l^*$ for all $l=1,\ldots,n$.
Let $\widetilde{\omega} =
(\widetilde{\omega}_1,\ldots,\widetilde{\omega}_n) \in \Rset^n$ be
the vector whose components verify that $0 < \widetilde{\omega}_l < 1/2$
and $\sin^2 \pi \widetilde{\omega}_l = c_1/c^*_l$.
Let $\epsilon = (\epsilon_1,\ldots,\epsilon_n) \in \Rset_+^n$ with
$\epsilon_l = c_{2l+1} - c_{2l}$.
Then
\begin{equation}\label{eq:omegatilde}
\omega = \widetilde{\omega} + \Or(\epsilon),\qquad
\epsilon \to (0^+,\ldots,0^+).
\end{equation}

Let $Q_l = \sqrt{c^*_l-c_1}\prod_{i \neq l} |c^*_i - c^*_l| > 0$.
Let $\{P_1(s),\ldots,P_n(s)\}$ be the basis of $\Rset_{n-1}[s]$
univocally determined by the interpolating conditions
\[
P_l(c^*_j) =
\cases{ Q_l & \mbox{if $j = l$}, \\
        0   & \mbox{otherwise}.}
\]
That is, $P_l(s) = (-1)^{l-1}\sqrt{c^*_l-c_1}\prod_{i \neq l} (c^*_i -s)$.
Using lemma~\ref{lem:R}, we get the estimates
\[
\fl
\mathcal{K}_0[P_l(s)] =
\int_0^{c_1} \left(\frac{(-1)^{l-1} \sqrt{c^*_l-c_1}}{(c^*_l -s)\sqrt{c_1-s}} + \Or(\epsilon)\right) \rmd s
        = 2 (-1)^{l-1} \atan \sqrt{\frac{c_1}{c^*_l-c_1}} + \Or(\epsilon),
\]
$\mathcal{K}_l[P_l(s)] = \pi + \Or(\epsilon)$,
and $\mathcal{K}_j[P_l(s)] = \Or(\epsilon)$ for $j \neq 0,l$.
Thus, the $n \times n$ linear system
\[
\mathcal{K}[P_l(s); \omega] = 0 \mbox{ for $l=1,\ldots,n$}
\]
is $\Or(\epsilon)$-close to the nonsingular decoupled linear system
\[
2(-1)^{l-1}\left( \atan \sqrt{\frac{c_1}{c^*_l-c_1}} - \pi \widetilde{\omega}_l \right) = 0
\mbox{ for $l=1,\ldots,n$},
\]
whose unique solution is given by $\tan^2 \pi \widetilde{\omega}_l = c_1/(c^*_l-c_1)$,
and so, by $\sin^2 \pi \widetilde{\omega}_l = c_1/c^*_l$.
Hence, the asymptotic formula~(\ref{eq:omegatilde})
follows from the first item in lemma~\ref{lem:LinearSystems}.

\subsection{Total singular collapse:
            $c_{2l-1},c_{2l} \to c^*_l$ for all $l=1,\ldots,n$}
\label{sap:TotalSingularCollapse}

Let $\widehat{\omega} = (1/2,\ldots,1/2) \in \Rset^n$,
$\epsilon = (\epsilon_1,\ldots,\epsilon_n) \in \Rset_+^n$,
and $\delta=(\delta_1,\ldots,\delta_n) \in \Rset_+^n$, where
$\epsilon_l = c_{2l} - c_{2l-1}$ and $\delta_l = |\log \epsilon_l|^{-1}$.
Then
\begin{equation}\label{eq:omegahat}
\omega = \widehat{\omega} + \Or(\delta),\qquad
\epsilon \to (0^+,\ldots,0^+).
\end{equation}

\begin{remark}
By applying repeatedly the result on simple singular collapses,
we see that
\[
\lim_{\epsilon_n \to 0^+}
\left( \cdots \lim_{\epsilon_2 \to 0^+}  \Big(\lim_{\epsilon_1 \to 0^+} \omega \Big) \right)
= \widehat{\omega}.
\]
In fact, these repeated limits can be taken in any order.
Nevertheless, this result is weaker than estimate~(\ref{eq:omegahat}),
so we need a formal proof of the estimate.
\end{remark}

We consider the constants
$Q_l = \sqrt{c_{2n+1}-c^*_l}\prod_{i \neq l} |c^*_i - c^*_l|$.
Let $\{P_1(s),\ldots,P_n(s)\}$ be the basis of $\Rset_{n-1}[s]$ univocally
determined by
\[
P_l(c^*_j) = \cases{Q_l & \mbox{if $j=l$}, \\ 0 & \mbox{otherwise.}}
\]
Now, using corollary~\ref{cor:S},
we get that there exists some constants $\zeta_{jl} \in \Rset$ such that
\[
\mathcal{K}_j[\delta_l P_l(s)] = \delta_l \mathcal{K}_j[P_l(s)] =
\cases{1 + \zeta_{ll}\delta_l + \order(\delta) & \mbox{if $j = l$}, \\
       \zeta_{jl}\delta_l + \order(\delta)     & \mbox{otherwise},}
\]
where $0 \le j \le n$ and $1 \le l \le n$.
Therefore, the $n \times n$ linear system
\[
\delta_l \mathcal{K}[P_l(s); \omega] = 0 \mbox{ for $l=1,\ldots,n$}
\]
is $\Or(\delta)$-close to the nonsingular linear system
\[
\left\{
\begin{array}{l}
1- 2\widehat{\omega}_1 = 0 \\
2 (-1)^{l-1} (\widehat{\omega}_{l-1} - \widehat{\omega}_l) = 0
\mbox{ for $l=2,\ldots,n$}
\end{array}
\right.
\]
whose unique solution is $\widehat{\omega} = (1/2,\ldots,1/2)$.
Thus, the asymptotic formula~(\ref{eq:omegahat})
follows from the first item in lemma~\ref{lem:LinearSystems}.

\begin{remark}
The vectorial estimate~(\ref{eq:omegahat}) can be refined in several ways.
For instance,
one can get the componentwise estimates $\omega_1 = 1/2 + \Or(\delta_1)$ and
$\omega_l = \omega_{l-1} + \Or(\delta_l)$ for $l > 1$.
In particular, $\omega_l = 1/2 + \Or(\delta_1,\ldots,\delta_l)$.
Even more, there exists a $n \times n$ constant
\emph{lower triangular} matrix $\mathbf{L}$ such that
\[
\omega = \widehat{\omega} + \mathbf{L} \delta + \order(\delta),
\qquad \epsilon \to (0^+,\ldots,0^+).
\]
We omit the proof,
since we do not need this result and the computations are cumbersome.
\end{remark}

\subsection{Asymptotic behavior of the function $\nu_\rmx$}
\label{sap:Nux}

The function $\nu_\rmx : (0,c) \cup (c,b) \to \Rset$ verifies that
$I(\lambda) + J(\lambda) \rho_\rmx(\lambda) + K(\lambda) \nu_\rmx(\lambda) = 0$,
where the coefficients $I,J,K : (0,c) \cup (c,b) \to \Rset$ were given by
\[
\fl
I(\lambda) = \int_0^{\underline{m}} \frac{\rmd s}{(a-s)\sqrt{T_\rmx(s)}},\quad
J(\lambda) = -2\int_{\overline{m}}^b \frac{\rmd s}{(a-s)\sqrt{T_\rmx(s)}},\quad
K(\lambda) = \frac{2\pi}{\sqrt{-T_\rmx(a)}},
\]
with $T_\rmx(s) = (\lambda-s)(c-s)(b-s)$, $\underline{m} = \min(\lambda,c)$,
and $\overline{m} = \max(\lambda,c)$.
Here, $\rho_\rmx(\lambda) = \rho(\lambda;c,b)$ is the rotation function
of the ellipse obtained by sectioning the ellipsoid $Q$ with the coordinate
plane $\{x=0\}$.
The asymptotic properties of rotation functions of billiards
inside ellipses were established in proposition~\ref{pro:2D}.

First, let us consider the case $\epsilon := \lambda \to 0^+$.
Using lemmas~\ref{lem:0} and~\ref{lem:G} we get:
\begin{itemize}
\item
$I(\epsilon) = I_0 \epsilon^{1/2} + \Or(\epsilon^{3/2})$, where
$I_0 = 2 a^{-1} (bc)^{-1/2}$;
\item
$J(\epsilon) = J_0 + \Or(\epsilon)$, where
$J_0 = -2\int_c^b (a-s)^{-1}(s(s-c)(b-s))^{-1/2}\rmd s$;
\item
$K(\epsilon) = K_0 + \Or(\epsilon)$, where
$K_0 = 2\pi (a(a-c)(a-b))^{-1/2}$;
\item
$\rho_\rmx(\epsilon) = \kappa^G \epsilon^{1/2} + \Or(\epsilon^{3/2})$,
where $\kappa^G = \kappa^G(b,c)$ can be found in proposition~\ref{pro:2D}; and
\item
$\nu_\rmx(\epsilon) =
-(I_0 + J_0 \kappa^G) K_0^{-1} \epsilon^{1/2} + \Or(\epsilon^{3/2}) =
\Or(\epsilon^{1/2})$.
It is possible to check that $(I_0 + J_0 \kappa^G) K_0^{-1} < 0$,
but we do not need it.
\end{itemize}

Next, let us consider the case $\epsilon := b - \lambda \to 0^+$.
We begin by computing the integral
\[
r(\beta,\alpha) := \int_0^\beta \frac{\rmd s}{(\alpha-s)\sqrt{\beta-s}} =
\frac{2}{\sqrt{\alpha - \beta}} \atan \sqrt{\beta/(\alpha-\beta)},
\]
for any $0 < \beta < \alpha$.
Then it is immediate to check that
\[
\int_0^\beta \frac{\rmd s}{(\alpha_+ -s) (\alpha_- -s) \sqrt{\beta -s}} =
\frac{r(\beta,\alpha_-) - r(\beta,\alpha_+)}{\alpha_+ - \alpha_-},
\]
for any $0 < \beta < \alpha_- < \alpha_+$.
We also need the formula
$r(\beta,\alpha) = 2 \pi (\alpha-\beta)^{-1/2}\varrho(\beta,\alpha)$,
where
$\varrho(\beta,\alpha) :=
 \lim_{\gamma \to \alpha^-} \rho(\gamma;\beta,\alpha)$ is one
of the limits of the rotation number described in proposition~\ref{pro:2D}.
Using these formulae, jointly with lemmas~\ref{lem:0} and~\ref{lem:R},
we see that:
\begin{itemize}
\item
$I(b-\epsilon) = I_\ast + \Or(\epsilon)$, where
$I_\ast =
2\pi(a-b)^{-1}\left( (b-c)^{-1/2}\varrho(c,b) -
                     (a-c)^{-1/2} \varrho(c,a) \right)$;
\item
$J(b-\epsilon) = J_\ast + \Or(\epsilon)$, where
$J_\ast = -2\pi (a-b)^{-1}(b-c)^{-1/2}$;
\item
$K(b-\epsilon) = K_\ast + \Or(\epsilon)$, where
$K_\ast = 2\pi (a-b)^{-1}(a-c)^{-1/2}$;
\item
$\rho_\rmx(b-\epsilon) =
\rho_\rmx(b) + \Or(\epsilon) =
\rho(b;c,b) + \Or(\epsilon) =
\varrho(c,b) + \Or(\epsilon)$; and
\item
$\nu_\rmx(b-\epsilon) =
\varrho(c,a) + \Or(\epsilon) =
\rho(a;c,a) + \Or(\epsilon) =
\rho_\rmy(a) + \Or(\epsilon)$.
\end{itemize}

The estimates in the limit $\epsilon := \overline{m} - \underline{m} \to 0^+$,
which means $\lambda \to c$, are:
\begin{itemize}
\item
$I(c \pm \epsilon) =
-(a-c)^{-1}(b-c)^{-1/2} \log \epsilon + \mu + \Or(\epsilon\log\epsilon)$,
where $\mu$ is a constant that can be exactly computed from lemma~\ref{lem:S};
\item
$J(c \pm \epsilon) =
-(a-c)^{-1}(b-c)^{-1/2} \log \epsilon + \eta + \Or(\epsilon\log\epsilon)$,
where $\eta$ is a constant that can be exactly computed from lemma~\ref{lem:S};
\item
$K(c \pm \epsilon) = 2\pi (a-c)^{-1}(a-b)^{-1/2} + \Or(\epsilon)$;
\item
$\rho_\rmx(c \pm \epsilon) =
 1/2 + \kappa^S \log^{-1} \epsilon + \Or(\log^{-2} \epsilon)$,
where $\kappa^S = \kappa^S(c,b) = \acosh (b/c)^{1/2}$ according
to proposition~\ref{pro:2D}; and
\item
$\nu_\rmx(c \pm \epsilon) =
(a-b)^{1/2}\left((a-c)(\eta-\mu) - 2 (b-c)^{-1/2} \kappa^S\right)/2\pi
 + \Or(\log^{-1} \epsilon)$.
After some tedious, but simple, computations,
one gets that
$\nu_\rmx(c) = \varrho(b,a) =\rho(a;b,a) = \rho_\rmz(a)$.
\end{itemize}

\subsection{Asymptotic behavior of the function $\nu_\rmy$}
\label{sap:Nuy}

The function $\nu_\rmy:(b,a) \to \Rset$ verifies that
\(
I(\lambda) + J(\lambda) \rho_\rmy(\lambda) + K(\lambda) \nu_\rmy(\lambda) = 0,
\)
where the coefficients $I,J,K: (b,a) \to \Rset$ were given by
\[
\fl
I(\lambda) = \int_0^c \frac{\rmd s}{(b-s)\sqrt{T_\rmy(s)}},\quad
J(\lambda) = 2\int_\lambda^a \frac{\rmd s}{(s-b)\sqrt{T_\rmy(s)}},\quad
K(\lambda) = -\frac{2\pi}{\sqrt{-T_\rmy(b)}},
\]
with $T_\rmy(s) = (c-s)(\lambda-s)(a-s)$.
Here, $\rho_\rmy(\lambda) = \rho(\lambda;c,a)$ is the rotation function
of the ellipse obtained by sectioning the ellipsoid $Q$ with the coordinate
plane $\{y=0\}$.

We begin with the limit $\epsilon := \lambda - b \to 0^+$.
Using lemmas~\ref{lem:0} and~\ref{lem:S} we see that:
\begin{itemize}
\item
$I(b + \epsilon) = I_0 + \Or(\epsilon)$, where
$I_0 = \int_0^c (b-s)^{-3/2}(c-s)^{-1/2}(a-s)^{-1/2} \rmd s$;
\item
$J(b + \epsilon) = 2\pi(a-b)^{-1/2}(b-c)^{-1/2} \epsilon^{-1/2} + \Or(1)$;
\item
$K(b + \epsilon) = -2\pi (a-b)^{-1/2}(b-c)^{-1/2} \epsilon^{-1/2}$;
\item
$\rho_\rmy(b + \epsilon) = \rho_\rmy(b) + \Or(\epsilon)$; and
\item
$\nu_\rmy(b + \epsilon) = \rho_\rmz(c) + \Or(\epsilon^{1/2})$,
since $\rho_\rmy(b) = \rho(b;c,a) = \rho(c;b,a) = \rho_\rmz(c)$.
\end{itemize}

Next, let us consider the case $\epsilon := a - \lambda \to 0^+$,
which is similar to the limit
$\lim_{\epsilon \to 0^+} \nu_\rmx(b-\epsilon)$ studied in the previous
subsection, so we need the same simple integrals.
Using them, jointly with lemmas~\ref{lem:0} and~\ref{lem:R}, we get:
\begin{itemize}
\item
$I(a - \epsilon) = I_\ast + \Or(\epsilon)$, where
$I_\ast =
2\pi(a-b)^{-1}\left( (b-c)^{-1/2}\varrho(c,b) -
                     (a-c)^{-1/2} \varrho(c,a) \right)$;
\item
$J(a - \epsilon) = J_\ast + \Or(\epsilon)$, where
$J_\ast = 2\pi (a-b)^{-1}(a-c)^{-1/2}$;
\item
$K(a - \epsilon) = K_\ast + \Or(\epsilon)$, where
$K_\ast = -2\pi (a-b)^{-1}(b-c)^{-1/2}$;
\item
$\rho_\rmy(a - \epsilon) = \varrho(c,a) + \Or(\epsilon)$; and
\item
$\nu_\rmy(a-\epsilon) =
\varrho(c,b) + \Or(\epsilon) =
\rho(b;c,b) + \Or(\epsilon) =
\rho_\rmx(b) + \Or(\epsilon)$.
\end{itemize}

\subsection{Asymptotic behavior of the function $\nu_\rmz$}
\label{sap:Nuz}

The function $\nu_\rmz:(0,c) \cup (c,b) \to \Rset$ verifies that
\begin{equation}\label{eq:nuz}
I(\lambda) + J(\lambda) \rho_\rmz(\underline{m}) +
K(\lambda) \nu_\rmz(\lambda) = 0,
\end{equation}
where the coefficients $I,J,K: (0,c) \cup (c,b) \to \Rset$ were given by
\[
\fl
I(\lambda) = \int_0^{\underline{m}} \frac{\rmd s}{(\overline{m}-s)\sqrt{T_\rmz(s)}},\quad
J(\lambda) = 2\int_b^a \frac{\rmd s}{(s-\overline{m})\sqrt{T_\rmz(s)}},\quad
K(\lambda) = -\frac{2\pi}{\sqrt{-T_\rmz(\overline{m})}},
\]
with $T_\rmz(s) = (\underline{m}-s)(b-s)(a-s)$, $\underline{m} = \min(\lambda,c)$,
and $\overline{m} = \max(\lambda,c)$.
Here, $\rho_\rmz(\lambda) = \rho(\lambda;b,a)$ is the rotation function
of the ellipse obtained by sectioning the ellipsoid $Q$ with the coordinate
plane $\{z=0\}$.

First, let us consider the case $\epsilon := \lambda \to 0^+$.
Using lemmas~\ref{lem:0} and~\ref{lem:G} we see that:
\begin{itemize}
\item
$I(\epsilon) = I_0 \epsilon^{1/2} + \Or(\epsilon^{3/2})$, where
$I_0 = 2 c^{-1} (ab)^{-1/2}$;
\item
$J(\epsilon) = J_0 + \Or(\epsilon)$, where
$J_0 = 2\int_b^a (s-c)^{-1}(s(s-b)(a-s))^{-1/2}\rmd s$;
\item
$K(\epsilon) = K_0 + \Or(\epsilon)$, where
$K_0 = -2\pi (c(b-c)(a-c))^{-1/2}$;
\item
$\rho_\rmz(\underline{m}) = \rho_\rmz(\min(\epsilon,c)) =
\rho_\rmz(\epsilon) = \kappa^G \epsilon^{1/2} + \Or(\epsilon^{3/2})$,
where the constant $\kappa^G = \kappa^G(b,a)$ can be found in
proposition~\ref{pro:2D}; and
\item
$\nu_\rmz(\epsilon) =
-(I_0 + J_0 \kappa^G) K_0^{-1} \epsilon^{1/2} + \Or(\epsilon^{3/2}) =
\Or(\epsilon^{1/2})$, with $(I_0 + J_0 \kappa^G) K_0^{-1} < 0$.
\end{itemize}

The estimates in the limit $\epsilon := \overline{m} - \underline{m} \to 0^+$,
which means $\lambda \to c$, are:
\begin{itemize}
\item
$I(c \pm \epsilon) =
\pi (a-c)^{-1/2}(b-c)^{-1/2} \epsilon^{-1/2} + \Or(1)$; see lemma~\ref{lem:S};
\item
$J(c \pm \epsilon) = \Or(1)$;
\item
$K(c \pm \epsilon) =
-2\pi (a-c)^{-1/2}(b-c)^{-1/2} \epsilon^{-1/2} + \Or(\epsilon^{1/2})$;
\item
$\rho_\rmz(\underline{m}) = \rho_\rmz(\min(c \pm \epsilon,c)) = \rho_\rmz(c) + \Or(\epsilon)$,
since $\rho_\rmz(\lambda)$ is analytic at $\lambda=c$; and
\item
$\nu_\rmz(c \pm \epsilon) = 1/2  + \Or(\epsilon^{1/2})$.
\end{itemize}

Next, we consider the case $\epsilon := b - \lambda \to 0^+$.
Using lemmas~\ref{lem:0} and~\ref{lem:R} we get:
\begin{itemize}
\item
$I(b - \epsilon) = \Or(1)$;
\item
$J(b - \epsilon) = J_\ast \epsilon^{-1/2} + \Or(1)$, where
$J_\ast = 2\pi (a-b)^{-1/2}(b-c)^{-1/2}$;
\item
$K(b - \epsilon) = K_\ast \epsilon^{-1/2}+ \Or(\epsilon^{1/2})$, where
$K_\ast = -2\pi (a-b)^{-1/2}(b-c)^{-1/2}$;
\item
$\rho_\rmz(\underline{m}) = \rho_\rmz(\min(b-\epsilon,c)) = \rho_\rmz(c)$; and
\item
$\nu_\rmz(b-\epsilon) = \rho_\rmz(c) + \Or(\epsilon^{1/2})$.
\end{itemize}

\section{A topological lemma}
\label{ap:Topological}

We recall that the complement of any Jordan curve $X$ in the plane $\Rset^2$
has two distinct connected components.
One of them is bounded and simply connected (the interior,
denoted by $\mathcal{B}_X$) and the other is unbounded
(the exterior, denoted by $\mathcal{U}_X$).

\begin{lem}\label{lem:LocalGlobal}
Let $X$ and $Y$ be two Jordan curves of $\Rset^2$.
If $f: \mathcal{B}_X \to \Rset^2$ is a bounded local homeomorphism
that has a continuous extension to the boundary $X$
such that $f(X) \subset Y$,
then $f: \mathcal{B}_X \to \mathcal{B}_Y$ is a global homeomorphism.
\end{lem}

\proof
We note that
$W = f(\mathcal{B}_X)$ is a non-empty open bounded subset of $\Rset^2$
such that
\[
\partial W = \partial f(\mathcal{B}_X) \subset
f(\partial \mathcal{B}_X) = f(X) \subset Y.
\]
Next, we are going to prove that $W = \mathcal{B}_Y$.
Using that $\partial W \subset Y$, we deduce that
the intersection $W \cap \mathcal{B}_Y$ (respectively, $W \cap \mathcal{U}_Y$)
is open and closed in $\mathcal{B}_Y$ (respectively, in $\mathcal{U}_Y$),
so it is either the empty set or the whole interior (respectively, exterior).
Therefore, we deduce that:
1) $W \cap \mathcal{U}_Y = \emptyset$, because $W$ is bounded;
2) $W \cap Y= \emptyset$, because $W$ is open; and
3) $W \cap \mathcal{B}_Y = \mathcal{B}_Y$, because $W$ is open and non-empty.
That is, $f(\mathcal{B}_X) = W = \mathcal{B}_Y$.

Once we know that $f:\mathcal{B}_X \to \mathcal{B}_Y$ is a surjective
local homeomorphism, we deduce from covering space theory that it is
a global homeomorphism.
It suffices to realize that $\mathcal{B}_X$ is connected and open,
and $\mathcal{B}_Y$ is simply connected.
\qed

In particular, if $f:\mathcal{B}_X \to \Rset^2$ is smooth or analytic,
then its inverse is also smooth or analytic.
This means that if $f$ is a local diffeomorphism whose image is bounded and
that has a continuous extension to the boundary $X$ such that $f(X) \subset Y$,
then $f: \mathcal{B}_X \to \mathcal{B}_Y$ is a global diffeomorphism.

\end{document}